\documentclass[12pt]{article}
\usepackage{amsmath,amsthm,amsfonts,amssymb}

\usepackage[usenames]{color} 

\newtheorem{thm}{Theorem}[section]

\newtheorem{lemma}[thm]{Lemma}

\newtheorem{cor}[thm]{Corollary}
\newtheorem{corollary}[thm]{Corollary}
\newtheorem{proposition}[thm]{Proposition}

\theoremstyle{definition}

\newtheorem{definition}[thm]{Definition} 

\newtheorem{notn}[thm]{Notation}

\newtheorem{ex}[thm]{Example}

\newtheorem{remark}[thm]{Remark}
\theoremstyle{remark}

\newcounter{remarks}

\newcounter{enumitemp}
\newenvironment{enumeratecontinue}{
 \setcounter{enumitemp}{\value{enumi}}
 \begin{enumerate}
 \setcounter{enumi}{\value{enumitemp}}
}
{
 \end{enumerate}
}

\DeclareMathOperator{\Fix}{Fix}
\DeclareMathOperator{\Per}{Per}
\DeclareMathOperator{\PA}{P}

\DeclareMathOperator{\PF}{PF}

\DeclareMathOperator{\Int}{int}
\DeclareMathOperator{\EF}{EF}

\newcommand{\Z}{{\mathbb Z}}

\newcommand{\T}{{\mathbb T}}
\def\H{{\mathbb H}}

\newcommand{\f}{F_n}

\newcommand{\E}{{\mathcal E}}

\newcommand{\oone}{\phi} 
\newcommand{\otwo}{\psi}

\newcommand{\aone}{\Phi}
\newcommand{\atwo}{\Psi}

\newcommand{\Out}{\mathsf{Out}}
\newcommand{\Aut}{\mathsf{Aut}}
\newcommand{\Inn}{\mathsf{Inn}}
\newcommand{\mcg}{\mathsf{MCG}}
\newcommand{\Stab}{\mathsf{Stab}}
\newcommand{\ffs}{free factor system}

\newcommand{\EG}{EG}

\newcommand{\F}{\mathcal F}
\newcommand{\rtt}{relative train track map}

\def\L{\mathcal L}

\newcommand{\fG} {f : G \to G}
\newcommand{\ti} {\tilde}
\newcommand{\iNp} {indivisible Nielsen path}

\newcommand{\filt}{\emptyset = G_0 \subset G_1 \subset \ldots  \subset G_N = G}

\newcommand{\irt}{improved relative train track}
\newcommand{\ofn}{\Out(F_n)}
\newcommand{\eg}{EG}
\newcommand{\noneg}{NEG}

\newcommand{\cs}{completely split}

\newcommand{\wt}{\widetilde}

\newcommand{\ipNp}{indivisible periodic Nielsen path}

\newcommand{\ct}{CT}
\newcommand{\csirt}{CT}

\newcommand{\comment}[1]{}

\title{The Recognition Theorem for $\Out(F_n)$} 
\author{Mark Feighn\thanks{This material is based upon work supported by the
National Science Foundation under Grant No.~DMS-0805440.}\and
Michael Handel\thanks{This material is based upon work supported by
the National Science Foundation under Grant No.~DMS-0706719 .}}
 
\begin{document}

 \maketitle

\begin{abstract} 
Our goal is to find dynamic invariants that completely determine
elements of the outer automorphism group $\Out(F_n)$ of the free group
$F_n$ of rank $n$. To avoid finite order phenomena, we do this for
{\it forward rotationless} elements. This is not a serious  
restriction. For example, there is $K_n>0$
depending only on $n$ such that, for all $\phi\in\Out(F_n)$,
$\phi^{K_n}$ is forward rotationless. An important part of our
analysis is to show that rotationless elements are represented by
particularly nice relative train track maps.
 \end{abstract}

\tableofcontents

\section{Introduction} \label{s:intro}  

 The Thurston classification theorem for  mapping class groups of  
surfaces inspired a surge in research on the outer automorphism group $\Out(F_n)$ of the  free group  of rank $n$.   One direction of this research is the development and use of relative train track maps which are the analog of the normal forms for mapping classes.   Thurston's normal forms give rise  to invariants that completely determine a mapping class, perhaps after passing to a finite power.  In this paper we provide similar invariants for elements of  $\Out(F_n)$ and we   add an important feature to relative train track maps.  

   We begin by recalling the invariants associated to a mapping class.   After passing to a  finite power,  there is a decomposition of the surface into invariant subsurfaces on which the action of the mapping class is   either  a Dehn twist  
in an annulus, trivial or pseudo-Anosov; in the pseudo-Anosov case, 
 each singular ray of the associated measured
foliations is invariant.  The mapping class is  
completely determined by the (isotopy classes of the) core curves of the annuli, the Dehn twist degrees,  
the  pseudo-Anosov measured foliations and the expansion factors on the pseudo-Anosov measured foliations.

To see how this might generalize to $\Out(F_n)$, consider the special case in which the mapping class $\psi$ is   a composition of disjoint Dehn twists and so is determined by the twisting circles and the degrees of twist.    The dual point of view is useful here.  The complementary components of the twisting curves are maximal subsurfaces on which $\psi$ acts trivially.  Viewing the mapping class group as the group of outer automorphism of the fundamental group $\pi_1(S)$ of the surface $S$, each complementary component   determines,  up to   conjugacy, a    subgroup of $\pi_1(S)$ of rank at least two that is the fixed subgroup $\Fix(\Psi)$ of  some automorphism $\Psi$ of $\pi_1(S)$ that represents the outer automorphism  $\psi$.   If $\Psi_1$ and $\Psi_2$ are two distinct such automorphisms,  corresponding to either the same or distinct subsurfaces, then   $\Fix(\Psi_1) \cap \Fix(\Psi_2)$ is either  trivial or a maximal cyclic group $ \langle a \rangle$.  In the latter case   $\Psi_2 = i_a^d \Psi_1$ for some $ d \ne 0$ where $i_a$ is the inner automorphism determined by $a$; moreover, the conjugacy class determined by $a$ represents a twisting curve for $\mu$ with twisting degree $\pm d$.  

This point of view, focusing on fixed subgroups of rank at least two and their intersections, is sufficient  \cite{cl:UL} if one restricts  to elements of   $\Out(F_n)$  that have linear growth.\footnote{$\phi\in\Out(F_n)$ has {\em
    linear growth} if, for all conjugacy classes of elements $a\in
  F_n$, the cyclically reduced word length of the conjugacy class of
  $\phi^m(a)$ is bounded by a linear function (depending on $a$) of
  $m$.}  
For general $\phi \in \Out(F_n)$,   we must also account for exponential and non-linear polynomial growth.  To do  
this, we consider the boundary of the free group.  An automorphism $ 
\Phi$ induces a homeomorphism $\hat\Phi\colon\partial F_n\to\partial F_n  
$.  In the linear case, the fixed point set $\Fix(\hat \Phi)$ of $\hat\Phi$ is equal to  
the boundary of the fixed subgroup $\Fix(\Phi)$ and contains no  
isolated points. In the general case, however, there are isolated  
points,
and these contain essential information about the automorphism.  For example,  the set of  attracting laminations  associated to $\phi$ can be recovered from the isolated attracting points.

The idea then is to replace $\Fix(\Phi)$ with the set $\Fix_N(\hat
\Phi)$ of non-isolated points and attractors in $\Fix(\hat \Phi)$; see
\cite{gjll}, \cite{ll:ends}, \cite{bfh:tits3} and \cite{ll:asymptotic}
where this same idea has been used effectively.  In
section~\ref{principal def}, we define the set $\PA(\oone) \subset \Aut(F_n)$ of {\em
principal automorphisms} representing $\phi$.    In the  case of linear growth, $\Phi$ is principal if and only if $\Fix(\Phi)$ has rank at least two.    The invariants that determine $\phi$, after possibly passing to a finite power,   are the sets $\Fix_N(\hat \Phi)$ as $\Phi$ varies
over $\PA(\oone)$,  the expansion factors for the attracting laminations of $\phi$ 
  and
twisting coordinates associated to  pairs of principal automorphisms whose fixed points sets intersect non-trivially.   

It is common when studying elements   $\phi \in \Out(F_n)$ to \lq stabilize\rq\ $\phi$ by  replacing  it with a power  $\phi^k$.       In section~\ref{rotOutAut}, we specify a subset of $\Out(F_n)$  whose elements require no stabilization.  These outer automorphisms are said to be {\em forward rotationless}.  In Lemma~\ref{uniform bound} we prove that there is $K_n>0$,
depending only on $n$,  so that   $\phi^{K_n}$ is forward rotationless for all  $\oone \in
\Out(F_n)$.   We also define what it means for a \rtt\ $\fG$ to be {\em rotationless} and
prove (Proposition~\ref{rotationless}) that $\oone$ is forward
rotationless if and only if some (every) \rtt\ representing it is
rotationless.  (There is no need to add \lq forward\rq\ to this
terminology because $f^k$ is only defined for $k\ge 1$.)   It is easy to check if $\fG$ is rotationless and if not to find the minimal $k$ such that $f^k$ is rotationless.

We can now state our main result.  Complete details and further motivation are supplied in  
section~\ref{recog}.   In addition to being of intrinsic interest this theorem is needed   in  \cite{fh:abeliansubgroups}.
  The set of attracting laminations for $\oone$ is denoted   $\L(\oone)$ and  the expansion factor for $\oone$ on $\Lambda \in \L(\oone)$ is denoted $\PF_{\Lambda}(\oone)$.  
   
\vspace{.1in}
\noindent {\em {\bf Theorem~\ref{t:recognition} (Recognition Theorem)}
Suppose that $\oone,\otwo\in\Out(\f)$ are forward rotationless and
that
\begin{enumerate}
\item $\PF_{\Lambda}(\oone)=\PF_{\Lambda}(\otwo)$, for all $\Lambda\in
\L(\oone)=\L(\otwo)$; and
\item there is bijection $B : \PA(\oone) \to \PA(\otwo)$ such that:  
\begin{itemize} 
\item [(i)] {\bf (fixed sets preserved)} $\Fix_N(\hat \aone) =
\Fix_N(\widehat{B(\aone)})$; and
\item [(ii)] {\bf (twist coordinates preserved)} if $ w \in
\Fix(\aone)$ and $\aone, i_w\aone \in \PA(\oone)$, then $B(i_w\aone)=
i_wB(\aone)$.
\end{itemize}
\end{enumerate}
Then  $\oone=\otwo$.}
 
\vspace{.1in}

In the course of proving Theorem~\ref{t:recognition}, we construct \rtt s   that are better than the   those constructed in \cite{bfh:tits1}.  We also reorganize elements of the theory to make future modification and referencing of results easier.   

	The idea behind
relative train track maps is that one can study the action of an outer automorphism   $\oone$ on conjugacy classes in $F_n$ or on $\partial
F_n$ by analyzing the action of a homotopy equivalence $\fG$ of a marked graph $G$
representing $\oone$ on paths, circuits and lines in $G$.  For
simplicity, suppose that $\sigma$ is a finite path in $G$, which by
convention is always assumed to be the immersed image of a compact interval.
The image $f(\sigma)$ of $\sigma$ is homotopic rel endpoints to a path
denoted $f_\#(\sigma)$.  Replacing $f(\sigma)$ with $f_\#(\sigma)$ is
called {\em tightening} and is analogous to replacing a word in a set
of generators for $F_n$ with a reduced word in those generators.  
A decomposition  of
$\sigma$ into subpaths $\sigma = \sigma_1\ldots \sigma_m$ is a {\em splitting} if $f^k_\#(\sigma) =
f^k_\#(\sigma_1)\ldots f^k_\#(\sigma_m) $ for all $k \ge
0$; i.e. if one can tighten the image of $\sigma$ under any
iterate of $f$ by tightening the images of the $\sigma_i$'s.     The more one can split $\sigma$ and the better one can understand the subpaths $\sigma_i$, the more effectively one can analyze the iterates $f^k_\#(\sigma)$.

   Relative train track maps were defined
and constructed in \cite{bh:tracks} with   exponentially growing
strata in mind.  Few restrictions were placed on the non-exponentially
growing strata.  This was rectified in \cite{bfh:tits1} where {\em
improved relative train tracks (IRTs)} are defined and shown to exist
for a sufficiently high, but unspecified, iterate of $\oone$.    For our current application, 
    IRTs  are inadequate.   In Theorem~\ref{comp sp exists}, we prove that every forward rotationless $\phi$ is represented by a  
\rtt\   $\fG$ that has all the essential properties of an IRT (see section~\ref{s:comparison}) and has the additional feature that, for
all $\sigma$ and all sufficiently large $k$, there is a canonical
splitting (called the {\em complete splitting}) of $f^k_\#(\sigma)$ into simple, explicitly described
subpaths.  (Splittings in an IRT  are  not  canonical and the subpaths $\sigma_i$ are understood
more inductively than explicitly.)  Such  $\fG$, called {\em \ct s}, are used
in the proof of the Recognition Theorem and in the classification of
abelian subgroups given in \cite{fh:abeliansubgroups}.  It is likely
that the existence of complete splittings  will be useful in other contexts as well.  For example,
complete  splittings are {\em hard splittings} as defined in
\cite{bg:hard}.

\section{Preliminaries}   \label{preliminaries}  

Fix $n \ge 2$ and let $F_n$ be the free group of rank $n$.  Denote the
automorphism group of $F_n$ by $\Aut(F_n)$, the group of inner
automorphisms of $F_n$ by $\Inn(F_n)$ and the group of outer
automorphisms of $F_n$ by $\Out(F_n) = \Aut(F_n)/\Inn(F_n)$.  We
follow the convention that elements of $\Aut(F_n)$ are denoted by
upper case Greek letters and that the same Greek letter in lower case
denotes the corresponding element of $\Out(F_n)$.  Thus $\aone \in
\Aut(F_n)$ represents $\oone \in \Out(F_n)$.

\subsection{Marked Graphs and Outer Automorphisms}    Identify $F_n$ with 
$\pi_1(R_n,*)$ where $R_n$ is the rose with one vertex $*$ and $n$
edges.  A {\em marked graph} $G$ is a graph of rank $n$, all of whose
vertices have valence at least two, equipped with a homotopy
equivalence $m : R_n \to G$ called a {\em marking}.  Letting $b = m(*)
\in G$, the marking determines an identification of $F_n$ with
$\pi_1(G,b)$.  It is often assumed that $G$ does not have valence two
vertices, but such vertices occur naturally in relative train track
theory so we allow them.

A homotopy equivalence $\fG$ and a path $\sigma$ from $b$ to $f(b)$
determines an automorphism of $\pi_1(G,b)$ and hence an element of
$\Aut(F_n)$.  If $f$ fixes $b$ and no path is specified, then we use
the trivial path.  This construction depends only on the homotopy
class of $\sigma$ and, as the homotopy class varies, the automorphism
ranges over all representatives of the associated outer automorphism
$\oone$.  We say that {\em$\fG$ represents $\phi$}.  We always assume
that $f$ maps vertices to vertices and that the restriction of $f$ to any edge is an immersion.

\subsection{Paths, Circuits and Edge Paths}

  Let $\Gamma$ be the universal cover of a marked graph $G$ and let
$pr : \Gamma\to G$ be the covering projection.  A proper map $\ti
\sigma : J \to \Gamma$ with domain a (possibly infinite) closed interval $J$
will be called a {\it path in $\Gamma$} if it is an embedding or if
$J$ is finite and the image is a single point; in the latter case we
say that $\ti \sigma$ is {\it a trivial path}.  If $J$ is finite, then
$\ti \sigma : J \to \Gamma$ is homotopic rel endpoints to a unique
(possibly trivial) path $[\ti \sigma]$; we say that {\it $[\ti
\sigma]$ is obtained from $\ti \sigma$ by tightening}. If $\ti f
:\Gamma \to \Gamma$ is a lift of a homotopy equivalence $\fG$, we
denote $ [\ti f(\ti \sigma)]$ by $\ti f_\#(\ti \sigma)$.

We will not distinguish between paths in $\Gamma$ that differ only by
an orientation preserving change of parametrization. Thus we are
interested in the oriented image of $\ti \sigma$ and not $\ti \sigma$
itself.  If the domain of $\ti \sigma$ is finite, then the image of
$\ti \sigma$ has a natural decomposition as a concatenation $\wt E_1
\wt E_2 \ldots \wt E_{k-1} \wt E_k$ where $\wt E_i$, $1 < i < k$, is
an edge of $\Gamma$, $\wt E_1$ is the terminal segment of an edge and
$\ti E_k$ is the initial segment of an edge. If the endpoints of the
image of $\ti \sigma$ are vertices, then $\wt E_1$ and $\wt E_k $ are
full edges.  The sequence $\wt E_1 \wt E_2\ldots \wt E_k$ is called
{\it the edge path associated to $\ti \sigma$.}  This notation extends
naturally to the case that the interval of domain is half-infinite or
bi-infinite.  In the former case, an edge path has the form $\wt
E_1\wt E_2\ldots$ or $\ldots \wt E_{-2} \wt E_{-1}$ and in the latter
case has the form $\ldots\wt E_{-1}\wt E_0\wt E_1\wt E_2\ldots$.

A {\it path in $G$} is the composition of the projection map $pr$ with
a path in $\Gamma$.  Thus a map $\sigma : J \to G$ with domain a
(possibly infinite) closed interval will be called a path if it is an
immersion or if $J$ is finite and the image is a single point; paths
of the latter type are said to be {\em trivial}. If $J$ is finite,
then $ \sigma : J \to G$ is homotopic rel endpoints to a unique
(possibly trivial) path $[ \sigma]$; we say that {\it $[ \sigma]$ is
obtained from $ \sigma$ by tightening}. For any lift $\ti \sigma : J
\to \Gamma$ of $\sigma$, $[\sigma] = pr[\ti \sigma]$.  We denote
$[f(\sigma)]$ by $f_\#(\sigma)$.  We do not distinguish between paths
in $G$ that differ by an orientation preserving change of
parametrization.  The {\it edge path associated to $\sigma$} is the
projected image of the edge path associated to a lift $\ti \sigma$.
Thus the edge path associated to a path with finite domain has the
form $ E_1 E_2 \ldots E_{k-1} E_k$ where $ E_i$, $1 < i < k$, is an
edge of $G$, $ E_1$ is the terminal segment of an edge and $ E_k$ is
the initial segment of an edge.  We will identify paths with their
associated edge paths whenever it is convenient.

We reserve the word {\it circuit} for an immersion $\sigma : S^1 \to
G$.  Any homotopically non-trivial map $\sigma : S^1 \to G$ is
homotopic to a unique circuit $[\sigma]$. As was the case with paths,
we do not distinguish between circuits that differ only by an
orientation preserving change in parametrization and we identify a
circuit $\sigma$ with a {\it cyclically ordered edge path}
$E_1E_2\dots E_k$.

A path or circuit \emph{crosses} or \emph{contains} an edge if that
edge occurs in the associated edge path.  For any path $\sigma$ in $G$
define $\bar \sigma$ to be \lq $\sigma$ with its orientation
reversed\rq.  For notational simplicity, we sometimes refer to the
inverse of $\ti \sigma$ by $\ti \sigma^{-1}$.

A decomposition of a path or circuit into subpaths is a {\em splitting} for 
$\fG$ and is denoted $\sigma = \ldots \cdot \sigma_1\cdot\sigma_2\cdot \ldots  
$ if $f^k_\#(\sigma) = \ldots  f^k_\#(\sigma_1) f^k_\#(\sigma_2)\ldots  $ 
for all $k \ge 0$.  In other words, a decomposition of $\sigma$ into subpaths 
$\sigma_i$ is a splitting if one can tighten the image of $\sigma$
under any iterate of $f_\#$ by tightening the images of the $\sigma_i$'s.

A path $\sigma$ is a {\em periodic Nielsen path} if $\sigma$ is
non-trivial and $f^k_\#(\sigma) = \sigma$ for some $k \ge 1$.  The minimal such $k$ is the {\em period} of $\sigma$ and if
the period is one then $\sigma$ is a {\em Nielsen path}.  A (periodic)
Nielsen path is {\em indivisible} if it does not decompose as a
concatenation of non-trivial (periodic) Nielsen subpaths. A path  or circuit  is
{\em root-free} if it is not equal to $\mu^k$ for some simpler path $\mu$ and some $k > 1$.

\subsection{Automorphisms and Lifts}  \label{autos and lifts} Section 1 of 
\cite{gjll} and section 2.1 of \cite{bfh:tits3} are good sources for
facts that we record below without specific references.  The universal
cover $\Gamma$ of a marked graph $G$ with marking $m :R_n \to G$ is a
simplicial tree.  We always assume that a base point $\ti b \in
\Gamma$ projecting to $b = m(*) \in G$ has been chosen, thereby
defining an action of $F_n$ on $\Gamma$.  The set of ends $\E(\Gamma)$
of $\Gamma$ is naturally identified with the boundary $\partial F_n$
of $F_n$ and we make implicit use of this identification throughout
the paper.

Each non-trivial $c \in F_n$ acts by a {\em covering translation} $T_c : \Gamma
\to \Gamma$ and each $T_c$ induces a homeomorphism $\hat T_c :
\partial F_n \to \partial F_n$ that fixes two points, a sink $T_c^+$
and a source $T_c^-$.  The line in $\Gamma$ whose ends converge to
$T_c^-$ and $T_c^+$ is called the {\em axis of $T_c$} and is denoted
$A_c$.  The image of $A_c$ in $G$ is the circuit corresponding to the
conjugacy class of $c$.

If $\fG$ represents $\phi \in \Out(F_n)$ then a path $\sigma$ from $b$
to $f(b)$ determines both an automorphism representing $\phi$ and a
lift of $f$ to $\Gamma$.  This defines a bijection between the set of
lifts $\ti f : \Gamma \to \Gamma$ of $\fG$ and the set of
automorphisms $\Phi : F_n \to F_n$ representing $\phi$.  Equivalently,
this bijection is defined by $\ti f T_c = T_{\Phi(c)} \ti f$ for all
$c \in F_n$.  We say that $\ti f$ {\em corresponds to} $\aone$ or {\em
is determined by} $\aone$ and vice versa.  Under the identification of
$\E (\Gamma)$ with $\partial F_n$, a lift $\ti f$ determines a
homeomorphism $\hat f$ of $\partial F_n$.  An automorphism $\Phi$ also
determines a homeomorphism $\hat \Phi$ of $\partial F_n$ and $\hat f =
\hat \aone$ if and only if $\ti f$ corresponds to $\aone$.  In
particular, $\hat i_c = \hat T_c$ for all $c \in F_n$ where $i_c(w) =
cwc^{-1}$ is the inner automorphism of $F_n$ determined by $c$.  We
use the notation $\hat f$ and $\hat \Phi$ interchangeably depending on
the context.

We are particularly interested in the dynamics of $\hat f  = \hat \aone$. The 
following two lemmas are contained in Lemmas 2.3 and 2.4 of 
\cite{bfh:tits3}  and in Proposition 1.1 of \cite{gjll}. 

\begin{lemma} \label{l: first from bk3} Assume that   $\ti f : \Gamma \to 
\Gamma$ corresponds to $\Phi \in \Aut(F_n)$.  Then the following are
equivalent:
\begin{description}
\item[(i)] $c \in \Fix(\Phi)$.
\item[(ii)] $T_c$ commutes with $\ti f$.
\item[(iii)] $\hat T_c$ commutes with $\hat f$.
\item[(iv)]  $\Fix(\hat T_c) \subset \Fix(\hat f) =  \Fix(\hat \Phi)$.
\item[(v)] $\Fix(\hat f) =  \Fix(\hat \Phi)$ is $\hat T_c$-invariant.
\end{description}
\end{lemma}

\begin{remark}
It is not hard to see that $T_c^+\in \Fix(\Phi)$ if and only if
$T_c^-\in\Fix(\Phi)$.
\end{remark}
 
A point $P \in \partial F_n$ is an {\em attractor} for $\hat \Phi$ if
it has a neighborhood $U$ such that $\hat \Phi(U) \subset U$ and such
that $\cap_{n=1}^\infty \hat \Phi^n(U) = P$.  If $Q$ is an attractor
for $\hat \Phi^{-1}$ then we say that it is a {\em repeller} for $\hat
\Phi$.

\begin{lemma} \label{l: second from bk3}  Assume that   $\ti f : \Gamma \to 
\Gamma$ corresponds to $\Phi \in \Aut(F_n)$ and that $\Fix(\hat \Phi)
\subset \partial F_n$ contains at least three points.  Denote
$\Fix(\Phi)$ by $\mathbb F$ and the corresponding subgroup of covering
translations of $\Gamma$ by $\T(\Phi)$.  Then
\begin{description}
\item[(i)] $\partial \mathbb F$ is naturally identified with the
closure of $\{T_c^{\pm} : T_c \in \T(\Phi)\}$ in $\partial F_n$.  None
of these points is isolated in $\Fix(\hat \Phi)$.
\item[(ii)] Each point in $\Fix(\hat \Phi) \setminus \partial \mathbb
F$ is isolated and is either an attractor or a repeller for the action
of $\hat \Phi$.
\item[(iii)] There are only finitely many $\T(\Phi)$-orbits in
$\Fix(\hat \Phi) \setminus \partial \mathbb F$.
\end{description}
\end{lemma} 

\subsection{Lines and Laminations}\label{laminations}
Suppose that $\Gamma$ is the universal cover of a marked graph $G$.
An unoriented bi-infinite  path in $\Gamma$ is called
a {\em line in $\Gamma$}.  The {\em space of lines in $\Gamma$} is
denoted $\ti {\cal B} (\Gamma)$ and is equipped with what amounts to
the compact-open topology.  Namely, for any finite path $\ti \alpha_0
\subset \Gamma$ (with endpoints at vertices if desired), define $N(\ti
\alpha_0) \subset \ti {\cal B}(\Gamma)$ to be the set of lines in
$\Gamma$ that contain $\ti \alpha_0$ as a subpath.  The sets $N(\ti
\alpha_0)$ define a basis for the topology on $\ti {\cal B}(\Gamma)$.

An unoriented bi-infinite   path in $G$ is called a {\em line
in $G$}.  {\em The space of lines in $G$} is denoted ${\cal B}(
G)$. There is a natural projection map from $\ti {\cal B} (\Gamma)$ to
${\cal B}(G)$ and we equip ${\cal B}( G)$ with the quotient topology.

A line in $\Gamma$ is determined by the unordered pair of its
endpoints $(P,Q)$ and so corresponds to a point in the {\em space of
abstract lines} defined to be $\ti {\cal B} := ((\partial F_n \times
\partial F_n) \setminus \Delta)/ \Z_2$, where $\Delta$ is the diagonal
and where $\Z_2$ acts on $\partial F_n \times \partial F_n$ by
interchanging the factors.  The action of $F_n$ on $\partial F_n$
induces an action of $F_n$ on $\ti {\cal B}$ whose quotient space is
denoted ${\cal B}$.  The \lq endpoint map\rq\ defines a homeomorphism
between $\ti {\cal B}$ and $\ti {\cal B}(\Gamma)$ and we use this
implicitly to identify $\ti {\cal B}$ with $\ti {\cal B}(\Gamma)$ and
hence $\ti {\cal B}(\Gamma)$ with $\ti {\cal B}(\Gamma')$ where
$\Gamma'$ is the universal cover of any other marked graph $G'$.
There is a similar identification of ${\cal B}(G)$ with ${\cal B}$ and
with ${\cal B}(G')$.  We sometimes say that the line in $G$ or
$\Gamma$ corresponding to an abstract line is {\em the realization} of
that abstract line in $G$ or $\Gamma$.

A closed set of lines in $G$ or a  closed $F_n$-invariant set of lines in
$\Gamma$ is called a {\em lamination} and the lines that compose it
are called {\em leaves}.  If $\Lambda$ is a lamination in $G$ then we
denote its pre-image in $\Gamma$ by $\ti \Lambda$ and vice-versa.
 
Suppose that $\fG$ represents $\oone$ and that $\ti f$ is a lift of
$f$.  If $\ti \gamma$ is a line in $\Gamma$ with endpoints $P$ and
$Q$, then there is a bounded homotopy from $\ti f(\ti \gamma)$ to the
line $\ti f_\#(\gamma)$ with endpoints $\hat f(P)$ and $\hat f(Q)$.
This defines an action $\ti f_\#$ of $\ti f$ on lines in $\Gamma$.  If
$\aone \in \Aut(F_n)$ corresponds to $\ti f$ then $\Phi_\# = \ti f_\#$
is described on abstract lines by $(P,Q) \mapsto (\hat \aone(P), \hat
\aone(Q))$.  There is an induced action $\oone_\#$ of $\oone$ on lines
in $G$ and in particular on laminations in $G$.

To each $\oone \in \Out(F_n)$ is associated a finite $\phi$-invariant
set of laminations $\L(\oone)$ called the set of {\em attracting
laminations} for $\oone$.  The individual laminations need not be
$\oone$-invariant.  By definition (see Definition 3.1.5 of
\cite{bfh:tits1}) $\L(\oone) = \L(\oone^k)$ for all $k \ge 1$ and each
$\Lambda \in \L(\oone)$ contains birecurrent leaves, called {\em
generic leaves}, whose weak closure is all of $\Lambda$.  Complete details
on $\L(\oone)$ can be found in section 3 of \cite{bfh:tits1}.

A point $P \in \partial F_n$ determines a lamination $\Lambda(P)$,
called {\em the accumulation set of $P$}, as follows.  Let $\Gamma$ be
the universal cover of a marked graph $G$ and let $\ti R$ be any ray
in $\Gamma$ converging to $P$.  A line $\ti \sigma \subset \Gamma$
belongs to $\widetilde{ \Lambda(P)}$ if every finite subpath of $\ti
\sigma$ is contained in some translate of $\ti R$.  Since any two rays
converging to $P$ have a common infinite end, this definition is
independent of the choice of $\ti R$.  The bounded cancellation lemma
implies (cf. Lemma 3.1.4 of \cite{bfh:tits1}) that this definition is
independent of the choice of $G$ and $\Gamma$ and that $\hat
\aone_\#(\widetilde{\Lambda(P)}) = \widetilde{\Lambda(\hat
\aone(P))}$.  In particular, if $P \in \Fix(\hat \aone)$ then
$\Lambda(P)$ is $\oone_\#$-invariant.

\subsection{Free Factor Systems} 
The conjugacy class of a free factor $F^i$ of $F_n$ is denoted
$[[F^i]]$.  If $F^1,\ldots,F^k$ are non-trivial free factors and if
$F^1 \ast \ldots \ast F^k$ is a free factor then we say that the
collection $\{[[F^1]],\ldots,[[F^k]]\}$ is a {\em free factor system}.
For example, if $G$ is a marked graph and $G_r \subset G$ is a
subgraph with non-contractible components $C_1,\ldots, C_k$ then the
conjugacy class $[[\pi_1(C_i)]]$ of the fundamental group of $C_i$ is
well defined and the collection of these conjugacy classes is a free
factor system denoted $\F(G_r)$; we say that $G_r$ {\em realizes}
$\F(G_r)$.

The image of a free factor $F$ under an element of $\Aut(F_n)$ is a
free factor.  This induces an action of $\Out(F_n)$ on the set of free
free systems.  We sometimes say that a free factor is
$\oone$-invariant when we really mean that its conjugacy class is
$\oone$-invariant.  If $[[F]]$ is $\oone$-invariant then $F$ is
$\aone$-invariant for some automorphism $\aone$ representing $\oone$
and $\aone|F$ determines a well defined element $\oone|F$ of
$\Out(F)$.

The conjugacy class $[a]$ of $a \in F_n$ is {\em carried by $[[F^i]]$}
if $F^i$ contains a representative of $[a]$. Sometimes we say that $a$
is carried by $F^i$ when we really mean that $[a]$ is carried by
$[[F^i]]$.  If $G$ is a marked graph and $G_r$ is a subgraph of $G$
such that $[[F^i]] = [[\pi_1(G_r)]]$, then $[a]$ is carried by
$[[F^i]]$ if and only if the circuit in $G$ that represents $[a]$ is
contained in $G_r$.  We say that {\em an abstract line $ \ell$ is
carried by $[[F^i]]$} if its realization in $G$ is contained in $G_r$
for some, and hence any, $G$ and $G_r$ as above.  Equivalently, $\ell$
is the limit of periodic lines corresponding to $[c_i]$ where each
$[c_i]$ is carried by $[[F^i]]$.  A collection $W$ of abstract lines
and conjugacy classes in $F_n$ is carried by a free factor system $\F
=\{[[F^1]],\ldots,[[F^k]]\}$ if each element of $W$ is carried by some
$F^i$.

There is a partial order $\sqsubset$ on free factor systems generated
by inclusion. More precisely, $[[F^1]] \sqsubset [[F^2]]$ if $F^1$ is
conjugate to a free factor of $F^2$ and $\F_1 \sqsubset \F_2$ if for
each $[[F^i]] \in \F_1$ there exists $[[F^j]] \in \F_2$ such that
$[[F^i]] \sqsubset [[F^j]]$.

The {\em complexity} of a \ffs\ is defined on page 531 of
\cite{bfh:tits1}.  Rather than repeat the definition, we recall the
three properties of complexity that we use.  The first is that if
$\F_1 \sqsubset \F_2$ for distinct free factor systems $\F_1$ and
$\F_2$ then the complexity of $\F_1$ is less than the complexity of
$\F_2$.  This is immediate from the definition.  The second is
Corollary~2.6.5. of \cite{bfh:tits1}.

\begin{lemma}\label{ffs lemma}
For any collection $W$ of abstract lines there is a unique free factor
system $\F(W)$ of minimal complexity that carries every element of
$W$.  If $W$ is a single element then $\F(W)$ has a single element.
\end{lemma}

The third is an immediate consequence of the uniqueness of $\F(W)$.

\begin{cor} \label{ffs invariance} If a collection $W$ of abstract lines and 
conjugacy classes in $F_n$ is $\oone$-invariant then $\F(W)$ is
$\oone$-invariant.
\end{cor} 
 
Further details on free factor systems can be found in section~2.6 of
\cite{bfh:tits1}.

\
\subsection{Relative Train Track Maps} \label{rtt review} 
In this section we review and set notation for relative train track
maps as defined in \cite{bh:tracks}.

Suppose that $G$ is a marked graph and that $\fG$ is a homotopy
equivalence representing $\oone \in \Out(F_n)$.  A {\it filtration} of
$G$ is an increasing sequence  $ \filt$ of subgraphs, each of whose components contains at least one edge.  If $f(G_i)
\subset G_i$ for all $i$ then we say that $\fG$ {\it respects the
filtration} or that the {\em filtration is $f$-invariant}.  A path or
circuit has {\em height} $r$ if it is contained in $G_r$ but not
$G_{r-1}$.  A lamination has height $r$ if each leaf in its
realization in $G$ has height at most $r$ and some leaf has height
$r$.

The $r^{th}$ {\em stratum} $H_r$ is defined to be the closure of $G_r
\setminus G_{r-1}$.  To each stratum $H_r$ there is an associated
square matrix $M_r$, called the {\em transition matrix for $H_r$},
whose $ij^{th}$ entry is the number of times that the $i^{th}$ edge
(in some ordering of the edges of $H_r$) crosses the $j^{th}$ edge in
either direction.  By enlarging the filtration, we may assume that
each $M_r$ is either irreducible or the zero matrix.  We say that
$H_r$ {\em is an irreducible stratum} if $M_r$ is irreducible and is
{\em a zero stratum} if $M_r$ is the zero matrix.

If $M_r$ is irreducible and the Perron-Frobenius eigenvalue of $M_r$
is $1$, then $M_r$ is a permutation matrix and $H_r$ is {\em
non-exponentially growing} or simply {\em NEG}.  After subdividing and
replacing the given NEG stratum with a pair of NEG strata if
necessary, the edges $\{E_1,\ldots,E_l\}$ of $H_r$ can be oriented and
ordered so that $f(E_i) = E_{i+1}u_{i}$ where $u_{i} \subset G_{r-1}$
and where indices are taken mod $l$.  We always assume that edges in
an NEG stratum have been so oriented and ordered.  If each $u_i$ is
trivial then $f^l(E_i) = E_{i}$ for all $i$ and we say that $E_i$
[$H_r$] is a {\em periodic edge [stratum]} with period $l$ or a {\em
fixed edge [stratum] } if $l=1$.  If each $u_i$ is a Nielsen path then
the combinatorial length of $f_\#^{k}(E_i)$ is bounded by a linear
function of $k$ for all $i$ and we say that $E_i$ [$H_r$] is a {\em
linear edge [stratum]}.

If $M_r$ is irreducible and if the Perron-Frobenius eigenvalue of
$M_r$ is greater than 1 then $H_r$ is an {\em exponentially growing
stratum} or simply an {\em \EG\ stratum}.   If $H_r$ is \eg\ and $\alpha \subset G_{r-1}$  is a non-trivial path 
with endpoints in $H_r \cap G_{r-1}$ then we say that $\alpha$ is a   {\em connecting path  for $H_r$}.   If  $H_r$ is \eg\ and $\sigma$ is a path with height $r$ then we sometimes say that {\em $\sigma$ has \eg\ height}.

A {\em direction} $d$ at $x \in G$ is the germ of an initial segment
of an oriented edge (or partial edge if $x$ is not a vertex) based at
$x$.  There is an $f$-induced map $Df$ on directions and we say that
$d$ is {\em a periodic direction} if it is periodic under the action
of $Df$; if the period is one then $d$ is a {\em fixed direction}.
Thus, the direction determined by an oriented edge $E$ is fixed if and
only if $E$ is the initial edge of $f(E)$.  Two directions with the
same basepoint belong to the same {\em gate} if they are identified by
some iterate of $Df$.  If $x$ is a periodic point then the number of
gates based at $x$ is equal to the number of periodic directions based
at $x$.

A \emph{turn} is an unordered pair of directions with a common base
point.  The turn is \emph{nondegenerate} if is defined by distinct
directions and is \emph{degenerate} otherwise.  If $E_1E_2 \ldots
E_{k-1} E_k$ is the edge path associated to a path $\sigma$, then we
say that {\it $\sigma$ contains the turns $(\bar E_i,E_{i+1})$ for $1
\le i \le k-1$}.  A turn is {\it illegal} with respect to $\fG$ if its
image under some iterate of $Df$ is degenerate; a {turn is \it legal}
if it is not illegal.  Equivalently, a turn is legal if and only if it
is defined by directions that belong to distinct gates.  A {\it path
or circuit $\sigma \subset G$ is legal} if it contains only legal
turns.  If $\sigma \subset G_r$ does not contain any illegal turns in
$H_r$, meaning that both directions correspond to edges of $H_r$, then
$\sigma$ is {\em $r$-legal}.  It is immediate from the definitions
that $Df$ maps legal turns to legal turns and that the restriction of
$f$ to a legal path is an immersion.

We recall the definition of relative train track map from page 38 of
\cite{bh:tracks}.
\begin{definition}  \label{rttDef} 
A homotopy equivalence $\fG$ representing $\oone$ is a \emph{\rtt} \
if it satisfies the following conditions for every \eg\ stratum $H_r$
of an $f$-invariant filtration $\F$.
\begin{itemize}
\item[(RTT-$i$)] $Df$ maps the set of directions in $H_r$ with basepoints
at  vertices to itself; in particular every turn with one direction in
$H_r$ and the other in $G_{r-1}$ is legal.
\item[(RTT-$ii$)] If $\alpha$ is a connecting path for $H_r$ then    $f_\#(\alpha)$ is  
a connecting path for $H_r$; in particular,  $f_\#(\alpha)$ is nontrivial.   
\item[(RTT-$iii$)] If $\alpha \subset G_r$ is $r$-legal then
$f_\#(\alpha)$ is $r$-legal.
\end{itemize}
\end{definition} 

\begin{remark} 
If $\fG$ is a relative train track map, then so is $f^k$ for $k>0$.
\end{remark}

A subgraph $C$ of $G$ is {\em wandering} if $f^k(C)$ is
contained in the closure of $G \setminus J$ for all $k > 0$; otherwise $C$ is {\em non-wandering}.   Each edge in a wandering subgraph is contained in a zero stratum.    If $C$ is a component of a filtration element then $C$ is non-wandering if  and only if $f^i(C) \subset C$ for some $i > 0$.

 \begin{remark} \label{rtt2} Suppose that $H_r$ is an \eg\ stratum and   that $\fG$ satisfies (RTT-$i$) for $H_r$.  Then $f(H_r \cap G_{r-1}) \subset  H_r \cap G_{r-1}$ and verifying (RTT-$ii$) for $H_r$ reduces to showing that $f_\#(\alpha)$ is nontrivial for each   connecting path $\alpha$ for $H_r$.  Suppose that the component $C$ of $G_{r-1}$ that contains $\alpha$ is non-wandering.     In checking (RTT-$ii$)  there is no loss in replacing $f$ by $f^i$ so we may assume that $f(C) \subset C$.    If $f$ permutes the elements of the finite set  $H_r \cap C$, then $f_\#(\alpha)$ is nontrivial for each      $\alpha \subset  C$.    If $f$ identifies two of these points then they can be connected by an arc $\alpha \subset C$ such that $f_\#(\alpha)$ is trivial (because $f|C :C \to C$ is a homotopy equivalence).   This proves that if $C$ is a non-wandering component of $G_{r-1}$ then  (RTT-$ii$) holds for all $\alpha \subset C$ if and only if $H_r \cap C \subset \Per(f)$.  
\end{remark}

The most common applications of the relative train track properties
are contained in following lemma.

\begin{lemma} \label{r legal} Suppose that $\fG$ is a \rtt\ and  that $H_r$ is 
an \EG\ stratum.
\begin{enumerate}
\item Suppose that a vertex $v$ of $H_r$ is contained in a component
$C$ of $G_{r-1}$ that is non-contractible or more generally satisfies
$f^i(C) \subset C$ for some $i > 0$.  Then $v$ is periodic and has at
least one periodic direction in $H_r$.
\item If $\sigma$ is an $r$-legal circuit or path of height $r$ with
endpoints, if any, at vertices of $H_r$ then the decomposition of
$\sigma$ into single edges in $H_r$ and maximal subpaths in $G_{r-1}$
is a splitting.
\end{enumerate} 
\end{lemma}    

\proof The first item  follows from Remark~\ref{rtt2}.    The second item is contained
in Lemma 5.8 of \cite{bh:tracks}.  \endproof

\begin{lemma} \label{two directions for EG}
If $H_r$ is an EG stratum of a \rtt \ $\fG$ and $x \in H_r $ is either
a vertex or a periodic point then there is a legal turn in $G_r$ that
is based at $x$.  In particular, there are at least two gates in $G_r$
that are based at $x$.
\end{lemma}

\proof There exists $j > 0$ and a point $y$ in the interior of an edge
$E$ of $H_r$ so that $x=f^j(y)$. By (RTT-$iii$), $x$ is in the
interior of an $r$-legal path.  Moreover the turn at $x$ determined by
this path is legal by Properties (RTT-$iii$) and (RTT-$i$).  \endproof

The following lemma describes \ipNp s with  EG height. 

\begin{lemma} \label{nielsen paths in egs}  Suppose that $\fG$ is a \rtt\ and 
that $H_r$ is an \EG\ stratum.
\begin{enumerate}
\item There are only finitely many  indivisible  periodic Nielsen paths   of 
height $r$.
\item If $\sigma$ is an indivisible periodic Nielsen path of height
$r$ then $\sigma = \alpha   \beta$ where $\alpha$ and $\beta$ are
$r$-legal paths that begin and end with directions in $H_r$ and the
turn $(\bar \alpha,   \beta)$ is illegal.   Moreover, if $\bar \alpha(k)$
and $\beta(k)$ are the initial segments of $\bar \alpha$ and $\beta$ that
are identified by $f_\#^k$ then the $\bar \alpha(k)   \beta(k)$'s are an
increasing sequence of subpaths whose union is the interior of
$\sigma$.

\item An indivisible periodic Nielsen path $\sigma$ of height $r$ has
period $1$ if and only if the initial and terminal directions of
$\sigma$ are fixed.
\end{enumerate}
\end{lemma} 
\proof  Item (1) and the first part of   (2) is  contained in the statement of Lemma~5.11 of \cite{bh:tracks}.    The moreover part of (2) is contained in the proof of that lemma.   We now turn our attention to (3).

Suppose that $\sigma$ is an indivisible periodic Nielsen path of
height $r$ and period $p$ and that $\sigma = \alpha  \beta$ as in
(2).  After subdividing at the endpoints of $f_\#^j(\sigma)$ for $0
\le j \le p-1$ (which clearly preserves the property of being a
relative train track map) we may assume that the endpoints of $\sigma$
are vertices.  Since $\alpha$ and $\bar \beta$  are $r$-legal and begin  with  edges in
$H_r$, Lemma~\ref{r legal}(2) implies that $Df$ maps the initial
directions of $\alpha$ and $\bar \beta$ to the initial direction of $f_\#(\alpha)$ and $f_\#(\bar \beta)$ respectively.    The
only if part of (3) therefore follows from the fact that
$f_\#(\sigma)$ is obtained from $f_\#(\alpha)$ and $f_\#(\bar \beta)$ by
cancelling their maximal common terminal segment.

Assume now that the first edge $E$ of $\alpha$ determines a fixed
direction and write $f_\#(\sigma) = \alpha_1  \beta_1$ as in
(2).   Then $E$   is  also the first edge in
$\alpha_1$.    The moreover part
of (2) implies that both $\alpha$ and $\alpha_1$ are initial segments
of $f^{Np}_\#(E)$ for all sufficiently large $N$.  Thus either
$\alpha$ is an initial segment of $\alpha_1$ or $\alpha_1$ is an
initial segment of $\alpha$.  The same argument applies to $\bar \beta$ and
$\bar \beta_1$.  The difference between the number of $H_r$
edges in $f_\#^{Np}(\alpha)$ and the number of $H_r$ edges in
$f_\#^{Np}(\bar \beta)$ is independent of $N$ and similarly for $\alpha_1$
and $\bar \beta_1$.  It follows that either $\alpha \subset \alpha_1$ and
$\bar \beta \subset \bar \beta_1$ or $\alpha_1 \subset \alpha$ and $\bar \beta_1
\subset \bar  \beta$.  For concreteness assume the former.

We claim that $\alpha_1 = \alpha$.  If not then $\alpha_1 =
\alpha\gamma$ for some non-trivial $\gamma$.  The path $\alpha_2 :=
\bar \beta\gamma$ is a subpath of $f^{Np}_\#(\bar \beta)$ for large $N$ and so
is $r$-legal.  The path $\alpha_2  \beta_1 = [(\bar \beta \bar
\alpha)(\alpha_1  \beta_1)]$ is a non-trivial periodic Nielsen
path with exactly one illegal turn in $H_r$ and is therefore
indivisible.  By construction, $\alpha_2$ and $\bar \beta$, and hence
$\alpha_2$ and $\bar \beta_1$, have a non-trivial common initial subpath in
$H_r$.  This implies as above that $\alpha_2$ and $\bar \beta_1$ are
initial subpaths of a common path $\delta$. They cannot be equal so
one is a proper initial subpath of the other.  But then the difference
between the number of $H_r$ edges in $f_\#^{Np}(\alpha_2)$ and the
number of $H_r$ edges in $f_\#^{Np}(\bar \beta_1)$ grows exponentially in
$N$ in contradiction to the fact that $\alpha_2   \beta_1$ is a
periodic Nielsen path.  This contradiction verifies the claim that
$\alpha = \alpha_1$.  The symmetric argument implies that $\beta =
\beta_1$ so $p =1$.  \endproof

We extend Lemma~\ref{nielsen paths in egs}(1) as follows.

\begin{lemma}\label{finitely many endpoints}  
For any \rtt\ $\fG$ there are only finitely many points that are
 endpoints of an \ipNp\ $\sigma$.
\end{lemma}

\proof Let $S_r$ be the set of endpoints of \ipNp s $\sigma$ with
height at most $r$.  We prove that $S_r$ is finite for all $r$ by
induction on $r$.  Since $S_0$ is empty we assume that $S_{r-1}$ is
finite and prove that $S_r$ is finite.  This is obvious if $H_r$ is a
zero stratum (because $S_r = S_{r-1}$) and follows from Lemma~\ref{nielsen paths in egs}(1) if
$H_r$ is EG.  We may therefore assume that $H_r$ is NEG and that there
is an indivisible Nielsen path $\sigma$ of height $r$.  Lemma 4.1.4 of
\cite{bfh:tits1} implies that $H_r$ is not periodic (because $\sigma$  is indivisible) and that edges of
$\sigma$, other than perhaps the first and last, are not contained in
$H_r$.  The only periodic points in $H_r$ are vertices so we may
assume that $\sigma = Eu$ where $E$ is an edge of $H_r$ and $u \subset
G_{r-1}$.  If $\sigma' = Eu'$ is another such \ipNp\ then $\bar u u'$
is an \ipNp\ with height less than $r$ and so its endpoint are contained in $S_{r-1}$.  Thus $S_r$ is obtained from
$S_{r-1}$ by adding at most two points for each edge in $H_r$.  This
completes the inductive step.  \endproof

An  irreducible matrix $M$ is {\em aperiodic} if $M^k$ has all positive
entries for some $k \ge 1$. For example, if some diagonal element of
$M$ is non-zero then $M$ is aperiodic.  If $H_r$ is an EG stratum of a
\rtt\ $\fG$ and if the transition matrix $M_r$ is aperiodic, then $H_r$
is said to be an {\em aperiodic EG stratum}.  For each aperiodic EG
stratum there is a unique  $\oone$-invariant attracting lamination $\Lambda_r\in \L(\oone)$ of height
$r$.  If every EG stratum is aperiodic then every element of $\L(\oone)$ is
related to an EG stratum in this way.  See Definition 3.1.12 of
\cite{bfh:tits1} and the surrounding material for details.

The following lemma produces rays and lines associated to aperiodic EG strata of
a \rtt.

\begin{lemma} \label{lam2}
Suppose that $H_r$ is an EG stratum of a \rtt\ $\fG$, that $\ti f:
\Gamma \to \Gamma$ is a lift and that $\ti v \in \Fix(\ti f)$.
\begin{enumerate}
\item 
If $E$ is an oriented edge in $H_r$ and $\ti E$ is a lift that
determines a fixed direction at $\ti v$, then there is an $r$-legal
ray $\ti R \subset \Gamma$ that begins with $\ti E$, intersects
$\Fix(\ti f)$ only in $\ti v$ and that converges to an attractor $P
\in \Fix(\hat f)$ whose accumulation set  is the  (necessarily $\oone$-invariant) element $
\Lambda_r$ of
$\L(\oone)$ whose realization in $G$ has height $r$. 
\item Suppose that $E' \ne E$ is another oriented edge in $H_r$, that
$\ti E'$ determines a fixed direction at $\ti v$ and that $\ti R'$ is
the ray associated to $\ti E'$ as in (1).  Suppose further that the
turn $(\bar E, E')$ is contained in the path $f^k_\#(E'')$ for some $k
\ge 1$ and some edge $E''$ of $H_r$.  Then the line $\ti R^{-1} \ti
R'$ is a generic leaf of $\ti \Lambda_r$.
\end{enumerate}
\end{lemma}

\proof    Lemma~\ref{r legal}(2) and (RTT-i) imply that $\ti f(\ti E) =
\ti E \cdot \ti \mu_1$ for some non-trivial $r$-legal subpath $\ti
\mu_1$ of height $r$ that ends with an edge of $H_r$.  Applying
Lemma~\ref{r legal}(2) again, we have $\ti f^2_\#(\ti E) = \ti E \cdot
\ti \mu_1\cdot \ti \mu_2$ for some $r$-legal subpath $\ti \mu_2$ of
height $r$ that ends with an edge of $H_r$.  Iterating this produces a
nested increasing sequence of paths $\ti E \subset \ti f(\ti E)
\subset \ti f^2_\#(\ti E) \subset \ldots$ whose union is a ray $\ti R$
that converges to some attractor $P \in \Fix_N(\hat f)$. 

The transition matrix for $H_r$ is aperiodic since $f_\#(E)$ contains $E$.  Let $\beta$ be    a generic leaf of $\ti \Lambda_r$.  Following Definition~3.1.7 of \cite{bfh:tits1}, a path of the form $f_\#^k(E_i)$ where $k \ge 0$ and $E_i$ is an edge of $H_r$ is called a {\em tile}.  By Lemma~\ref{r legal} and by construction, every tile occurs infinitely often in $\ti R$.  By Corollary~3.1.11 of \cite{bfh:tits1}, each subpath of $\beta$ is contained in some tile and each tile occurs as a subpath of $\beta$. It follows that the accumulation set of $P$ is equal to the weak limit of $\beta$ which is $\Lambda^+_r$ because $\beta$ is generic.  This proves (1).

 Assuming now the notation of (2), each finite subpath of $\ti R^{-1} \ti
R'$ is contained in  a tile.   This implies, as in the previous case, that $\ti R^{-1} \ti
R'$  is contained in $\Lambda_r$ and that    $\ti R^{-1} \ti
R'$ is birecurrent.     Lemma~3.1.15 of \cite{bfh:tits1} implies that  $\ti R^{-1} \ti R'$ is generic.   
\endproof

As noted in section~\ref{laminations}, $\Out(F_n)$ acts on the set of
laminations in $G$.  The stabilizer $\Stab(\Lambda)$ of a lamination
$\Lambda$ is the subgroup of $\Out(F_n)$ whose elements leave $\Lambda$ invariant.
We recall Corollary 3.3.1 of \cite{bfh:tits1}.

\begin{lemma} \label{stabilizer} For each $\Lambda \in \L(\oone)$, there 
is a
homomorphism $PF_{\Lambda}:$
Stab($\Lambda) \to \mathbb
Z$ such that  $\psi \in$ Ker($PF_{\Lambda}$) if 
and
only if $\Lambda \not \in {\cal L}({\psi})$ and  
$\Lambda \not \in  {\cal
L}({ \psi}^{-1})$.  
\end{lemma}

We refer to $PF_{\Lambda}$ as the {\em expansion factor homomorphism}
associated to $\Lambda$.  The notation is chosen to remind readers
that the expansion factor is realized as the logarithm of a Perron
Frobenius eigenvalue in a natural way.
       
\subsection{Modifying \rtt s} \label{s:modification}  

To simplify certain arguments in section~\ref{rotOutAut} and as a
 step toward our ultimate existence theorem (Theorem~\ref{cs
exists}), we add properties to the \rtt s produced in
\cite{bh:tracks}.  We can not simply quote results from
\cite{bfh:tits1} because, unlike in \cite{bfh:tits1}, here we do not
allow iteration and we make no assumptions on $\oone$.

We need some further notation.  The {\em set of periodic points of
$f$} is denoted $\Per(f)$.  A path $\alpha$ is {\em pre-trivial} if
$f_\#^k(\alpha)$ is trivial for some $k > 0$.   We say that a {\em nested
sequence $\cal C$ of free factor systems $\F^1 \sqsubset \F^2
\sqsubset \ldots \sqsubset \F^m$ is  realized} by a \rtt\ $\fG$
and filtration $\filt$ if  each $\F^j$ is realized by some  
$G_{l(j)}$.      For any finite graph $K$, the {\em core
of $K$} is the subgraph of $K$ consisting of edges that are crossed by
some circuit in $K$.  If $K_0 = K$ and $K_{i+1}$ is inductively
obtained from $K_i$ by removing all edges with at least one valence
one endpoint, then $K_i$ is the core of $K$ for all sufficiently large
$i$.

The transition matrix $M_r$ for an  \eg\ stratum $H_r$ of a topological representative  $\fG$ of $\phi$ has a Perron-Frobenius eigenvalue $\lambda_r > 1$.    The set $\{\log(\lambda_r) : H_r \mbox{ is \eg}\}$, listed in non-increasing order, is denoted $\PF(f)$.  (In \cite{bh:tracks} this set is denoted $\Lambda(f)$ but $\Lambda$ is now usually reserved for laminations.)     The  set $\{ \PF_{\Lambda}(\phi): \Lambda \in  \L(\phi)\} $ of expansion factors for $\phi$, listed in non-increasing order, is denoted $\EF(\phi)$.  Then $\PF(f) \ge \EF(\phi)$ in the lexicographical order for all $f$ representing $\phi$ and  $\PF(f) = \EF(\phi)$ if $\fG$ is a \rtt\  representing $\phi$ by Lemma~3.3 of \cite{bfh:tits1}.    

The number of \iNp s for $\fG$ with height $r$ is denoted $N_r(f)$.

\begin{remark}  \label{just needRTT2}  If $\fG$ satisfies (RTT-$i$) and if $\PF(f) = \EF(\phi)$ then $\fG$ satisfies  (RTT-$iii$) by Lemma~5.9 of \cite{bh:tracks}.    Thus any topological representative that satisfies  (RTT-$i$),  (RTT-$ii$) and   $\PF(f) = \EF(\phi)$   is a \rtt.
\end{remark}

\begin{lemma}\label{unchanged Nielsen paths}  Suppose that $\fG$ and $f' : G' \to G'$ are \rtt s with \eg\ strata $H_r$ and $H'_s$ respectively and that $p:G \to G'$ is a homotopy equivalence such that 
\begin{enumerate}
\item $p(G_r) = G'_s$ and $p$ induces a bijection between the edges of $H_r$ and the edges of $H'_s$.
\item $p_\# f_\#(\sigma)  =  f'_\#p_\#(\sigma)$ for all paths $\sigma \subset G_r$ with endpoints at vertices.
\end{enumerate}
Then $p_\#$ induces  a bijection between the indivisible periodic Nielsen paths in $G$ with height $r$ and the indivisible periodic Nielsen paths in $G'$ with height $s$.
\end{lemma}

\proof     If $\sigma$ is a periodic Nielsen path for $f$
  then $  f'_\#p_\#(\sigma)=  p_\# f_\#(\sigma) = p_\# (\sigma)$ which proves that  $\sigma' =  p_\# (\sigma)$ is a  periodic Nielsen path for $f'$.     If $\sigma$ has height $r$ and is indivisible   then
$\sigma$ begins and ends in $H_r$ and has exactly one illegal turn in
$H_r$ by Lemma~\ref{nielsen paths in egs}.   Our hypotheses imply that the same is true of
$\sigma'$ which implies that $\sigma'$ is indivisible.  For the
converse note that for any path \ $\sigma' \subset G'_s$ that begins with $p(E_i) \subset H'_s$ and ends with $p(E_j) \subset H'_s$  there is a unique path $\sigma$ that begins with $E_i$, ends with $E_j$   and that satisfies $p_\#(\sigma) = \sigma'$.
By the same reasoning $\sigma$ is an \iNp\ if $\sigma'$ is.  
\endproof

We next recall the {\em sliding} operation from \cite{bfh:tits1}.
Suppose that $H_s$ is a non-periodic NEG stratum with edges
$\{E_1,\ldots,E_m\}$ satisfying $f(E_l) = E_{l+1}u_l$ for paths $u_l
\subset G_{s-1}$ where $1 \le l \le m$ and indices are taken mod $m$.
Choose $1 \le i \le m$ and let $\tau$ be a path in $G_{s-1}$ from the
terminal endpoint $v_i$ of $E_i$ to some vertex $w_i$.  Roughly
speaking, we use $\tau$ to continuously change the terminal endpoint
of $E_i$ from $v_i$ to $w_i$ and to mark the new graph.  

More
precisely, define a new graph $G'$ from $G$ by replacing $E_i$ with an
edge $E_i'$ that has terminal vertex $w_i$ and that has the same
initial vertex as $E_i$.  
 There are homotopy equivalences $p : G \to
G'$ and $p': G' \to G$ that are the identity on the common edges of
$G$ and $G'$ and that satisfy $p(E_i)=E'_i\bar \tau$ and $p'(E_i') =
E_i \tau$.  Use $p$ to define the marking on $G'$ and define $f' : G'
\to G'$ on edges by tightening $pfp'$.  Complete details can be found
in section~5.4 of \cite{bfh:tits1}.

\begin{lemma}\label{sliding prelim} With notation as above. 
\begin{itemize}
\item $f':G' \to G'$ is a \rtt.
\item $f'|G_{s-1} =  f|G_{s-1}$
 \item If $m = 1$ then $f'(E'_1) = E'_1[\bar \tau u_1 f(\tau)]$ 
\item If $m \ne 1$ then $f'(E_{i-1}) = E'_i[\bar \tau u_{i-1}],
f'(E'_i) = E_{i+1}[u_i f(\tau)]$ and $f'(E_j) = E_j u_j $ for $j \ne
i-1,i$.
\item For each \eg\ stratum $H_r$, $p_\#$ defines  a bijection between the set of  the indivisible periodic Nielsen paths in $G$ with height $r$ and the indivisible periodic Nielsen paths in $G'$ with height $r$.
\end{itemize}
\end{lemma}

\proof If $m = 1$ then this is contained in Lemma 5.4.1 of
\cite{bfh:tits1}.  The argument for $m > 1$ is a straightforward
extension of $m = 1$ case and we leave the details to the reader.
\endproof

  \begin{notn} Suppose that $u < r$ and that 
  \begin{enumerate}
  \item $H_u$ is irreducible.
  \item $H_r$ is \eg\ and  each component of $G_r$ is non-contractible.
  \item   for each $u < i < r$,   $H_i$ is a zero stratum that is a component of $G_{r-1}$  and each vertex of $H_i$ has valence at least two in $G_r$. 
  \end{enumerate}
We say that each $H_i$ is {\em enveloped by $H_r$} and write $H_r^z = \cup_{k=u+1}^r H_k$.  It is often  convenient to treat $H_r^z$ as a single unit.   
  \end{notn}

\begin{thm} \label{rtt existence}  
For every $\oone \in \Out(F_n)$  there is
a \rtt\ $\fG$ and filtration that represents $\oone$  
 and satisfies the following properties.  
 \begin{enumerate}
  \item [(V)] The endpoints of all indivisible periodic Nielsen paths are
vertices.
\item [(P)] If  a stratum $ H_m \subset \Per(f)$ is a forest then there exists a filtration element $G_j$ such  $\F(G_j) \ne \F(G_l \cup H_m)$ for any $G_l$.  (See also items (1) and (5) of Lemma~\ref{formerly remarks}.)
 \item [(Z)]  Each zero stratum $H_i$ is enveloped by an  \eg\ stratum $H_r$.    Each vertex in $H_i$ is contained in $H_r$ and has link  contained in $H_i \cup H_r$.
 \item [(\noneg)] The terminal endpoint   of an edge in a non-periodic NEG stratum  $H_i$  is periodic and is contained in a filtration element of height less than $i$ that is its own core.
  \item [(F)]   The core of each filtration element is a filtration element.
  \end{enumerate}
  Moreover, if  $\cal C$ is a   nested sequence 
of non-trivial $\oone$-invariant free factor systems then we may assume that $\fG$  realizes $\cal C$.
\end{thm}

Before proving  Theorem~\ref{rtt existence} we record some useful observations.

\begin{lemma}  \label{formerly remarks}Suppose that $\fG$ is a \rtt\   representing   $\phi$ with filtration $\filt$.   
\begin{enumerate} 
\item Suppose that $\fG$ satisfies (P),  that  $ H_m \subset \Per(f)$ is a forest and that $v$  has valence one in $H_m$.  Then $v \in G_j$  for some $j < m$.     If in addition,  $\fG$ satisfies (F) then we may choose  $G_j$ to be  is its own core.    
\item If  $\fG$ satisfies (F) and $H_r$ is \eg\ then $G_r$ is its own core.
\item  If $\fG$ satisfies (P) and (\noneg)  then    each contractible component of a filtration element is composed of zero strata. 
\item   If $\fG$ satisfies (Z) and (\noneg)  and if $G_k$ is a filtration element that is its own core, then every vertex in $G_k$ has at least two gates in $G_k$.
\item  If $\fG$ satisfies (P)  then no component of a filtration element $G_m$ is  a  forest in $ \Per(f)$. 
\end{enumerate}
\end{lemma}

\proof     Suppose that $f:G\to G$ satisfies (P) and  that $H_m\subset \Per(f)$ is
a forest.        By (P) there    exists  $j$ so that $\F(G_j) \ne \F(G_l \cup H_m)$ for any $l$;  if  $f:G\to G$ also satisfies (F) then we may assume that $G_j$ is its own
core.    In particular,  $\F(G_j\cup H_m) \ne  \F(G_j)$.  Thus some, and hence every, vertex of valence one in $H_m$
must have valence at least two in $G_j\cup H_m$. This proves (1).
 
 Item (5) follows from the first part of (1).
 
If $H_r$ is \eg\ then $G_r$ is the smallest filtration element that contains the attracting lamination associated to $H_r$.  This proves (2).

  For (3), suppose that $C$ is a contractible component of some $G_i$.       If $C$ contains an edge in an irreducible stratum then it is non-wandering and by (\noneg) the lowest stratum  $H_j$  in $C$ is either periodic or \eg.   But  $H_j$ can not be periodic by item (1) of this lemma,  and cannot be \eg\ by  Lemma~\ref{two directions for EG}.  Thus every edge in $C$ is contained in a zero stratum. 
  
   The proof of  (4) is by induction on $k$.   Suppose that $G_k$ is a filtration element that is its own core and that $v$  is a vertex in $G_k$.   If   no illegal turns  in $G_k$  are based at $v$ then the number of gates based at $v$ is at least  the valence of $v$ in $G_k$, and so is at least two.  We may therefore assume that   there is an   illegal turn $(d_1,d_2)$ in $G_k$ based at $v$.  By (Z), either both   $d_1$ and $d_2$ belong to the same \eg\ stratum or one of the $d_i$'s is the terminal end of a non-fixed \noneg\ edge.   Lemma~\ref{two directions for EG},   (\noneg)    and the inductive hypothesis imply that  in both of these cases $v$ has at least  two gates in $G_k$.  
   \endproof

\noindent{\em Proof of Theorem~\ref{rtt existence}.} By Lemma~2.6.7 of
\cite{bfh:tits1}, there is a \rtt\ $\fG$ and filtration
$\filt$ that represents $\oone$ and realizes $\cal C$.     (In the statement of Lemma~2.6.7, $\cal
C$ is replaced by a single invariant free factor system $\F$.  The
more general case that we use is explicitly included in the proof of
that lemma.)     Lemma~\ref{finitely many
endpoints} implies that (V) can be arranged via a finite subdivision. 
For convenient reference we divide the rest of the proof into steps.  Changes are made to $\fG$ in these steps but we start each step by referring to the current \rtt\ as $\fG$.  

   If $\cal C$ has not already been specified, let $\cal C$ be the nested sequence of   free factor systems determined by the $\F(G_r)$'s.     For possible future application   we will prove the following statement in place of (P).
    \begin{itemize}
\item [(P$_{\cal C}$)] If  a stratum $ H_m \subset \Per(f)$ is a forest then there exists  $\F^j \in \cal C$ that is not realized by  $G_l \cup H_m$ for any $G_l$. 
\end{itemize}
If $\cal C$ is determined by the $\F(G_r)$'s then (P) and (P$_{\cal C}$) are the same but otherwise the latter is  stronger than the former.

\begin{remark}  \label{same Nielsen paths} For reference in the proof of Theorem~\ref{comp sp exists} we record the following  property of the remainder of our construction.  If $\fG$ is the   \rtt\  as it is now  and  $g:G' \to G'$ is the ultimate modified \rtt\ produced by the six steps listed below, then there is a bijection $H_r \to H'_{s}$ between the \eg\ strata of $\fG$ and the \eg\ strata of  $g:G' \to G'$ such that 
 \begin{enumerate}
 \item[(a)]$H_r$ and $H'_{s}$ have the same number of edges.
 \item [(b)]  $N_r(f) = N_s(g)$. \end{enumerate}
The five moves used in our construction are valence two homotopies away from \eg\ strata,  sliding (which is defined following Lemma~\ref{unchanged Nielsen paths}), reordering of strata, tree replacements (see step 2) and   collapsing forests in \noneg\ and zero strata.   Item (a) will be obvious as will   (b) for the valence two homotopies and  the reordering of strata.    For the remaining three moves    (b) will follow from  Lemmas~\ref{sliding prelim} and \ref{unchanged Nielsen paths}.
\end{remark}

\vspace{.1in} \noindent{\bf(Step 1 :  A   weak form of (\noneg))}    Suppose that $H_s$ is
a non-periodic NEG stratum with edges $\{E_1,\ldots,E_m\}$ satisfying
$f(E_i) = E_{i+1}u_i$ for paths $u_i \subset G_{s-1}$ where indices
are taken mod $m$.    Our goal in this step is to arrange 
\begin{enumerate}
\item The endpoints of an \noneg\ edge are periodic. 
\end{enumerate}
Care is  taken that  no vertices with valence one in $G$  are created.

We first arrange that   the  terminal endpoint $v_1$ of $E_1$ is either periodic or has valence at least three in $G$.  If this is not already the case, let $E$  be the unique oriented edge  of $G$, other than $E_1$,  whose initial endpoint is  $v_1$.   Since $H_s$ is not a periodic stratum,   $E \subset G_{s-1}$.  Lemma~\ref{two directions for EG}  implies that $E$ does not belong to an \eg\ stratum.  Perform a valence two homotopy as defined on page 13 of  \cite{bh:tracks}.  There are two steps to this operation.  The first is to modify  $f$ so that no vertex is mapped to $v_1$.  The second is to amalgamate $E_1E$ into a single edge named $E_1$ by removing $v_1$ from the list of vertices.      If $\sigma$ is a path, neither of  whose endpoints  are  mapped to $v_1$ by the original map, then occurrences of $E_1$ and $E$ in the original edge path for $f_\#(\sigma)$  always occur as $E_1E$ or its inverse;   the new edge path for $f_\#(\sigma)$ is obtained from the original  by removing all occurrences of $E$.    Since $v_1$ is not incident to any \eg\ stratum, this applies to any path $\sigma$ with endpoints in   \eg\ strata.       It follows that the new map (still called $\fG$) is still a \rtt,     satisfies (V) and realizes  $\cal C$.     It also follows that this operation does not change the number of edges in any \eg\ stratum and does not change the number of \iNp s with height corresponding to any \eg\ stratum.  After finitely many valence two homotopies we may assume that the  terminal endpoint $v_1$ of $E_1$ is either periodic or has valence at least three in $G$. 
 
 The component $C_1$ of $G_{s-1}$ that contains $v_1$  does not wander and so contains a periodic vertex $w_1$. Choose a path $\tau \subset G_{s-1}$   from $v_1$ to   $w_1$ in $C$
and slide to change the terminal endpoint of $E_1$ to $w_1$.  Because of our previous move, $v_1$ still has valence at least two in $G$.

After repeating these operations finitely many times we have arranged  (1) and not created any valence one vertices.

\vspace{.1in} \noindent{\bf(Step 2 :  A weak form of (Z))}  In this step we prove most of (Z).  The one missing item is that we only show that the components of $G_r$ are non-wandering instead of showing that they are non-contractible. 

The union of the non-wandering components of  a filtration element  $G_i$ is  $f$-invariant.  Thus the strata, if any, that are contained in wandering components of $G_i$ can be moved up the filtration to be above the strata contained in the non-wandering components of $G_i$.    After
finitely many such changes  and after amalgamating zero strata if necessary, we may assume that 

\begin{enumeratecontinue}
\item   If $G_i$ has  wandering components then $H_i$ is a wandering component of $G_i$.
\end{enumeratecontinue}

Suppose that $K$ is a component of the union  of all zero strata,  that $H_i$ is the highest stratum that contains an edge in  $K$ and that  $H_u$ is the highest irreducible stratum below $H_i$.  We prove that  $K \cap G_u = \emptyset$ by assuming that    $K \cap G_u \ne \emptyset$ and arguing to a contradiction.  Since each component of $G_u$ is non-wandering and since  some iterate of $f$ maps $K$ into $G_u$, there is a unique component $C$ of $G_u$ that intersects $K$.    If each  vertex $v \in K$ has valence at least two in $C \cup K$ then each edge of $K$ is    contained in a path in $K$ with  endpoints in $C$, and so by the connectivity of   $C$,  is contained in a circuit in  $K \cup C$.  This contradicts the fact that  some iterate of $f$ maps $K \cup C$ into $C$, and we conclude that   some vertex $v$ of $K$ has valence one in $K \cup C$.  In particular,  $v$  is non-periodic and so is not the endpoint of an \noneg\ edge by (1) and is not the endpoint of an \eg\ edge in a stratum above $G_i$ by  Remark~\ref{rtt2}.   By construction,   $v$ is  not the endpoint of an edge in a zero stratum above $G_i$.  But then $v$ has valence one in $G$.      This contradiction proves that $K \cap G_u = \emptyset$.    After reorganizing the edges in zero strata, we may assume that $H_i = K$.  Note in particular that no vertex in $H_i$ is periodic.

Let $H_r$ be the first irreducible stratum above $H_i$.   The component of $G_r$ that contains $H_i$ is non-wandering by (2) and so must intersect $H_r$.   Since no vertex of $H_i$ is periodic, (1) implies that $H_r$ is \eg.  Moreover,  the argument used in the previous paragraph proves that the link in $G$ of each vertex in $H_i$ is contained in $H_i \cup H_r \subset G_r$ and that $H_i$ is contained in the core of $G_r$.

 We arrange that each vertex in $H_i$ is contained in $H_r$ by the following {\em tree replacement} move.  Replace  $H_i$ with a tree $H_i'$ whose vertex set is exactly $H_i \cap H_r$.  Do this for each zero stratum $H_i$ and call the resulting graph $G'$.    We view the  union $X$ of the irreducible strata as a subgraph of both $G$ and  $G'$.  The set of vertices in $X$ is $f$-invariant by  (RTT-$i$) and (1).   There is a homotopy equivalence $p':G' \to G$ that is the identity on $X$ and that sends an edge $E'$ of $H'_i$ to the unique path in $H_i$ connecting the endpoints of $E'$.  The homotopy inverse $p:G \to  G'$ is the identity on $X$ and sends $H_i$ to $H_i'$; the exact definition depends on choosing $p(w) \in H_i'$ for those vertices $w \in H_i$, if any, that are not contained in $X$.    Note that $p_\#$ defines a bijection, with inverse $p'_\#$, between paths in $G$ with endpoints at vertices in $X$ and  paths in $G'$ with endpoints at vertices in $X$.  The homotopy equivalence  $f' :G' \to G'$ defined on edges by $f'(E') =  (pfp')_\#(E')$ is independent of the choices made in defining $p$ and represents $\phi$.    To verify (RTT-$ii$) for $f'$ it suffices to show that if $\alpha' \subset H_i'$ is non-trivial then     To verify (RTT-$ii$) for $f'$ it suffices to show that if $\alpha' \subset H_i'$ is non-trivial then      $f'_\#(\alpha') =  p_\#f_\#(p'_\#(\alpha'))$ is non-trivial.    This follows from (RTT-$ii$) for $f$  and the fact that both $p_\#$ and $p'_\#$ preserve non-triviality for paths with endpoints in $X$.     It is easy to see that (RTT-$i$) for $f$ implies (RTT-$i$) for $f'$ and that $\PF(f) = \PF(f')$.  Remark~\ref{rtt2} implies that $f' :G' \to G'$ is a \rtt.   None of the moves in this step change the free factor systems represented by filtration elements so $f' :G' \to G'$ still realizes $\cal C$.

\vspace{.1in} \noindent{\bf(Step 3 : Property (P$_{\cal C}$))} If (P$_{\cal C}$) fails then there is a stratum $H_m \subset \Per(f)$ that is a forest with the property that for each $\F^j \in \cal C$  there is  a filtration element  $G_{l(j)}$ such that $G_ {l(j)}\cup H_m$ realizes $\F^j$.      We will  construct a new \rtt\ $f' :G' \to G'$ with one fewer \noneg\ stratum that still realizes $\cal C$ and satisfies  (1)  and the weak form of (Z).   After repeating this finitely many times we will have achieved (P$_{\cal C}$).

Let $Y$ be  
 the set of edges   in $G \setminus H_m$ that are mapped entirely into $H_m$ by some iterate of $f$.    Then each edge of $Y$ is contained in a zero stratum and $H_m \cup Y$
is a forest that is mapped into itself by $f$ and into $H_m$ by some iterate of $f$.  We next arrange that
\begin{itemize}
\item[($\ast$)]     if $\alpha$ is a path in a zero stratum with endpoints at vertices and if $\alpha$ is not contained in $Y$ then $f_\#(\alpha)$ is not contained in  $ H_m\cup Y$. 
\end{itemize}
Suppose to the contrary that  $\alpha \subset H_k$ violates ($\ast$).  Choose  an edge $E_i \subset H_k$ that is crossed by $\alpha$ and  is not contained in $Y$.    Perform a tree replacement move on $H_k$ as in step 2, replacing $E_i$ by an edge connecting the endpoints of $\alpha$.  The new edge is mapped entirely into $H_m$ by some iterate of $f$ and we add it to $Y$.  After finitely many such moves $(\ast)$ is satisfied.

Let $G'$ be the marked graph obtained by collapsing each component of
$H_m \cup Y$ to a point and let $p : G \to G'$ be the corresponding quotient
map.    Identify the edges of $G'$ with those of $G \setminus (H_m \cup Y)$ and
define $f': G' \to G'$ by $f'(E) = [pf(E)]$. As an edge path, $f'(E)$ is obtained from $f(E)$
by removing all occurrences of edges in $H_m \cup Y$.   It follows that the strata $H_r$ and $p(H_r)$ (if the latter is non-empty) have the same type (zero, \eg, \noneg),   that  $f':G' \to G'$ has one fewer \noneg\ stratum than $\fG$,  that  $f': G' \to G'$  satisfies (1) and the weak form of (Z),  that $\PF(f') = \PF(f)$ and that  $f': G' \to G'$ satisfies (RTT-$i$).   Lemma~5.9 of \cite{bh:tracks} implies that $f': G' \to G'$  satisfies (RTT-$iii$).   

   To verify  (RTT-$ii$), suppose that $H'_s = p(H_r)$  is \eg\  and that $\alpha' \subset G'_{s-1}$ is a connecting path  for  $H'_s$.   If    $\alpha'$ is contained in a zero stratum then   $f'_\#(\alpha')$ is non-trivial  by $(\ast)$.      We may therefore assume     that   the  component $C'$ of $G_{s-1}'$ that contains $\alpha'$  is not a zero stratum and hence is non-wandering.   To prove that  $f'_\#(\alpha')$ is non-trivial,  it suffices, by Remark~\ref{rtt2}, to show that  each $v' \in H'_{s} \cap C'$ is a periodic point.  
   
   Since $v'$ is incident to an edge in $H'_s$, there is a vertex $v \in H_r$ such that $p(v) = v'$.  If $v$ is periodic, we are done.  We may therefore assume that $v \not \in H_m$.  If $v \in Y$ then the component of $G_{r-1}$ that contains $v$ is a zero stratum by the weak form of (Z) contradicting the assumption that $v'$ is contained in a non-wandering component of $G'_{s-1}$.  It follows that $v = p^{-1}(v')$.  By the same reasoning, $v$ is contained in  a non-wandering component of $G_{r-1}$ and so is periodic by Remark~\ref{rtt2}.  Thus $v'$ is periodic and we have verified (RTT-$ii$) for $f'$. This completes the proof  that $f':G'\to G'$ is a \rtt. 
   
    There exists $k > 0$ so that each non-contractible component of $G_{l(j)} \cup H_m$ is $f^k$-invariant and so that $f^k$ induces a  rank preserving bijection between  the non-contractible components of  $ G_{l(j)} \cup H_m  \cup  Y$  and the  non-contractible components of  $ G_{l(j)} \cup H_m$.  Thus $G_{l(j)} \cup H_m \cup Y$ and hence $  p(G_{l(j)})$ realizes $\F^j$, proving that $f':G' \to G'$ and its filtration realize $\cal C$.

    \vspace{.1in} \noindent{\bf(Step 4 : Property (Z))}    The lowest stratum in a non-wandering component   $C$ of a filtration element is either \eg\ or periodic.  Lemma~\ref{two directions for EG} and  item (5) of Lemma~\ref{formerly remarks}  imply that  $C$ is not   contractible.   (Z) therefore follows from the weak form of (Z).     

\vspace{.1in} \noindent{\bf(Step 5 : Property (NEG))}   Suppose that $H_s$ is
a non-periodic NEG stratum with edges $\{E_1,\ldots,E_m\}$ satisfying
$f(E_i) = E_{i+1}u_i$ for paths $u_i \subset G_{s-1}$ where indices
are taken mod $m$.  The component $C_i$ of $G_{s-1}$ that contains the
terminal endpoint $v_i$ of $E_i$ does not wander.  The lowest stratum
$H_t$ in $C_i$ is either EG or periodic.    In the former case, every
vertex in $H_t$ has at least two gates in $H_t$ by Lemma~\ref{two
directions for EG} and so $H_t$ is its own core.  In the latter case, the same result follows from
(P$_{\cal C}$) and item (1) of Lemma~\ref{formerly remarks}.

Choose a path $\tau \subset G_{s-1}$   from $v_i$ to a periodic vertex $w_i$ in $H_t$
and slide to change the terminal endpoint of $E_i$ to $w_i$.  After
performing this sliding operation finitely many times,  working up through the filtration, (NEG) is satisfied.   The resulting homotopy equivalence is a  \rtt\  by Lemma~\ref{sliding prelim},     still realizes $\cal C$ and still  satisfies (Z).     

Sliding may have introduced valence one vertices to $G$.   But no such vertex is the image of a vertex with valence greater than one by   (NEG), (Z) and  Lemma~\ref{two directions for EG}.   We may therefore remove all vertices of valence one and the edges that are incident to them.  After repeating this finitely many times $G$ has no valence one vertices.

  If (P$_{\cal C}$) is no longer satisfied then return to step 3.  Since this reduces the number of \noneg\ stratum the process stops.

\vspace{.1in} \noindent{\bf(Step 6 : Property (F))}  If $H_l$ is a zero stratum then $G_l$ and $G_{l-1}$ realize the same free factor system.  We may therefore assume that $G_l$ is not a zero stratum and hence that every component of $G_l$ is non-contractible.   If $w$ is a valence one vertex of $G_l$ then by item (1) of Lemma~\ref{formerly remarks},  (\noneg) and  Lemma~\ref{two
directions for EG},  $w$ must be the initial  endpoint of a non-periodic \noneg\ edge in some $H_k$ with $k \le l$ and no vertex with valence at least  two in $G_l$ maps to $w$.   The initial  endpoint of  each edge  in $H_k$ has valence one in $G_l$ and $G_l \setminus H_k$ is $f$-invariant.   We may therefore reorder   the strata to move $H_k$ above $G_l \setminus H_k$.  After finitely many such moves, $G_l$ is a core graph.  Working our way up the filtration  we arrange that (F) is satisfied.  
\qed

\vspace{.1in}

We conclude this section by recalling an operation from page 46 of \cite{bh:tracks} and  Definitions~5.3.2 of \cite{bfh:tits1}.  

 Suppose that  $H_r$ is an \eg\ stratum of a  \rtt\  $f:G \to G$  that satisfies item (Z) of Theorem~\ref{rtt existence}    and that   $\rho$ is an \iNp\ of height $r$.   Decompose
$\rho = \alpha \beta$ into a concatenation of 
maximal $r$-legal subpaths as in Lemma~\ref{nielsen paths in egs} 
and let $E_1 \subset H_r$ and $E_2
\subset H_r$ be the initial  edges of $\bar
\alpha$ and $\beta$ respectively.     If  one of the
edge paths $f(E_i)$, i=1 or 2, is an initial subpath of
the other then we say that {\em the fold at the illegal turn of $\rho$ is a  full fold}; otherwise it is a {\em  partial fold}.
There are two kinds of full folds.      If $f(E_1) \ne  f(E_2)$ then the full fold is  {\em proper}; otherwise it is {\em improper}.  

  Suppose that   the fold at the  illegal turn of $\rho$ is proper, say that     $f(E_1)$ is  a proper initial subpath of $f(E_2)$.   Write $\bar \alpha = E_1 b E_3\ldots $ where $b$ is  a (possibly trivial)   subpath of $G_{r-1}$ and $E_3$ is an edge in $H_r$.   Since both the initial  edge of $f(E_3)$ and the first edge of $f(\beta)$ that is not cancelled when $f(\alpha) f(\beta)$ is tightened to   $\alpha \beta$ belong to $H_r$, there is a decomposition $E_2 = E_2''E_2'$ into subpaths such that  $f(E_2'') = f_\#(E_1b)$ and such that the first edge in $f(E_2')$ is contained in $H_r$.  Subdivide $E_2$ into two edges $E_2''$ and  $E_2'$ and then identify $E_2''$ with $E_1b$ to form a new graph $G'$.   The quotient map $F: G \to G'$  is called the {\it extended fold determined by  $\rho_r$}.       
 
   We think of $G \setminus E_2$ as a subgraph of $G'$ on which $F$ is the identity.  By construction $F(E_2) =  E_1bE_2'$.     The filtration on $G'$ is defined by   $H_i' = H_i$ for $i \ne r$ and $H_r' = (H_r \setminus E_2) \cup E_2'$.     There is a map $g : G' \to G$ such that $gF = f$.      We  refer to  $g : G' \to G$ as {\it map induced by the extended fold}.

The following lemma states that the map $f'   : G'\to G'$  obtained from $Fg: G'\to G'$ by  tightening the images of edges is a \rtt\ that satisfies item (Z) of Theorem~\ref{rtt existence}.       We say that  $f' : G'
\to G'$  is obtained from $\fG$ by {\it folding  $\rho_r$} and that $\rho_r' = F_\#(\rho_r)$ is {\em the \iNp\ determined by $\rho_r$}.     If the illegal turn of   $\rho'_r$ is   proper  then this process can be repeated.   This is referred to as {\em iteratively folding $\rho$}.

\begin{lemma}  Assuming notation as above, $f'   : G'\to G'$  is a  \rtt\  that satisfies item (Z) of Theorem~\ref{rtt existence}.
\end{lemma}

\proof   By construction, $f|'G_{r-1} = f|G_{r-1}$.    If $E$ is an edge in $H_r$ then $f(E)$ does not cross the illegal turn in $\rho_r$.  If $E \ne E_2$ then $f'(E)$ is obtained from $f(E)$ by replacing each occurrence of $E_2$ with $E_1bE_2''$.    Similarly, $f'(E_2')$ is obtained from $f(E_2')$ by  replacing each occurrence of $E_2$ with $E_1bE_2''$.  It follows that $H'_r$ satisfies (RTT-$i$) - (RTT-$iii$). 

  If $H_k$ is a zero stratum above $H_r$ then each edge $E_k$ in $H_k$ is a connecting path for some \eg\ stratum $H_s$ above $H_k$ by item (Z) of Theorem~\ref{rtt existence}.   Thus $f^i_\#(E_k)$ is  non-trivial for all $i \ge 0$.  Since $F$ does not identify points that are not identified by $f$,   and since $F|E_k$ is the  identity, $(Fg)_\#(E_k) =( Ff)_\#(E_k)$ is non-trivial.   This shows that no edges are collapsed when  $Fg$ is tightened to $f'$.   The same argument shows that if $\sigma \subset H_k$ is any path with endpoints at vertices then $f'_\#(\sigma)$ is non-trivial.
  
  If $H_l$ is \noneg\ then $H_l'$ is \noneg.

   Suppose that $E_m$ is an edge in an \eg\ stratum $H_m$ above $H_r$ and that $f(E_m) = \mu_1\nu_1 \mu_2 \ldots \nu_l \mu_{l+1}$ is the decomposition into subpaths  $\mu_j \subset H_m$ and subpaths $\nu_j \subset G_{m-1}$.   Then $f'(E_m) = (Fg)_\#(E_m) =(Ff)_\#(E_m) = \mu_1\nu'_1 \mu_2 \ldots \nu'_l \mu_{l+1}$ where  $\nu_j' = F_\#(\nu_j)$ is non-trivial because $f_\#(\nu_j)$ is non-trivial.   This proves that $H_m'$
 satisfies (RTT-$i$) and   (RTT-$iii$).    

To verify (RTT-$ii$) for $H_m$ suppose that $\sigma'$ is a connecting path for $H'_m$.  If $\sigma'$ is contained in a zero stratum $H_k'$ then it is disjoint from $G_r$ and so is  identified with a connecting path $\sigma \subset H_k$.  By our previous argument, $f'_\#(\sigma')$ is non-trivial.  If $\sigma'$ is contained in non-contractible component of $G'_{m-1}$ then there is a connecting path $\sigma$ for $H_m$ in a  non-contractible component of $G_{m-1}$ such that $F_\#(\sigma) = \sigma'$.  The endpoints of $\sigma$ are periodic for $f$ by  Remark~\ref{rtt2}.  It follows that the endpoints of $\sigma'$ are periodic for $f'$ and another application of  Remark~\ref{rtt2} proves that $H_m'$ satisfies (RTT-$ii$).  
This completes the proof that  $f'   : G'\to G'$  is a  \rtt.

  Item (Z) of Theorem~\ref{rtt existence} for $f'$ therefore follows from item (Z) of Theorem~\ref{rtt existence} for $f$.   \endproof

\section{Forward Rotationless Outer Automorphisms} \label{rotOutAut}
To avoid issues raised by finite order phenomenon, one often replaces $\oone
\in \Out(F_n)$ with an iterate $\oone^k$.  In this section we explain
how this can be done canonically by exhibiting the natural class of
outer automorphisms that require no iteration.  We also define
principal automorphisms in the context of $\Out(F_n)$.  These
automorphisms play a central role in both the definition of forward
rotationless outer automorphisms (Definition~\ref{rotOutDef}) and in
the formulation of the Recognition Theorem
(Theorem~\ref{t:recognition}).

In section~\ref{nielsen approach}, we recall how principal
automorphisms occur in the context of the mapping class group.
Examples and definitions for $\Out(F_n)$ are given in
section~\ref{principal def}.  An equivalent definition is then given
in terms of relative train track maps and the Nielsen classes of their
fixed points.  Finally, in section~\ref{rotationless properties} we
record some properties of forward rotationless outer automorphisms
that justify their name; for example, we show that a $\oone$-periodic
free factor is $\oone$-invariant.

\subsection{The Nielsen Approach to the Mapping Class Group} \label{nielsen 
approach} To provide historical context and motivation for our
techniques and results, we briefly recall Nielsen's point of view on
the mapping class group.  Further details and proofs can be found, for
example, in \cite{ht:surfaces}.

Let $S$ be a closed orientable surface of negative Euler
characteristic and let $h : S \to S$ be a homeomorphism representing
an element $\mu \in \mcg(S)$. A choice of complete hyperbolic
structure on $S$ identifies the universal cover $\ti S$ of $S$ with
the hyperbolic plane $\H$.  Using the Poincare disk model for $\H$,
there is an induced compactification of $\ti S$ by adding a
topological circle $S_{\infty}$.
 
To avoid cumbersome superscripts we use $g$ to denote a positive
iterate $g:=h^k$ of $h$. Any lift $\ti g : \ti S \to \ti S$ of $g$
extends to a homeomorphism of the compactification.  The restriction
of this extension to $S_{\infty}$, denoted $\hat g : S_{\infty}\to
S_{\infty}$, depends only on $k$, the isotopy class of $h$ and the
choice of lift. More precisely, $h$ induces an outer automorphism of
$\pi_1(S)$ and $\hat g = \widehat{\Phi_{\ti g}}$ where $\Phi_{\ti g}$
is the automorphism of $\pi_1(S)$ corresponding to $\ti g$ and
$\widehat{\Phi_{\ti g}}:S_{\infty} \to S_{\infty}$ is the
homeomorphism determined by the identification of $S_{\infty}$ with
$\partial \pi_1(S)$ \cite{floyd}. 

Denote the set of non-repelling fixed points of $\hat g$ by
$\Fix_N(\hat g)$.  If $\Fix_N(\hat g)$ contains at least three
non-repelling points then we say that $\tilde g$ is a {\em principal
lift of $g$} and that $\Phi_{\ti g}$ is a principal automorphism
representing $\mu^k$.  The sets $\Fix_N(\hat g)$ determined by the
principal lifts of iterates of $h$ are central to Nielsen's
investigations.
 
The mapping class $\mu$ determined by $h$ is rotationless as defined
in the Section~\ref{s:intro} if and only for all $k$, each principal
lift $\ti g$ of $g = h^k$ has the form $\ti h ^k$ where $\ti h$ is a
principal lift of $h$ and where $\Fix_N(\hat g) = \Fix_N(\hat h)$.
Thus from the point of view of principal lifts and their $\Fix_N$
sets, nothing changes if $\mu$ is replaced by an iterate.  For the
remainder of this discussion we assume that $\mu$ is rotationless and
that $k=1$.
 
The intersection $\ti \Delta(\ti g)$ of the convex hull of
$\Fix_N(\hat g)$ with $\ti S = \H$ is called the {\em principal
region} for $\ti g$ and its image in $S$ is denoted $\Delta(\ti g)$.
Thus $\ti g$ is principal if and only if $\Delta(\ti g)$ has non-empty
interior.

Assume that $\ti g$ is principal.  If no point in $\Fix_N(\hat g)$ is
isolated then $\Delta(\ti g)$ is a compact subsurface and there is a
homeomorphism $f : S \to S$ representing $\mu$ whose restriction to
$\Delta(\ti g)$ is the identity.  If $\Fix_N(\hat g)$ is finite, or
more generally, is finite up to the action of a single covering
translation that commutes with $\ti g$, then the interior of
$\Delta(\ti g)$ is a component of the complement in $S$ of one of the
pseudo-Anosov laminations $\Lambda$ associated to $\mu$.  The boundary
of $\Delta(\ti g)$ is a finite union of leaves of $\Lambda$ and
perhaps one reducing curve.  These are the only cases that occur if
there is non-trivial twisting along each reducing curve in the
Thurston normal form for $\mu$.  In the general case, $\Delta(\ti g)$
is a finite union of the two types.
 
A proof of (most of) the Thurston classification theorem from this
point of view is contained in \cite{ht:surfaces} and
\cite{miller:nielsen}.

\subsection{Principal Automorphisms}  \label{principal def} 
Suppose that $\fG$ is a relative train track map representing $\oone
\in \Out(F_n)$.  Recall from section~\ref{autos and lifts} that there
is a bijection between lifts $\ti f: \Gamma \to \Gamma$ to the
universal cover and automorphisms $\Phi \in \Aut(F_n)$ representing
$\phi$.  A fixed point

\begin{definition} \label{d:principal auto} For $\aone \in \Aut(F_n)$  representing $\oone$, let $\Fix_N(\hat \aone)\subset \Fix(\hat 
\aone)$ be the set of non-repelling fixed points of $\hat \aone$.  We
say that $\aone$ is a {\em principal automorphism} and write {\em
$\Phi \in \PA(\oone)$} if either of the following hold.
\begin{itemize}
\item 
$\Fix_N(\hat \aone)$ contains at least three points. 
\item $\Fix_N( \hat \aone)$ is a two
point set that is neither the set of endpoints of an axis $A_c$ nor
the set of endpoints of a lift $\ti \lambda$ of a generic leaf of an
element of $\L(\oone)$.
\end{itemize} 
The corresponding lift $\ti f: \Gamma \to \Gamma$ is a {\em principal lift}.
\end{definition}

\begin{remark} \label{at least one}
For all $\oone \in \Out(F_n)$ there exists, by Lemma 5.2 of
\cite{bfh:tits3} or Proposition I.5 of \cite{ll:ends}, $k \ge 1$ such
that $\PA(\phi^k) \ne \emptyset$.  Moreover, if the conjugacy class of
$a \in F_n$ is invariant under $\phi^k$, then one may
choose $\aone \in \PA(\phi^k)$ to fix $a$.
\end{remark}

\begin{remark}  \label{rank two fixed subgroup} If $\Fix(\aone)$ has rank at 
least two then $\aone$ is a principal automorphism by Lemma~\ref{l:
second from bk3}.
\end{remark}

\begin{remark}  \label{distinct fixed point sets}
If $\aone_1$ and $\aone_2$ are distinct representatives of $\oone$
then $\Fix_N(\hat \aone_1) \cap \Fix_N(\hat \aone_2)$ is contained in
$\Fix(\hat \aone_1^{-1}\hat \aone_2) = \{T_c^{\pm}\}$ for some
non-trivial covering translation $T_c$.  It follows that if $\aone_1$
and $\aone_2$ are principal then $\Fix_N(\hat \aone_1) \ne \Fix_N(\hat
\aone_2)$.
\end{remark}

\begin{remark}
The second item in our definition of principal automorphism does not
occur in the context of mapping class groups.  It arises in
$\Out(F_n)$ to account for nonlinear \noneg\ strata as illustrated by
Example~\ref{example1} below.
\end{remark}

\begin{remark} \label{leaf lifts}
Each $\Lambda \in \L(\oone)$ has infinitely many generic leaves that
are invariant by an iterate of $\oone_\#$.  If $\{P,Q\}$ is the
endpoint set of a lift of such a leaf then (Lemma~\ref{isolated for
egs}) there exists $\aone$ representing an iterate of $\oone$ such
that $P$ and $Q$ are attracting fixed points for $\hat \aone$ .
Remark~\ref{finitely many isogredience classes} and Lemma~\ref{l:
second from bk3} imply that for all but finitely many such leaves,
$\Fix_N(\hat \aone) = \{P,Q\} $ and $\aone$ is not principal.
\end{remark}

\begin{remark} \label{gjll} 
If $\Phi$ has positive index in the sense of \cite{gjll}, then $\Phi$
is a principal automorphism.  The converse fails for the principal
automorphism $\Phi_2$ of Example~\ref{example1}.
\end{remark}

We say that $x, y \in \Fix(f)$ are {\em Nielsen equivalent} or belong
to the same {\em Nielsen class} if they are the endpoints of a Nielsen
path for $f$.  Each Nielsen class is an open subset of $\Fix(f)$
because every sufficiently short path with endpoints in $\Fix(f)$ is a
Nielsen path.  In particular, there are only finitely many Nielsen
classes.

If $\ti f: \Gamma \to \Gamma$ is a lift of $f:G \to G$, then any path
$\ti \alpha \subset \Gamma$ with endpoints in $\Fix(\ti f)$ projects
to a Nielsen path $\alpha \subset G$ for $f$.  Conversely, if $\alpha$
is a Nielsen path for $f$ and $\ti f$ fixes one endpoint of a lift
$\ti \alpha$ of $\alpha$ then $\ti f$ also fixes the other endpoint of
$\ti \alpha$.  Thus $\Fix(\ti f)$ is either empty or projects onto a
single Nielsen class in $\Fix(f)$.

A pair of automorphisms $\Phi_1$ and $\Phi_2$ {\em are equivalent} if
there exists $c \in F_n$ such that $\Phi_2 = i_c \Phi_1 i_c^{-1}$.
Translating this into the language of lifts, $\ti f_1$ is equivalent
to $\ti f_2$ if $\ti f_2 = T_c \ti f_1 T_c^{-1}$.  This equivalence
relation is called {\em isogredience}.

\begin{lemma}  \label{equivalence classes}
Suppose that $\fG$ represents $\oone \in \Out(F_n)$ and that $\ti f_1$
and $\ti f_2$ are lifts of $f$ with non-empty fixed point sets.  Then
$\ti f_1$ and $\ti f_2$ belong to the same isogredience class if and
only if $\Fix(\ti f_1)$ and $\Fix(\ti f_2)$ project to the same
Nielsen class in $\Fix(f)$.
\end{lemma}

\proof If $\ti f_2 = T_c \ti f_1 T_c^{-1}$ then $\Fix(\ti f_2) = T_c
\Fix(\ti f_1)$ and $\Fix(\ti f_2)$ and $\Fix(\ti f_1)$ project to the
same Nielsen class in $\Fix(f)$.  Conversely, if $\Fix(\ti f_2)$ and
$\Fix(\ti f_1)$ have the same non-trivial projection then there exists
$\ti x \in \Fix(\ti f_2)$ and a covering translation $T_c$ such that
$T_c(\ti x) \in \Fix(\ti f_1)$ which implies that $\ti f_2$ and $T_c
\ti f_1 T_c^{-1}$ agree on a point and hence are equal.  \endproof

\begin{remark} \label{finitely many isogredience classes}
We show below (Corollary~\ref{non-trivial fixed set}) that principal
lifts have non-trivial fixed point sets in $\Gamma$.  Since there are
only finitely many Nielsen classes in $\Fix(f)$, it follows that there
are only finitely many isogredience classes of principal lifts for
$\phi$.
\end{remark} 

In the following examples, $G$ is the rose $R_3$ with basepoint $v$ at
the unique vertex.  We use $A,B$ and $C$ to denote both the oriented
edges of $G$ and the corresponding generators of $F_3$.  Our examples
are all positive automorphisms $\aone$, meaning that they are defined
by $A \mapsto w_A$, $B \mapsto w_B$ and $C \mapsto w_C$ where
$w_A,w_B$ and $w_C$ are words in the letters $A,B$ and $C$ (and not
the inverses $\bar A, \bar B$ and $\bar C$).  These words also define a homotopy
equivalence $\fG$.  Since $w_A,w_B$ and $w_C$ use only $A,B$ and $C$, and not  $\bar A, \bar B$ and $\bar C$,
the homotopy equivalence is a \rtt\ for $\oone$.

The universal cover of $\ti G$ is denoted $\Gamma$ and we assume that
a basepoint $\ti v$ has been chosen.  Some statements in the examples
are left for the reader to verify or follow from results we establish
later in this section; none of these statements are ever quoted.

\begin{ex} \label{example1}
Let $\aone_1 \in \Aut(F_3)$ be determined by $w_A = A, w_B = BA$ and
$w_C =BCB^2$.  Then $\Fix(\aone_1) = \langle A, BA\bar B\rangle$ and
$\Fix(\hat \aone_1) = \partial(\Fix(\aone_1))$.    The lift $\ti f_1$ that fixes
$\ti v$ is the principal lift corresponding to $\aone_1$.

The unique fixed point $x$ of $f$ in the interior of $C$ is not
Nielsen equivalent to $v$.  Let $\ti C$ be the lift of $C$ whose
initial endpoint is $\ti v$ and let $\ti f_2$ be the lift that fixes
the unique lift of $\ti x$ of $x$ in $\ti C$.  Then $\ti f_2$ is
principal and $\Fix_N(\hat f_2)$ is a pair of attractors which bound
the line that is the union of the increasing sequence $\ti C \subset
(\ti f_2)_\#(\ti C) \subset (\ti f_2)^2_\#(\ti C)\subset \ldots$.  If
$\aone_2$ is the principal automorphism corresponding to $\ti f_2$
then $\Phi_2 = i_B^{-1} \Phi_1$.  \end{ex}

\begin{ex} \label{example2}
Let $\aone_1 \in \Aut(F_3)$ be determined b $w_A = A, w_B = BA$ and
$w_C =CB^2$.  Then $\Fix(\aone) = \langle A, BA\bar B\rangle$ and
$\Fix(\hat \aone)$ is the union of $\partial(\Fix(\aone))$ with the
$\Fix(\aone)$-orbit of a single attractor $P$. The lift $\ti f$ that
fixes $\ti v$ is the principal lift corresponding to $\aone$ and $P$
is the endpoint of the ray that is the union of the increasing
sequence $\ti C \subset (\ti f_2)_\#(\ti C) \subset (\ti f_2)^2_\#(\ti
C)\subset \ldots$.
\end{ex}

\begin{ex} \label{example3}
Let $\aone \in \Aut(F_3)$ be determined by $w_A = ACBA, w_B = BA$ and
$w_C =CBA$, let $\ti f$ be the lift that fixes $\ti v$ and let $\ti A,
\ti B, \ti C$ and $\ti A^{-1}$ be the lifts of the oriented edges $A,
B, C$ and $\bar A$ with $\ti v$ as initial vertex.  The directions
determined by the initial edges of $\ti A,\ti B, \ti C$ and $\ti
A^{-1}$ are fixed by $D \ti f$.  Lemma~\ref{lam2} produces attractors
$P_A,P_B,P_C$ and $P_{\bar A}$ in $\Fix_N(\hat f)$ such that lines
connecting $P_{\bar A}$ to the other three points are generic leaves
of an attracting lamination.  Lemma~\ref{just one ray} implies that
$\Fix_N(\hat f) =\{P_A,P_B,P_C,P_{\bar A}\}$.
\end{ex}

We now come to the second main definition of this section. Note that
if $\aone$ is a principal lift of $\oone$ then $\aone^k$ is a
principal lift for $\oone^k$ and $\Fix_N(\hat \aone) \subset
\Fix_N(\hat \aone^k)$ for all $k \ge 1$.  The set of non-repelling
periodic points in $\Per(\hat \Phi)$ is denoted $\Per_N(\hat \Phi)$.
By iterating $\oone$ we might pick up more principal lifts and
principal lifts might pick up more non-repelling fixed points. If this
doesn't happen, then we say that $\oone$ is {\it forward
rotationless}. Here is the precise definition.

\begin{definition} \label{rotOutDef}
An outer automorphism $\oone$ is {\em forward rotationless} if
$\Fix_N(\hat \aone) = \Per_N(\hat \aone)$ for all $\aone \in
\PA(\oone)$ and if for each $k \ge 1$, $\Phi \mapsto \Phi^k$ defines a
bijection (see Remark~\ref{surjection}) between $\PA(\oone)$ and
$\PA(\oone^k)$.  Our standing assumption is that $n \ge 2$. For
notational convenience we also say that the identity element of
$\Out(F_1)$ is forward rotationless.
\end{definition} 

\begin{remark} \label{surjection}
By Remark~\ref{distinct fixed point sets} there is no loss in
replacing the assumption that $\Phi \mapsto \Phi^k$ defines a
bijection with the a priori weaker assumption that $\Phi \mapsto
\Phi^k$ defines a surjection.
\end{remark}

\subsection{Rotationless Relative Train Track Maps and Principal Periodic Points}
We now want to characterize those \rtt s $\fG$ that represent forward
rotationless $\phi \in \Out(F_n)$ and to determine which lifts of such
$f$ are principal.  We precede our main definitions by showing that
principal lifts have fixed points.

Suppose that $\fG$ represents $\oone$ and that $\ti f: \Gamma \to
\Gamma$ is a lift of $f$.  We say that $\ti z \in \Gamma$ {\em moves
toward} $P \in \Fix(\hat f)$ under the action of $\ti f$ if the ray
from $\ti f(\ti z)$ to $P$ does not contain $\ti z$.  Similarly, we
say that {\it $\ti f$ moves $\ti y_1$ and $\ti y_2$ away from each
other} if the path in $\Gamma$ connecting $\ti f(\ti y_1)$ to $\ti
f(\ti y_2)$ contains $\ti y_1$ and $\ti y_2$ and if $\ti f(\ti y_1)
<\ti y_1 < \ti y_2 < \ti f(\ti y_2)$ in the order induced by the
orientation on that path.

The following lemma relates the action of $\hat f$ to the action of
$\ti f$ and gives a criterion for elements of $\Fix(\hat f)$ to be
contained in $\Fix_N(\hat f)$.  Recall that $\partial F_n$ is
identified with the set of ends of $\Gamma$.  It therefore makes sense
to say that points in $\Gamma$ are close to $P \in \partial F_n$ or
that $P$ is the limit of points in $\Gamma$.

\begin{lemma} \label{not a repeller}
Suppose that $P \in \Fix(\hat f)$ and that there does not exist $c \in
F_n$ such that $\Fix(\hat f) = \{T_c^{\pm}\}$.
\begin{enumerate}
\item If $P$ is an attractor for the action of $\hat f$ on
$\partial\Gamma$ then $\ti z$ moves toward $P$ under the action of
$\ti f$ for all $\ti z \in \Gamma$ that are sufficiently close to $P$.
\item If $P$ is an endpoint of an axis $A_c$ or if $P$ is the limit of
points in $\Gamma$ that are either fixed by $\ti f$ or that move
toward $P$ under the action of $\ti f$, then $P \in \Fix_N(\hat f)$.
\end{enumerate} 
\end{lemma}

\proof The lemma is an immediate consequence of Proposition 1.1 of
\cite{gjll} if $P$ is not the endpoint $T_c^{\pm}$ of an axis $A_c$.
If $P$ is $T^+_c$ or $T^- _c$, then $\Fix(\hat f)$ contains
$\{T_c^{\pm}\}$ and at least one other point.  Lemma~\ref{l: second
from bk3} implies that $P$ is not isolated in $\Fix(\hat f)$ and is
therefore neither an attractor nor a repeller for the action of $\hat
f$.  \endproof

The next lemma is based on Lemma 2.1 of \cite{bh:tracks}.

\begin{lemma} \label{l:maps over}
If $\ti f$ moves $\ti y_1$ and $ \ti y_2$ away from each other, then
$\ti f$ fixes a point in the interval bounded by $\ti y_1$ and $\ti
y_2$.
\end{lemma} 

\begin{proof}  
Denote the oriented paths connecting $\ti y_1$ to $\ti y_2$ and $\ti
f(\ti y_1)$ to $\ti f(\ti y_2)$ by $\ti \alpha_0$ and $\ti \alpha_1$
respectively.  Let $r : \Gamma \to \ti \alpha_1$ be retraction onto
the nearest point in $\ti \alpha_1$ and let $g = r \ti f : \ti
\alpha_0 \to \ti \alpha_1$.  By hypothesis, $\ti \alpha_0$ is a proper
subpath of $\ti \alpha_1$ and $g$ is a surjection. If $\ti y$ is the
first point in $\ti \alpha_0$ such that $\ti g(\ti y) = \ti y$ then
$\ti y_1 < \ti y < \ti y_2 $ and $\ti g(\ti z) < \ti g(\ti y)$ for
$y_1 < z < y$.  It follows that $\ti f(\ti y) \in \ti \alpha_1$ and
hence that $\ti y$ is fixed by $\ti f$.
\end{proof} 

\begin{corollary} \label{non-trivial fixed set}
If $\ti f$ is a principal lift then $\Fix(\ti f) \ne \emptyset$.
\end{corollary}

\proof Suppose that there is a non-trivial covering translation $T_c$
that has its endpoints in $\Fix_N(\hat f)$ and so commutes with $\ti
f$.  Assuming without loss that $A_c$ is fixed point free, there is a
point in $A_c$ that moves toward one of the endpoints of $A_c$, say
$P$.  Since $\ti f$ commutes with $T_c$, there are points in $\Gamma$
that are arbitrarily close to $P$ and that move toward $P$.  The same
property holds for an attractor $P \in \Fix(\hat f)$ by Lemma~\ref{not
a repeller}.  One may therefore choose distinct $P_1,P_2 \in
\Fix_N(\hat f)$ and $\ti x_1,\ti x_2 \in \ti \Gamma$ such that $\ti
x_i$ is close to and moves toward $P_i$.  It follows that $\ti x_1$
and $\ti x_2$ move away from each other.  Lemma~\ref{l:maps over}
produces the desired fixed point.  \endproof

There are two cases in which a lift $\ti f$ corresponding to a Nielsen
class in $\Fix(f)$ is not principal.  The first arises from a \lq
non-singular\rq\ leaf of an attracting lamination as noted in
Remark~\ref{leaf lifts}; in this case $\Fix(\ti f)$ is a single point.
In the second case, there is a circle component of $\Fix(f)$ with no
outward pointing periodic directions and $\Fix(\ti f)$ is an axis
$A_c$.  The second type could be eliminated by adding properties to
Theorem~\ref{rtt existence}.  We allow the circle components for now
and defer the additional properties until section~\ref{section cs
exists}.

Periodic points for $f$ are {\em Nielsen equivalent} if they are
Nielsen equivalent as fixed points for some iterate of $f$.

\begin{definition}\label{rotRttDef}
 We say that $x$ is {\em principal} if
neither of the following conditions are satisfied.
\begin{itemize}
\item $x$ is the only element of $\Per(f)$ in its Nielsen class and
there are exactly two periodic directions at $x$, both of which are
contained in the same EG stratum.
\item $x$ is contained in a component $C$ of $\Per(f)$ that is
topologically a circle and each point in $C$ has exactly two periodic
directions.
\end{itemize}
Lifts to $\Gamma$ of principal periodic points in $G$ are said to be
{\em principal}.  If each principal vertex and each periodic direction
at a principal vertex has period one then we say that $\fG$ is {\em
rotationless}.
\end{definition}

In practice, we only apply these definitions to $\fG$  that satisfy the conclusions of Theorem~\ref{rtt
existence}.  In particular,  by  item (4) of Lemma~\ref{formerly remarks}, there are at least two periodic directions at each $x \in \Per(f)$. 

Principal periodic points are either contained in periodic edges or
are vertices.  Thus every $\fG$ has a rotationless iterate.  Any
endpoint of an \ipNp\ is principal as is the initial endpoint of any
non-periodic NEG edge.  The latter implies that each NEG stratum in a
rotationless \rtt\ is a single edge.

The following lemma shows that an $EG$ stratum has at least one
principal vertex.

\begin{lemma} \label{some essential}
Assume that $\fG$ satisfies the conclusions of Theorem~\ref{rtt
existence}.  For every \eg\ stratum $H_r$ there is a principal vertex
whose link contains a periodic direction in $H_r$.
\end{lemma} 
\proof If some vertex $v \in H_r$ belongs to a non-contractible
component of $G_{r-1}$ then $v$ is periodic and there is at least one
periodic direction in $H_r$ by Lemma~\ref{r legal}.  There is also at
least one periodic direction at $v$ determined by an edge of $G_{r-1}$
so $v$ is principal.  If there is no such vertex, then $H_r$ is a union of 
components of $G_r$.  Lemma~5.2 of \cite{bfh:tits3} states there is a
principal lift of some iterate of $f|H_r$.  In the course of proving
this lemma, it is shown that either there is a vertex with three
periodic directions or there is an indivisible periodic Nielsen path
in $H_r$. In either case there is a vertex that is principal for
$f|H_r$ and hence also for $f$.  \endproof

\begin{remark}
Lemma~\ref{some essential} implies that the transition matrix $M_r$ of
an EG stratum of a rotationless $\fG$   satisfying the conclusions of Theorem~\ref{rtt
existence} has at least one non-zero
diagonal entry and so is aperiodic.  For each $\Lambda \in \L(\oone)$
there is an EG stratum $H_r$ such that $\Lambda$ has height $r$ and
this defines a bijection (see Definition 3.1.12 of \cite{bfh:tits1})
between $\L(\oone)$ and the set of EG strata.
\end{remark}    

The next lemma relates an attractor in $\Fix_N(\hat f)$ to a fixed
direction of $D \ti f$.

\begin{lemma} \label{just one ray}
Suppose that $\ti f$ is a principal lift of a \rtt\ $\fG$.
\begin{enumerate}
\item For each attractor $P \in \Fix_N(\hat f)$ there is a (not
necessarily unique) $\ti x \in \Fix(\ti f)$ such that the interior
of the ray $\ti R_{\ti x, P}$ that starts at $\ti x$ and that
converges to $P$ is fixed point free.
\item If $P \in \Fix_N(\hat f)$ is an attractor, if $\ti x \in \Fix(\ti
f)$ and if the interior of $\ti R_{\ti x,P}$ is fixed point free then
no point in the interior of $\ti R_{\ti x,P}$ is mapped by $\ti f$ to
$\ti x$; in particular, the initial direction determined by $\ti
R_{\ti x,P}$ is fixed.
\item If $P$ and $Q$ are distinct attractors in $\Fix_N(\hat f)$, if
$\ti x \in \Fix(\ti f)$ and if the interiors of both $\ti R_{\ti x,P}$
and $\ti R_{\ti x,Q}$ are fixed point free then the directions
determined by $\ti R_{\ti x,P}$ and $\ti R_{\ti x,Q}$ are distinct.
\end{enumerate}
\end{lemma}

\proof To find $\ti x \in \Fix(\ti f)$ and $\ti R_{\ti x,P}$ as in
(1), start with any ray $\ti R'$ whose initial point is in $\Fix(\ti
f)$ and that converges to $P$ and let $\ti R_{\ti x,P}$ be the subray
of $\ti R'$ that begins at the last point $\tilde x$ of $\Fix(\ti f)$
in $\ti R'$.  If $\ti R_{\ti x,P}$ and $\ti R_{\ti x,Q}$ are as in (3)
and have the same initial edge then their \lq difference\rq\ would be
a fixed point free line whose ends converge to attractors in
contradiction to Lemma~\ref{not a repeller} and Lemma~\ref{l:maps
over}.  This verifies (3).  By the same reasoning, no points in the
interior of $\ti R_{\ti x,P}$ can map to $\ti x$, which implies that
the initial edge of $\ti R_{\ti x,P}$ determines a fixed direction at
$\ti x$.  This proves (2). \endproof

\begin{corollary}  \label{principal is essential} Assume that $\fG$ satisfies the conclusions of 
Theorem~\ref{rtt existence}.  If $\ti f$ is a principal lift then each
element of $\Fix(\ti f)$ is principal.
\end{corollary}   

\proof Let $\Phi$ be the automorphism corresponding to $\ti f$.  If
$\Fix(\Phi)$ has rank at least two then $\Fix(\ti f)$ is neither a
single point nor a single axis and we are done (see Corollary~3.17 and
Lemma~2.1).  If $\Fix(\Phi)$ has rank one then $\Fix(\ti f)$ is
infinite and $\Fix_N(\hat f)$ contains an attractor by the definition
of principal lift and by Lemma~\ref{l: second from
bk3}. Lemma~\ref{just one ray} implies that some $\ti x \in \Fix(\ti
f)$ has a fixed direction that does not come from a fixed edge and
again we are done.  In the remaining case, $\Fix(\hat f)$ is a finite
set of attractors and does not contain the endpoints of any
axis. Obviously $\Fix(\ti f)$ is not an axis.  Suppose that $\Fix(\ti
f)$ is a single point $\ti x$, that there are only two periodic
directions at $\ti x$ and these two directions are determined by lifts
$\ti E_1$ and $\ti E_2$ of oriented edges of the same EG stratum
$H_r$. Lemma~\ref{just one ray} and Lemma~\ref{lam2} imply that
$\Fix(\hat f)$ is the endpoint set of a generic leaf of an element of
$\L(\oone)$ in contradiction to the assumption that $\ti f$ is
principal.  We conclude that $\ti x$ is principal as desired.
\endproof

To prove the converse we use fixed directions of $\ti f$ to find
elements of $\Fix_N(\hat f)$.  The following lemma is from
\cite{bfh:tits1}; the proof is short and is repeated for the readers
convenience.

\begin{lemma}\label{sliding from book1}  If $\Fix(\ti f) = \emptyset$ then  
there is a ray $\ti R \subset \Gamma$ converging to an element $P \in
\Fix(\hat f)$ and there are points in $\ti R$ arbitrarily close to $P$
that move toward $P$.
\end{lemma}

\proof For each vertex $\ti y$ of $\Gamma$, we say that the initial
edge of the path from $\ti y$ to $h(\ti y)$ is {\em preferred} by $\ti
y$.  Starting with any vertex $\ti y_0$, inductively define $\ti
y_{i+1}$ to be the other endpoint of the edge preferred by $\ti y_i$.
If $\ti E$ is preferred by both of its endpoints then $\ti f$ maps a
proper subinterval of $\ti E$ over all of $\ti E$ (reversing
orientation) in contradiction to the assumption that $\Fix(\ti f) =
\emptyset$.  It follows that the $\ti y_i$'s are contained in a ray
that converges to some $P \in \Fix(\hat f)$ and that $\ti y_i$ moves
toward $P$.  \endproof

We isolate the following notation and lemma for reference throughout
the paper.

\begin{notn} \label{defining h} 
Suppose that $\fG$ satisfies the conclusions of Theorem~\ref{rtt
existence}, that $H_r$ is a single edge $E_r$ and that $f(E_r) = E_r u$
for some non-trivial path $u \subset G_{r-1}$.  Let $\ti E_r$ be a
lift of $E_r$ and let $\ti f : \Gamma \to \Gamma$ be the lift of $f$
that fixes the initial endpoint of $\ti E_r$.  By (NEG), the component $C$ of $G_{r-1}$ that contains the
terminal endpoint $w$ of $E_r$ is not contractible. Denote the copy of
the universal cover of $C$ that contains the terminal endpoint of $\ti
E_r$ by $\Gamma_{r-1}$ and the restriction $\ti f|\Gamma_{r-1}$ by $h
: \Gamma_{r-1} \to \Gamma_{r- 1}$.

The covering translations that preserve $\Gamma_{r-1}$ define a free
factor $F(C)$ of $F_n$ such that $[[F(C)]] = [[\pi_1(C)]]$.  The
closure in $\partial F_n$ of $\{T_c^{\pm} : c \in F(C)\}$ is naturally
identified with $\partial F(C)$ and with the space of ends of
$\Gamma_{r-1}$.  Moreover, $\hat h = \hat f|\partial F(C): \partial
F(C) \to \partial F(C)$.
\end{notn}

\begin{lemma}\label{sliding corollary}
Assume that $\ti f$ and $h$ are as in Notation~\ref{defining h}. If
$\Fix(h) = \emptyset$ then there is a ray $\ti R \subset \Gamma_{r-1}$
converging to an element $P \in \Fix(\hat h)$ and there are points in
$\ti R$ arbitrarily close to $P$ that move toward $P$.
\end{lemma}

\proof This follows from Lemma~\ref{sliding from book1} applied to $h
: \Gamma_{r-1} \to \Gamma_{r-1}$.  \endproof

Our next result is an extension of Lemma~\ref{lam2}.

\begin{lemma} \label{new iterates to}
Suppose that $\fG$  satisfies the conclusions of Theorem~\ref{rtt
existence} and is   rotationless, that $\ti f: \Gamma \to \Gamma$ is
a lift of $f$, that $\ti v \in \Fix(\ti f)$ and that $D\ti f$ fixes
the direction at $\ti v$ determined by a lift $\ti E$ of an edge $E
\subset H_r$.  Then there exists $P \in \Fix(\hat f)$ so that the ray
$\ti R$ from the initial endpoint of $\ti E$ to $P$ contains $\ti E$
and satisfies the following properties.
\begin{enumerate}
\item There are points in $\ti R$ arbitrarily close to $P$ that are
either fixed or move toward $P$.  If there does not exist $c \in F_n$
such that $\Fix(\hat f) = \{T_c^{\pm}\}$ then $P \in \Fix_N(\hat f)$.
\item If $H_r$ is EG then $P$ is an attractor whose accumulation set
is the unique attracting lamination of height $r$, the interior of
$\ti R$ is fixed point free and $\ti R$ projects to an $r$-legal ray
in $G_r$.
\item  If $H_r$ is  NEG and non-fixed  then $\ti R \setminus \ti E$ 
projects into $G_{r-1}$. 
\item No point in the interior of $\ti R$ is mapped to $\ti v$ by any
iterate of $\ti f$.
\end{enumerate}
\end{lemma} 

\proof The second part of (1) follows from the first part of (1) and
Lemma~\ref{not a repeller}.

The proof is by induction on $r$, starting with $r=1$.  If $G_1 \subset \Fix(f)$ then  we
may choose $P$ to be the endpoint of any ray $\ti R$ that begins with
$\ti E$ and projects into $G_1$; the existence of such a ray follows from  the fact (Theorem~\ref{rtt
existence}(F)) that $G_1$ is its own core.  If $G_1$ is EG then the existence of
$P$ follows from Lemma~\ref{lam2} and Lemma~\ref{just one ray}.  This
completes the $r=1$ case so we may now assume that the lemma holds for
edges with height less than $r$.

If $H_r$ is EG then the existence of $P$ follows from Lemma~\ref{lam2}
and Lemma~\ref{just one ray}.  We may therefore assume that $H_r$ is
NEG.  Let $h :\Gamma_{r-1} \to \Gamma_{r-1}$ be as in
Notation~\ref{defining h}.  If $\Fix(h) \ne \emptyset$, then the
initial endpoint of $\ti E$ and some $\ti x \in \Fix(h)$ cobound an
\iNp.  Thus $\ti x$ is principal, there is a fixed direction in
$\Gamma_{r-1}$ at $\ti x$ and the existence of an appropriate $P \in
\Fix(\hat h)$ follows from the inductive hypothesis.  The case that
$\Fix(h) = \emptyset$ follows from Lemma~\ref{sliding corollary}.
\endproof

We now can prove the converse to Corollary~\ref{principal is
essential} under the assumption that $\fG$ is rotationless.

\begin{corollary} \label{new essential is principal}
Suppose that $\fG$ satisfies the conclusions of Theorem~\ref{rtt
existence} and is   rotationless.  If some, and hence every, $\ti x
\in \Fix(\ti f)$ is principal then $\ti f$ is principal.
\end{corollary}

\proof   Assume that $\Fix(\ti f)$ consists of principal points.  Since
$\fG$ is rotationless, periodic directions based in $\Fix(\ti f)$ are
fixed.  Lemma~\ref{formerly remarks}(4) implies that each  $\ti x \in \Fix(\ti f)$ has at least two  fixed directions.   If some $\ti x \in \Fix(\ti f)$ has at least three fixed
directions, then Lemma~\ref{new iterates to} produces at least three
points in $\Fix_N(\hat f)$ and we are done.    We may therefore assume
that there are exactly two fixed, and hence exactly two periodic,
directions at each $\ti x \in \Fix(\ti f)$.    If $\Fix(\ti f)$ contains an edge, then by Definition~\ref{rotRttDef} there must be such an edge with a valence one vertex in $\Fix(f)$.   This contradicts  items (1) and (4) of Lemma~\ref{formerly remarks} and we conclude that 
that there are no fixed edges.    Choose an edge $E \subset
H_r$ and a lift $\ti E$ whose initial direction is fixed and based at
some $\ti x \in \Fix(\ti f)$.  Let $\ti R$ be the ray that begins with
$\ti E$ and ends at some $P \in \Fix(\hat f)$ as in Lemma~\ref{new
iterates to}.

If $H_r$ is EG then the accumulation set of $P$ is an attracting
lamination which implies by Lemma~3.1.16 of \cite{bfh:tits1} that $P$
is not the endpoint of an axis.  If $H_r$ is NEG then the accumulation
set of $P$ is contained in $G_{r-1}$ which implies that $P$ is not the
endpoint of an axis that contains $\ti E$.  It follows that the line
composed of $\ti R$ and the ray determined by the second fixed
direction at $\ti x$ is not an axis.  We have now shown that
$\Fix(\hat f)$ is not the endpoint set of an axis and hence that every
point in $\Fix(\hat f)$ produced by Lemma~\ref{new iterates to} is
contained in $\Fix_N(\hat f)$.  Thus $\Fix_N(\hat f)$ contains at
least two points and is not the endpoint set of an axis.

To complete the proof we assume that $\Fix_N(\hat f)$ is the endpoint
set of a lift $\ti \ell$ of a generic leaf of an attracting lamination
and argue to a contradiction.  Since $\ell$ is birecurrent and
contains $E$, $H_r$ is EG and the second fixed direction based at
$\ti x$ comes from an edge in $H_r$. Lemma~\ref{new iterates to}(2)
implies that $\ell \subset G_r$ is $r$-legal and hence does not
contain any \iNp s of height $r$.  But then $\ti x$ must be the only
fixed point in $\ti \ell$.  Since $\Fix(\ti f)$ is principal it
must contain a point other than $\ti x$ and that point would have a
fixed direction that does not come from the initial edge of a ray
converging to an endpoint of $\ti \ell$.  This contradiction completes
the proof.  \endproof

\subsection{Rotationless is  Rotationless}
We prove in this section that rotationless \rtt s represent forward
rotationless outer automorphisms and vice-versa.

\begin{lemma} \label{period one} Suppose that $\fG$  satisfies the conclusions of Theorem~\ref{rtt
existence} and is   rotationless.
   Every periodic Nielsen path
$\sigma$ with principal endpoints has period one.
\end{lemma}

\proof There is no loss in assuming that $\sigma$ is either a single
edge or an indivisible periodic Nielsen path.  In the former case,
$\sigma$ is a periodic edge with an principal endpoint and so is
fixed.  We may therefore assume that $\sigma$ is indivisible.

The proof is by induction on the height $r$ of $\sigma$ with the $r=0$
case being vacuously true.  Let $p$ be the period of $\sigma$ and let
$v \in \Fix(f)$ be an endpoint of $\sigma$.  The case that $H_r$ is
\EG\ follows from Lemma~\ref{nielsen paths in egs}(3).

We may therefore assume that $H_r$ is a single non-fixed NEG edge
$E_r$.  Lemma 4.1.4 of \cite{bfh:tits1} implies that after reversing
the orientation on $\sigma$ if necessary, $\sigma = E_r \mu$ or
$\sigma= E_r \mu \bar E_r$ for some path $\mu \subset G_{r-1}$.  Let
$\ti E_r$ be a lift of $E_r$ with initial endpoint $\ti v$, let $\ti
f$ be the lift that fixes $\ti v$ and let $h : \Gamma_{r-1}\to
\Gamma_{r- 1}$ be as in Notation~\ref{defining h}.  By Lemma~\ref{new
essential is principal}, $\ti f$ is principal. Denote the terminal
endpoint of the lift $\ti \sigma$ that begins at $\ti v$ by $\ti w$.

If $\sigma=E_r\mu$ then $\ti w \in \Gamma_{r-1}$.  If $p \ne 1$ then
the path $\ti \tau$ connecting $\ti w$ to $h(\ti w)$ projects to a
non-trivial periodic Nielsen path $\tau \subset G_{r- 1}$ that is
closed because $\ti w$ projects to $w \in \Fix(f)$.  Since $\ti w$ is
principal, the inductive hypothesis implies that $\tau$ has period one
and hence that the projection of the closed path $\ti \tau h(\ti \tau)
\dots h^{p-1}(\ti\tau)$ to $G_{r-1}$ is homotopic to $\tau^p$.  This
contradicts the fact that $\tau$ and hence $\tau^p$ determines a
non-trivial conjugacy class in $F_n$.  Thus $p = 1$ in the case that
$\sigma=E_r\mu$.

Suppose now that $\sigma= E_r \mu \bar E_r$. If $\Fix(h^p) \ne
\emptyset$ then the path $\ti \sigma_1$ connecting $\ti v$ to $\ti x
\in \Fix(h^p)$ and the path $\ti \sigma_2$ connecting $\ti x$ to $\ti
w$ are periodic Nielsen paths.  By the preceding case $\sigma_1$ and
$\sigma_2$, and hence $\sigma$, has period one.  We may therefore
assume that $\Fix(h^p) = \emptyset$.

Let $T_c : \Gamma \to \Gamma$ be the covering translation satisfying
$T_c(\ti v) = \ti w$.  Then $T_c$ commutes with $\ti f^p$, the axis
$A_c$ is contained in $\Gamma_{r-1}$ and $ T_c^{\pm} \in \Fix(\hat
h^p)$.  If $\Phi$ is the principal automorphism corresponding to $\ti
f$ then $\hat T_{\Phi(c)} = \hat f \hat T_c \hat f^{-1}$, which
implies that $T^{\pm}_{\Phi(c)} =\hat h(T^{\pm}_c) \in \Fix(\hat h^p)$
and that $A_{\Phi(c)} \subset \Gamma_{r-1}$.  If $\{T^{\pm}_c\}$ is
not $\hat h$-invariant, then $\Fix_N(\hat h^p)$ contains the four
points $\{T^{\pm}_c\} \cup \{\hat h(T^{\pm}_c)\}$ and $h^p$ is a
principal lift of $f|C$ where $C$ is the component of $G_{r-1}$ that
contains the terminal endpoint of $E_r$.  This contradicts
Corollary~\ref{non-trivial fixed set} and the assumption that $\Fix(h^p) =
\emptyset$. Thus $\{T^{\pm}_c\}$ is $\hat h$-invariant.  If $\hat h$
interchanges $T_c^{\pm}$ then $\Fix_N(\hat h^p)$
contains $\{T_c^{\pm}\}$ and at least one point in $\Fix(\hat h)$ by
Lemma~\ref{sliding corollary}.  This contradicts
Corollary~\ref{non-trivial fixed set} and we conclude that $T_c^{\pm} \in
\Fix(\hat h)$. It follows that $\ti f$ commutes with $T_c$ and hence
that $\ti w \in \Fix(\ti f)$.  This proves that $p = 1$ and so
completes the inductive step.  \endproof

\begin{proposition} \label{rotationless}
Suppose that $\fG$ represents $\oone$ and satisfies the conclusions of
Theorem~\ref{rtt existence}.  Then $\fG$ is rotationless if and only
$\oone$ is forward rotationless.
\end{proposition}

\begin{proof}
Suppose that $\fG$ is rotationless, that $k \ge 1$ and that $\ti g :
\Gamma \to \Gamma$ is a principal lift of $g:=f^k$.
Corollary~\ref{non-trivial fixed set} and Corollary~\ref{principal is
essential} imply that $\Fix(\ti g)$ is a non-empty set of principal
fixed points.  Since $f$ is rotationless, for each $\ti v \in \Fix(\ti
g)$ there is a lift $\ti f : \Gamma \to \Gamma$ that fixes $\ti v$ and
all periodic directions at $\ti v$.  To prove that $\oone$ is forward
rotationless it suffices to show that $\Fix_N(\hat f) = \Fix_N(\hat
g)$ and hence (Remark~\ref{distinct fixed point sets}) that $\ti f^k =
\ti g$.

The path connecting  $\ti v$ to another 
point in $\Fix(\ti g)$  projects to a Nielsen path  for $g$ and hence by 
Lemma~\ref{period one}, a Nielsen path for $f$.  Thus $\Fix(\ti f) = \Fix(\ti 
g)$.  It follows that $\ti g$ and $\ti f$ commute with the same covering 
translations and  Lemma~\ref{l: second from bk3} implies that $\Fix_N(\hat f)$ 
and $\Fix_N(\hat g)$ have the same non-isolated points.  

Each isolated point $P \in \Fix_N(\hat g)$ is an attractor for $\hat
g$.  It suffices to show that $P \in \Fix_N(\hat f)$.  By
Lemma~\ref{just one ray} there is a ray $\ti R$ that terminates at
$P$, that intersects $\Fix(\ti g)$ only in its initial endpoint and
whose initial direction is fixed by $D\ti g$, and hence also by $D\ti
f$.  We may assume that the height $r$ of the initial edge $\ti E$ of
$\ti R$ is minimal among all choices of $\ti R$.  By Lemma~\ref{new
iterates to}, $\ti E$ extends to a ray $\ti R'$ that converges to some
$P' \in \Fix(\hat f)$.  It suffices to show that $P' = P$ since a
repeller for $\hat f$ could not be an attractor for $\hat g$.  If
$H_r$ is EG this follows from Lemma~\ref{new iterates to}(2) and
Lemma~\ref{just one ray}(3) applied to $g$.  We may therefore assume
that $H_r$ is NEG.  If there exists $\ti x \in \Fix(\ti g) \cap \ti
R'$ then $\ti x \in \ti R'\setminus \ti E$ and the ray connecting $\ti
x$ to $P$ is contained in $G_{r-1}$ in contradiction to our choice of
$r$. We may therefore assume that $\Fix(\ti g) \cap \ti R' =
\emptyset$.  Lemma~\ref{new iterates to}(1) implies that there exists
$\ti x \in \ti R'$ that is moved toward $P'$ by $\ti f$ and
Lemma~\ref{l:maps over} then implies that $P = P'$.  This completes
the proof of the only if direction of the proposition.

For the if direction, assume that $\oone$ is forward rotationless and
choose $k > 0$ so that $g:=f^k$ is rotationless.  For each principal
$v \in \Fix(g)$, there exist a lift $\ti v$ of $v$ and a principal
lift $\ti g$ of $g$ that fixes $\ti v$.  Since $\oone$ is forward
rotationless, there is a lift $\ti f$ of $f$ such that $\ti f^k = \ti
g$ and such that $\Fix_N(\hat f) = \Fix_N(\hat g)$.  It suffices to
show that $\ti v \in \Fix(\ti f)$ and that each $D\ti g$-fixed
direction $\ti d_1$ at $\ti v$ is $D\ti f$-fixed.

The edge determined by $\ti d_1$ extends to a ray $\ti R_1$ that
converges to some $P_1 \in \Fix_N(\hat g) = \Fix_N(\hat f)$.  Define
$P_2$ and $\ti R_2$ similarly using a second $D\ti g$-fixed direction
$\ti d_2$ based at $\ti v$ and denote the line connecting $P_1$ to
$P_2$ by $\ti \gamma$.  Thus $\ti f_\#(\ti \gamma) =\ti \gamma$ and
the turn $(\ti d_1,\ti d_2)$ is legal for $\ti g$ and hence for $\ti
f$.  If $\ti f(\ti v) \not \in \ti \gamma$ then there exists $\ti y
\in \ti \gamma$ not equal to $\ti v$ such that $\ti f(\ti y) = \ti
f(\ti v)$.  But then $\ti f^k(\ti y) = \ti f^k(\ti v) = \ti v$ which
contradicts Lemma~\ref{new iterates to}(4) applied to $\ti g$.  This
proves that $\ti f(\ti v) \in \ti \gamma$.  Suppose that $\ti f(\ti v)
\ne \ti v$.  Denote $\ti v$ by $\ti v_0$ and orient $\ti \gamma$ so
that $P_1$ is the negative end and $P_2$ is the positive end.
Assuming without loss that $\ti v < \ti f(\ti v)$ in the order induced
from the orientation, there exist $\ti v_i \in \ti \gamma$ for $1 \le
i \le k$ such that $\ti v_i < \ti v_{i-1}$ and such that $\ti f(\ti
v_i) = \ti v_{i-1}$.  But then $\ti f^k(\ti v_k) = \ti v$ in
contradiction to Lemma~\ref{new iterates to}(4).  We conclude that
$\ti f(\ti v) \ne \ti v$.  A third application of Lemma~\ref{new
iterates to}(4) implies that the directions $\ti d_i$ are fixed by
$D\ti f$.
\end{proof}  

 \subsection{Properties of Forward Rotationless $\phi$} \label{rotationless properties}

\begin{lemma} \label{no iterates necessary}
The following   hold for each forward rotationless $\oone \in
\Out(F_n)$.
\begin{enumerate}
\item Each periodic conjugacy class is fixed and each representative
of that conjugacy class is fixed by some principal automorphism
representing $\phi$.
\item Each $\Lambda \in \L(\oone)$ is $\oone$-invariant.
\item  A free factor that is invariant under an iterate of $\oone$ is $\oone$-invariant. 
\end{enumerate}
\end{lemma}

\proof If the conjugacy class of $c$ is fixed by $\oone^k$ for $k \ge
1$ then by Remark~\ref{at least one} there exists a principal
automorphism $\aone_k \in \PA(\oone^k)$ that fixes $c$.  By
Lemma~\ref{l: first from bk3} this is equivalent to $T_c^{\pm} \in
\Fix_N(\hat \aone_k)$.  Since $\oone$ is forward rotationless, we may
assume that $k=1$.  This completes the proof of the first item.

If $\fG$ is a rotationless \rtt\ that represents $\phi$ and that  satisfies the conclusions of Theorem~\ref{rtt
existence}   then  each
$\Lambda \in \L(\oone)$ is associated to a unique EG stratum $H_r$; see Definition~3.1.12 of \cite{bfh:tits1} for details.  Lemma~\ref{some essential} and
Lemma~\ref{new iterates to}(2) imply that $\Lambda$ is the
accumulation set of an attractor $P \in \Fix_N(\hat \Phi)$ for some
principal automorphism representing $\oone$.  This implies (2).

For the third item, suppose that the free factor $F$ is
$\oone^k$-invariant for some $k \ge 1$.  If $F$ has rank one then it
is $\oone$-invariant by the first item of this lemma.  We may
therefore assume that $F$ has rank at least two.  Let ${\mathcal C}$
be the set of bi-infinite lines $ \gamma$ that are carried by $F$ and
for which there exist a principal lift $\aone$ of an iterate of
$\oone$ and a lift $\ti \gamma$ of $\gamma$ whose endpoints are
contained in $\Fix_N(\hat \Phi)$.  Since $\oone$ is forward
rotationless, each $\gamma$ is $\oone$-invariant, so ${\mathcal C}$ is
$\oone$-invariant.  Obviously ${\mathcal C}$ is carried by $F$ so to
prove that $F$ is $\oone$-invariant it suffices, by Corollary~\ref{ffs
invariance}, to show that no proper $\oone$-invariant free factor
system $\F$ of $F$ carries ${\mathcal C}$.

Suppose to the contrary that such an $\F$ exists.  By Theorem~\ref{rtt
existence} there is a relative train track map $g : G' \to G'$
representing $\oone^k|F$ in which $\F$ is represented by a proper
filtration element $G'_r \subset G'$.  After replacing $\oone^k|F$ and
$f$ by iterates we may assume that they are (forward) rotationless.
There is an principal vertex $v \in G'$ whose link contains an edge
$E$ of $G' \setminus G'_r$ that determines a fixed direction.  This
follows from Lemma~\ref{some essential} if there is an EG stratum in
$G' \setminus G'_r$ and from the definition of principal otherwise.
Lemma~\ref{new iterates to} and the fact that there are at least two
periodic directions based at $v$ imply that there is a principal lift
$\ti f' : \Gamma' \to \Gamma'$ and a line $\ti \gamma$ whose endpoints
are contained in $\Fix_N(\hat f')$ and whose projected image $\gamma$
crosses $E$ and so is not carried by $G'_r$.  The automorphism $\aone'
\in \PA(\oone|F)$ determined by $\ti f'$ extends to an element $\aone
\in \PA(\oone^k)$ with $\Fix_N(\hat \aone') \subset \Fix_N(\hat
\aone)$.  Thus $\gamma \in {\mathcal C}$ in contradiction to our
choice of $\F$ and $G'_r$.  \endproof

\begin{cor} \label{restriction of rotationless}
If $\oone$ is forward rotationless and $F$ is a $\oone$-invariant free
factor, then $\theta:= \oone|F \in \Out(F)$ is forward rotationless.
\end{cor}
\proof Lemma~\ref{no iterates necessary}(1) handles the case that $F$
has rank one so we may assume that $F$ has rank at least two.  Choose
a \rtt\ $\fG$ and filtration $\filt$ satisfying the conclusions of
Theorem~\ref{rtt existence} and representing $\oone$ such that the
conjugacy class of $F$ is represented by $G_l$ for some $l$.
Proposition~\ref{rotationless} implies that $\fG$ is rotationless.
The restriction of $\fG$ to $G_l$ is a rotationless \rtt\ representing
$\theta$ and satisfying the conclusions of Theorem~\ref{rtt
existence}.  A second application of Proposition~\ref{rotationless}
implies that $\theta$ is forward rotationless.  \endproof

\section{Completely Split Relative Train Track Maps}  \label{section cs exists} 
For every $\oone \in \Out(F_n)$ there exists $k > 0$ such that
$\oone^k$ is represented by an improved relative train track map (IRT)
$\fG$ as defined in Theorem 5.1.5 of \cite{bfh:tits1}.  In this
section we update this theorem, replacing IRTs with \ct  s, by 
controlling   the iteration index $k$, adding a very useful property  called complete splitting, and by making  small changes to previous definitions.    Section~\ref{s:def} contains all the necessary definitions.  In section~\ref{s:hard} we show that completely splittings   in a \ct\ are hard splittings in the sense of  \cite{bg:hard}.  A detailed comparison of IRTs and \ct s is given in section~\ref{s:comparison}.  There is one new move needed for the construction of \ct s.  It is defined and in section~\ref{s:new move} and the existence theorem is stated and proved in section~\ref{s:existence}.    A few additional properties of \ct s are presented in section~\ref{s:extra}

\subsection {Definitions and Notation}  \label{s:def}   For $a \in F_n$, we let $[a]_u$ be the {\em unoriented conjugacy class determined by $a$} . Thus, $[a]_u = [b]_u$ if and only if $a$ is conjugate to either of $a$ or $\bar a$.  If $\sigma$ is a closed path then we let $[\sigma]_u$ be the {\em unoriented conjugacy class determined by $\sigma$}, thought of as a circuit. 

Suppose that $\fG$ is a rotationless \rtt\ with filtration $\filt$.      Each \noneg\ stratum $H_i$ is a  single edge $E_i$ satisfying  $f(E_i) = E_i u_i$ for some (necessarily  closed by (\noneg)) path $u_i \subset G_{i-1}$ that is  sometimes called the {\em suffix} for $E_i$.     If $u_i$ is a non-trivial Nielsen path,  then we say that $E_i$ is a {\em linear edge}.   In the linear case, we define the {\em axis} or {\em twistor} for $E_i$ to be $[w_i]_u$ where $w_i$ is root-free and $u_i = w_i^{d_i}$ for some $d_i \ne 0$.

\begin{definition} \label{def:exceptional}  If $E_i$ and $E_j$ are linear edges and if there are $m_i,m_j > 0$ and   a closed root-free Nielsen path $w$ such that $u_i =w^{m_i}$ and $u_j =w^{m_j}$  then a path of the form $E_i w^p \bar E_j$ with $p \in \Z$ is called an {\em exceptional path}.   
\end{definition}

\begin{remark}   If $E_i w^p \bar E_j$ is an exceptional path  then $f^k_\#(E_i w^p \bar E_j) = E_i w^{p+k(m_i-m_j)}\bar E_j$ for all $k \ge 0$.  It follows that $E_i w^p \bar E_j$ is a Nielsen path if and only if $m_i = m_j$, that $f_\#$ induces a height preserving bijection on the   set of exceptional paths and that the interior of $E_i w^p \bar E_j$ is an increasing union of pre-trivial paths.   
\end{remark}
 
 \begin{definition}   A filtration $\filt$ that satisfies the following
property is said to be {\em reduced (with respect to $\phi$)} : if a
free factor system $\F'$ is $\phi^k$-invariant for some $k > 0$ and if
$\F(G_{r-1}) \sqsubset \F' \sqsubset \F(G_{r})$ then either $\F' =
\F(G_{r-1})$ or $\F' = \F(G_{r})$.   
 \end{definition}

 \begin{definition} \label{def:cs path}   If $E$ in an edge in an irreducible 
stratum $H_r$  and  $k >0$ then a maximal subpath $\sigma$ of $f^k_\#(E)$ in a zero stratum $H_i$  is said to be {\em $r$-taken}  or just {\em taken}  if $r$ is irrelevant.  Note that if $H_i$ is enveloped by an \eg\  stratum $H_s$ then   $\sigma$ has endpoints in  $H_s$ and so is  a connecting path.  A non-trivial path or circuit $\sigma$ is {\em completely split} if it has a splitting,  called 
a {\it complete splitting}, into subpaths, each of which is either a
single edge in an irreducible stratum, an \iNp, an exceptional path or
a  taken connecting path in a zero stratum.
\end{definition}

\begin{definition} \label{def:cs} A \rtt\ is {\em completely split} if 
\begin{enumerate}
\item  $f(E)$ is completely split for each edge $E$ in each irreducible 
stratum.  
\item If $\sigma $ is a taken connecting path  in a zero stratum   then  $f_\#(\sigma)$ is completely split.   \end{enumerate}
\end{definition}

The next lemma states that if $\fG$ is   completely split then  $f_\#$ maps   completely split paths to completely split paths.

\begin{lemma} \label{stays cs}  If $\fG$ is completely split and $\sigma$ is a 
\cs\ path or circuit then $f_\#(\sigma)$ is \cs.  Moreover if $\sigma
= \sigma_1 \cdot \ldots \cdot \sigma_k$ is a complete splitting then
$f_\#(\sigma)$ has a complete splitting that refines $f_\#(\sigma) =
f_\#(\sigma_1) \cdot \ldots \cdot f_\#(\sigma_s)$
\end{lemma}

\proof This is immediate from the definitions and the fact that $f_\#$ carries 
\iNp s to \iNp s and exceptional paths to exceptional paths.
\endproof 

We now come to our main definition.   When equivalent descriptions of a property are available, for example in    (\eg\ Nielsen Paths),  we have chosen the one that is easiest to check.   

\begin{definition} \label{def:ct} A   \rtt\ $\fG$ and filtration $\F $ given by $ \filt$\ is said to be a {\em \ct}\ (for completely split improved relative train track map) if it satisfies the following properties.   
\begin{enumerate}

\item{\bf(Rotationless)} $\fG$ is rotationless.   (See Remark~\ref{ct satisfies 2.19}.)
\item{\bf(Completely Split)} $\fG$ is completely split. 

\item {\bf (Filtration)}  $\F$ is reduced.   
  The core of each filtration element is a filtration element.

\item {\bf (Vertices)} The endpoints of all indivisible periodic (necessarily fixed) Nielsen paths are (necessarily principal) vertices.   The terminal endpoint of each non-fixed \noneg\ edge is principal (and hence fixed). (See Remark~\ref{attaching vertex}  and Lemma~\ref{neg edges}.)

\item {\bf(Periodic Edges)} Each periodic edge is fixed and each endpoint of a fixed edge is principal.    If the unique edge $E_r$ in a fixed stratum $H_r$  is not a loop then $G_{r-1}$ is a core graph  and both ends of $E_r$ are contained in $G_{r-1}$. 

\item {\bf (Zero Strata)}  If $H_i$ is a zero stratum, then  $H_i$ is enveloped by an  \eg\ stratum $H_r$, each edge in $H_i$ is $r$-taken and each vertex in $H_i$ is contained in $H_r$ and has link  contained in $H_i \cup H_r$.

  \item {\bf(Linear Edges)}  For each linear $E_i$ there is a closed root-free Nielsen path $w_i$  such that $f(E_i) = E_i w_i^{d_i}$ for some $d_i \ne 0$.     If $E_i$ and $E_j$ are distinct linear edges  with the same axes  then $w_i = w_j$ and  $d_i \ne d_j$.   (See Remark~\ref{exceptional}.)

\item {\bf (\noneg\ Nielsen Paths)} 
If the highest edges in an \iNp\ $\sigma$  belong to an \noneg\ stratum then there is a  linear edge $E_i$ with $w_i$  as in (Linear Edges) and there exists $k \ne 0$ such that $\sigma = E_i w_i^k \bar E_i$.   

\item {\bf (\eg\ Nielsen Paths)} (See also Lemmas~\ref{still works bh} and \ref{still works bfh} and Corollaries~\ref{eg irt} and  \ref{condition for noneg NP})  If $H_r$ is \eg\  and $\rho$ is an \iNp\ of height $r$, then   $f|G_r = \theta\circ f_{r-1}\circ f_{r}$  where :
\begin{enumerate}
\item $f_r : G_r \to G^1$ is a composition of proper extended folds  defined by iteratively folding $\rho$.
\item $f_{r-1} : G^1 \to  G^2$ is a composition of   folds   involving edges in $G_{r-1}$.
 \item $\theta : G^2 \to G_r$ is a homeomorphism.
 \end{enumerate}
\end{enumerate}
\end{definition}

 \begin{remark}   \label{ct satisfies 2.19} A \ct\ satisfies the conclusions of Theorem~\ref{rtt existence}.  This is immediate from the definitions and from Lemma~\ref{neg edges}.
\end{remark}

\begin{remark}\label{attaching vertex}  It is an immediate consequence of (Vertices), Remark~\ref{rtt2} and the definitions that a vertex whose link contains edges in more than one irreducible stratum is principal.  
\end{remark}

 \begin{remark}  The tree replacement move introduced in  step 2 of the proof of Theorem~\ref{rtt existence}  shows that no taken path in zero strata  is more intrinsically an edge than any other.  This explains why we do not distinguish between such paths in Definition~\ref{def:cs path}.
 \end{remark}

\begin{remark}\label{exceptional}   Exceptional paths associated to an axis $[a]_u$ have the form $E_iw^p\bar E_j$   where $E_i$ and  $E_j$ are linear edges associated to $[a]_u$ such that $d_i$ and $d_j$ have the same sign,  where $w = w_i = w_j$ and where $ p \in \Z$.
\end{remark} 

\subsection{Hard Splittings}  \label{s:hard}

The first item in the following lemma establishes the  uniqueness of complete splittings; the second item (see also Corollary~\ref{induced splitting})  shows  that  a complete splitting  is a hard splitting as
defined in \cite{bg:hard}.

\begin{lemma}\label{cs is unique}
Suppose that $\fG$ is a \ct, that $\sigma$ is a circuit or path and that $\sigma = \sigma_1   \ldots \sigma_m$ is a decomposition into
subpaths, each of which is either a single edge in an irreducible
stratum, an indivisible Nielsen path, an exceptional path or a taken
connecting path in a zero stratum.  Suppose also that  each turn $(\bar \sigma_i,  \sigma_{i+1})$ is legal.       Then \begin{enumerate}
\item $\sigma = \sigma_1 \cdot
\ldots \cdot \sigma_m$ is the unique  complete splitting
of $\sigma$.
\item each pre-trivial subpath $\tau$ of $\sigma$ is contained in a single  $\sigma_i$.
\item  a subpath of $\sigma$ that has the same height as $\sigma$ and is
either a fixed edge or an \iNp\ equals $\sigma_i$ for some
$i$.
\end{enumerate}
\end{lemma}
\proof      Let $\ti \sigma = \ti \sigma_1\cdot \ldots \cdot \ti \sigma_m$
be a lift of $\sigma$ and let $\ti f : \Gamma \to \Gamma$ be a lift of
$\fG$.   The main step in the proof is to establish the following property.
\begin{enumeratecontinue} 
\item   If   $\sigma_i$ is not a taken connecting path in a zero stratum then for each $k \ge 0$ there exist non-trivial initial and terminal subpaths $\ti \alpha_{i,k}$ and $\ti \beta_{i,k}$ of $\ti \sigma_i$ such that  $ \ti f^{-k}(\ti f^k(\ti x) )\cap \ti \sigma
= \{\ti x\}$ for all $\ti x \in \Int (\ti \alpha_{i,k}) \cup \Int(\ti \beta_{i,k})$.
\end{enumeratecontinue}  
The proof of (4) is    by double induction, first on  $k$  and then on $m$.    The $k=0$ case is obvious  so we may assume that (4) holds  for any iterate less than $k$.    

To establish the second base case, assume that $m=1$ or equivalently  that $\sigma = \sigma_i$.    If  $\sigma_i$ is exceptional then (4)  is clear (cf. Lemma~4.1.4 of \cite{bfh:tits1}).  If   $\sigma_i$ is an \iNp\   then (4) follows from (\noneg\ Nielsen Paths)   and  Lemma~\ref{nielsen paths in egs}(2).  The remaining possibility is that   $\sigma_i$ is an   edge $E$ in  an irreducible stratum and in this case we make use of the inductive hypothesis that (4) holds in general for any iterate less than $k$.       The first and last terms in the complete splitting of  $f(E)$    are not   connecting paths in   zero strata.  By the inductive hypothesis there exist   initial and terminal subpaths $\ti \alpha' $ and $\ti \beta' $ of $\ti f( \ti E)$ such that $\ti f^{-(k-1)}(\ti f^{k-1}(\ti x)) \cap \ti f(\ti E)= \{\ti x\}$ for all $\ti x \in \Int(\ti \alpha' ) \cup \Int(\ti \beta')$.   Since $\ti f|\ti E$ is an embedding, we can pull $\ti \alpha'$ and $\ti \beta'$ back to initial and terminal subpaths  $\ti \alpha_{i,k}$  and $\ti \beta_{i,k}$ of $\ti E$ that  satisfy (4).  This completes    the $m =1$ case.

Suppose now that (4) holds for $k$ and  $\sigma$ if the  complete splitting of $\sigma$  has fewer than $m\ge 2$ terms.   As a first case suppose that $\sigma_1$ is a   taken connecting path in a zero stratum $H_p$.  By (Zero Strata),    $\sigma_2$ is an edge in a \eg\ stratum $H_r$ with $r > p$. Define $\ti \alpha_{i,k}$ and $\ti \beta_{i,k}$ using $\ti \sigma_2   \ldots \ti \sigma_m$ in place of $\ti \sigma$.    Then  $\ti f^k(\ti \sigma_1) \cap \Int(\ti f^k(\ti \alpha_{2,k}))  = \emptyset$ because $\ti f^k(\ti \sigma_1)$ has height $< r$ and      $\ti f^k(\ti \alpha_{2,k})$ is an embedded path whose   initial direction  has height $r$.    Since  $\Int(\ti f^k(\ti \alpha_{2,k}))$ separates $\ti f^k(\ti \sigma_1)$ from  each $\ti f^k (\ti \beta_{i,k}) $ and from $\ti f^k(\ti \alpha_{i,k})$ for $i > 2$, $\ti f^k(\ti \sigma_1)$ is disjoint from each of these sets.  This proves that  each 
 $\ti \alpha_{i,k}$ and $\ti \beta_{i,k}$ satisfies (4) with respect to $\sigma$.  
 
 As a second case,   suppose that $\sigma_2$ is a taken connecting path in a zero stratum $H_p$.  By (Zero Strata),  $\sigma_1$  and $\sigma_3$ (if $m \ge 3$) are edges  in an \eg\ stratum $H_r$ with $r > p$.   Define  $\ti \alpha_{1,k}$ and $\ti \beta_{1,k}$ using $\ti \sigma_1$ in place of $\ti \sigma$.  For $i > 2$, define $\ti \alpha_{i,k}$ and $\ti \beta_{i,k}$ using $\ti \sigma_2   \ldots \ti \sigma_m$ in place of $\ti \sigma$.    As in the previous case, $\Int(\ti f^k(\beta_{1,k})) \cap  \ti f^k(\ti \sigma_2) = \emptyset$ and   $\Int(\ti f^k(\alpha_{3,k})) \cap  \ti f^k(\ti \sigma_2) = \emptyset$.  Also as in the previous case, this implies that  each 
 $\ti \alpha_{i,k}$ and $\ti \beta_{i,k}$ satisfies (4) with respect to $\sigma$.  
  
 The final case is that neither $\sigma_1$ nor  $\sigma_2$ is a taken connecting path in a zero stratum.  Define  $\ti \alpha_{1,k}$ and $\ti \beta_{1,k}$ using $\ti \sigma_1$ in place of $\ti \sigma$.  For $i > 2$, define $\ti \alpha_{i,k}$ and $\ti \beta_{i,k}$ using $\ti \sigma_2   \ldots \ti \sigma_m$ in place of $\ti \sigma$.   Since the turn $(\bar \sigma_1,  \sigma_{2})$ is legal,  the interiors of $\ti f^k(\ti \alpha_{1,k})$ and  $\ti f^k(\ti \beta_{2,k})$ are disjoint.  The proof now concludes as in the previous two cases.    This completes the induction step and so the proof of (4).

If $\tau$ is a pre-trivial path then there exists $k > 0$ such that $f^k_\#(\tau)$ is trivial.  For each $\ti x \in
\ti \tau$ there exists $\ti y \ne \ti x$ in $\ti \tau$ such that $\ti
f^k(\ti x) = \ti f^k(\ti y)$.  If $\sigma_i$ is not a connecting path in a zero stratum  and if $\tau$ intersects $\Int(\sigma_i)$ then $\tau \subset \Int(\sigma_i)$ by (4).  Since this applies to at least one of any pair of consecutive $\sigma_i$'s we have proved (2).  It follows that       $\sigma = \sigma_1 \cdot
\ldots \cdot \sigma_m$ is a splitting, and hence a  complete splitting, 
of $\sigma$.

Suppose that  $\sigma = \sigma_1' \cdot \ldots \cdot \sigma'_q$ is also a complete splitting.  If $\sigma'_i$ is an exceptional path or an \iNp\  then (by inspection in the \noneg\ case and by  Lemma~\ref{nielsen paths in egs}(2) in the \eg\ case) the interior of $\sigma'_i$ is the increasing union of pre-trivial subpaths.    Item (2) implies that $\sigma'_i$ is contained in some $\sigma_j$.  Since $\sigma_j$    is  not a single edge and is not contained in a zero stratum, it must be  an \iNp\ or an exceptional path.  By symmetry,  $\sigma'_i = \sigma_j$.     The terms  that are taken connecting paths in zero strata  are then
characterized as the maximal subpaths, in the complement of the \iNp s
and exceptional paths, that are contained in zero strata.  All
remaining edges are terms in the complete splitting.  This proves that complete splittings are unique and so completes the proof of (1).

A fixed edge of maximal height in $\sigma$ is not contained in a taken
connecting path of a zero stratum, an \iNp\ or an exceptional path in $\sigma$ and so
must be a term in the complete splitting of $\sigma$.  An \iNp\ in
$\sigma$ must be contained in a single $\sigma_i$ by (2).  If it has
maximal height then, by inspection of the four possibilities for
$\sigma_i$ it must be all of $\sigma_i$.  This proves (3).  \endproof

\begin{cor} \label{induced splitting}
Assume that $\fG$ is a \ct\ and that $\sigma=\sigma_1 \cdot \ldots
\cdot \sigma_s$ is the complete splitting of a path $\sigma \subset
G$.  If $\tau$ is an initial segment of $\sigma$ with terminal
endpoint in $\sigma_j$ then $\tau = \sigma_1 \cdot \ldots \cdot
\sigma_{j-1}\cdot \mu_j$ is a splitting where $\mu_j$ is the initial
segment of $\sigma_j$ that is contained in $\tau$.  In particular if
$\tau$ is a non-trivial Nielsen path then $\sigma_i$ is a Nielsen path
for all $i \le j$ and if $\sigma_j$ is not a single fixed edge then
$\mu_j = \sigma_j$.
\end{cor}

\proof The main statement follows immediately from  
Lemma~\ref{cs is unique}(2).  The statement about Nielsen paths then
follows from the fact that no proper non-trivial initial segment of a
non-fixed term in a complete splitting of any path is a Nielsen path.
\endproof

\subsection{\ct\ versus IRT} \label{s:comparison}

Theorem~5.1.5 of \cite{bfh:tits1} is both the definition of, and the existence theorem for, an \irt\ map.  There are 
  eight bulleted items in the statement of the theorem, the  last seven of which should be considered the definition.   For notational convenience, we refer to these as (IRT-1) through (IRT-7).    In this section we discuss the extent to which a \ct\ $\fG$ satisfies these seven items.    By the end of this section we will have verified that \ct s satisfy  all of the important properties of IRTs.
  
     (IRT-1)  is  that $\F$ is reduced, which is  part of (Filtration).    The following lemma states that every \ct\ satisfies (IRT-2).  

  \begin{lemma} If $\fG$ is a \ct\ then every periodic Nielsen path has period one. 
  \end{lemma}
 
 \proof   Each periodic Nielsen path is a concatenation of periodic edges and indivisible periodic Nielsen path.  The former  has period one by  (Periodic Edges)  and the latter has period one by (Vertices) and Lemma~\ref{period one}. 
 \endproof
 
The next lemma shows that a \ct\ satisfies most of  (IRT-3).  The exception is that  there may be some vertices $v$ for which $f(v)$ is not  fixed.   

\begin{lemma} If $\fG$ is a \ct\  then every vertex $v \in G$ has at least two gates.  If the  link of $v$ is not contained in $H^z_r$ for some  \eg\ stratum $H_r$ then $v$ is principal and hence fixed. 
\end{lemma}

\proof      If $(d_1,d_2)$ is an illegal turn then either one of the $d_i$'s is the terminal end of a non-fixed \noneg\ edge  or both   $d_1$ and $d_2$ belong to $H^z_r$ for some  \eg\ stratum $H_r$.  (Vertices), (Zero Strata) and  Lemma~\ref{two directions for EG}    imply that the vertex in both of these case has two gates.  At any other vertex the number of gates equals the valence.  This proves the first statement of the lemma. 

  It follows from (Periodic Edges) and the definition of principal vertex that if $v$ is periodic and the link of $v$ is not contained in a single \eg\ stratum then $v$ is principal.    If $v$ is not periodic then  its link is contained in  some $H^z_r$ by Remark~\ref{rtt2} and   (Zero Strata).
\endproof

       The difference between (IRT-4) and the conclusion of the next lemma is that a zero stratum in a IRT can be the union of contractible components.
    
    \begin{lemma} \label{ct contractible components} Assume that $\fG$ is a \ct.  Then $G_i$ has  a contractible component   if and only if $H_i$  is a zero stratum.
\end{lemma}

\proof  The if direction follows from (Zero Strata).  For   the only if direction we assume that $H_i$ is not a zero stratum and prove, by induction up the filtration, that  every component of $G_{i}$ is non-contractible.        For the base case, $H_1 = G_1$ is  either \eg\  or periodic and so is connected and not contractible, by   Lemma~\ref{two directions for EG} in the former case and (Periodic Edges) in the latter. We now consider the inductive step.       If some component of $G_{i-1}$ is  contractible then $H_{i-1}$ is a zero stratum and  (Zero Strata) implies that every component of $G_{i}$ is non-contractible.   If every component of $G_{i-1}$ is non-contractible then (Periodic Edges) and Lemma~\ref{two directions for EG} complete the proof.
\endproof 

There are two differences between (IRT-5) and (Zero Strata).  The first is that  an IRT can have a vertex whose link is contained in a zero stratum but a \ct\ cannot.   The second is that the restriction of an IRT to   a  zero stratum   is always an immersion but this need not be true for a \ct.     We have replaced the immersion condition with  Definition~\ref{def:cs}(2) and the assumption that every edge in a zero strata is $r$-taken.   The primary motivation for removing the immersion condition is that it  lacks robustness.  For example, it need not hold for $f^2|\sigma$.  Also, since the main  application of relative train track maps is in analyzing the action of the induced map $f_\#$  on paths with endpoints at vertices, it makes sense to make definitions that focus on $f_\#$ and not on $f$.
      
    Corollary~\ref{eg irt} below implies that a  \ct\  satisfies  (IRT-7).  In the definition of a \ct, we have replaced a list of properties satisfied by \iNp s corresponding to \eg\ strata (see the statement of  Corollary~\ref{eg irt})  with the underlying property (\eg\ Nielsen Paths)  from which these properties were derived.    One advantage of this is that it is easier to deduce additional properties as in Lemma~\ref{no zero strata}.        

 \begin{lemma}  \label{stable}  Suppose that   $H_r$ is an \eg\ stratum of   a \rtt\ $\fG$, that $\rho$ is an \iNp\ of height $r$ and that   $\rho$ and     $H_r$ satisfy the conclusions of  (\eg\ Nielsen Paths).  Then  
 the illegal turn at each \iNp\ obtained by iteratively folding  $\rho$ is   proper. 
         \end{lemma}
  
  \proof   Assume the notation of (\eg\ Nielsen Paths).   Since $f$ and $f_r$ makes exactly the same identifications on edges in $H_r$, each $f$-illegal turn in $H_r$ is folded by $f_r$.  Since $f_r$ is defined by iteratively folding $\rho$,  say $K$ times, the illegal turn of $\rho$  is the only  illegal turn in $H_r$.   Let $\rho(k)$ be the \iNp\  determined by iteratively folding $\rho = \rho(0)$ $k$ times.   By hypothesis, the fold at the illegal turn of $\rho(k)$ is proper for $k < K$.   It suffices to show that this holds for all $k$.   
    
Let $\rho(k) = a_{1,k}b_{1,k} \ldots b_{m-1,k}a_{m,k}$ be the decomposition of $\rho(k)$ into subpaths $a_{i,k} \subset H_r$ and maximal subpaths $b_{i,k} \subset G_{r-1}$.   Denote the   ordered sequence of the $a_{i,k}$'s   by $S_k$.   Then $S_k$ determines the type (partial, proper, improper) of the extended fold of $\rho(k)$ and, assuming that the fold is full, also determines $S_{k+1}$.    Since $f|G_r = \theta\circ f_{r-1}\circ f_{r}$,  $S_K = S_0$.    It follows that the sequence of $S_k$'s  is periodic with period $K$   and hence that  the fold at the illegal turn of $\rho(k)$ is proper for all $k$.  
  \endproof

   \begin{lemma}  \label{still works bh} Theorem~5.15 of \cite{bh:tracks}  remains true if the hypothesis   that $\fG$ is stable is replaced by the hypothesis that 
  for each \eg\ stratum $H_r$  there is an \iNp\  $\rho$ of height $r$  such  that   $\rho$ and     $H_r$ satisfy the conclusions of  (\eg\ Nielsen Paths). 
   \end{lemma}

\proof     The proof of Theorem~5.15  has two parts.  The first is a reduction to the case that  the illegal turn at each \iNp\ obtained by iteratively folding  $\rho$ is full and proper.   The second    is the observation that in this case    the proof  for the special case that $\fG$ is irreducible   given in Lemma~3.9  of \cite{bh:tracks}  applies to the general case as well.  This lemma therefore follows from Lemma~\ref{stable}.
\endproof

\begin{lemma}  \label{still works bfh}  Proposition~5.3.1 of \cite{bfh:tits1}  remains true if the hypothesis  $\fG$ is $\F$-Nielsen minimized is replaced by the hypothesis that  $H_r$ satisfies   (\eg\ Nielsen Paths).
\end{lemma}

\proof   The proof of    Proposition~5.3.1 of \cite{bfh:tits1} makes use of Lemmas~5.3.6, 5.3.7, 5.3.8 and 5.3.9 of that paper.     Lemma~5.3.6 states that if $\fG$ is $ \F$-Nielsen minimized and if $\rho_r$  crosses every edge of $H_r$ exactly twice then $H_r$ satisfies   (\eg\ Nielsen Paths).  The remaining three lemmas use (\eg\ Nielsen Paths) but do not refer   directly to being $\F$-Nielsen minimized.  
\endproof

The next corollary refers to {\em geometric strata}; complete details can be found in Definition~5.1.4 of \cite{bfh:tits1}. 

\begin{corollary} \label{eg irt} Suppose that $\fG$ is a \rtt\ and that   (\eg\ Nielsen Paths) holds for the \eg\ stratum  $H_r$. 
 Then the following properties are satisfied.
\begin{description}

\item[eg-(i)] There is at most one \iNp\ $\rho_r
\subset G_i$ that intersects $H_r$
non-trivially.  The initial edges of $\rho_r$ and 
$\bar
\rho_r$ are distinct   edges in
$H_r$.   
\item[eg-(ii)]  If $\rho_r
\subset G_r$ is an \iNp\ that intersects $H_r$ non-trivially 
and if $H_r$ is
not geometric, then   there is an
edge  $E$ of $H_r$ that $\rho_r$ crosses exactly
once.
\item[eg-(iii)]  If $H_r$ is geometric then there is an 
\iNp\ $\rho_r
\subset G_r$ that intersects $H_r$ non-trivially and 
satisfies the following
properties : (i)
$\rho_r$  is a closed path
with  basepoint  not contained in $G_{r-1}$;  (ii)
the circuit determined
by $\rho_r$ corresponds to the unattached
peripheral curve $\rho^*$ of $S$;
and (iii) the surface $S$ is connected.
\end{description}
In particular, $H_r$ satisfies the  \eg\ properties of an \irt.
\end{corollary}

 \proof   Theorem~5.15 of \cite{bh:tracks}, which applies here by Lemma~\ref{still works bh},  implies that there is at most one \iNp\ $\rho_r
\subset G_i$ that intersects $H_r$
non-trivially and if such  a $\rho_r$ exists then it either crosses every edge in $H_r$ exactly twice or crosses some edge of $H_r$ exactly once.   Lemma~5.1.7 of \cite{bfh:tits1} implies that if $\rho_r$ crosses some edge of $H_r$ exactly once then $\rho_r$ is not a closed path.  This completes the proof of eg-(i).  The remaining two items follow from   Proposition~5.3.1 of \cite{bfh:tits1} which applies here by Lemma~\ref{still works bfh}.   
 \endproof
 
 \begin{remark}  \label{distinct endpoints}Item eg-(ii) of Corollary~\ref{eg irt} and \cite[Lemma~5.1.7]{bfh:tits1} imply that if $H_r$ is an \eg\ stratum of a \ct\ that  is not geometric and if $\rho$ is an \iNp\ of height $r$ then $\rho$ has distinct endpoints.
 \end{remark}
 
 The remaining item (IRT-6) concerns \noneg\ strata and has    three parts.  The first two  statements of the next lemma shows  that a \ct\ satisfies the first two parts of  (IRT-6).  Corollary~\ref{noneg properties} shows that a \ct\ satisfies the third  part of  (IRT-6).

 \begin{lemma} \label{neg edges} If   $\fG$ is a \ct\ and $H_i$ is \noneg\ then $H_i$ is a single edge $E_i$.  If $E_i$ is not contained in $\Fix(f)$ then   there is a non-trivial closed path $u_i \subset G_{i-1}$ such that $f(E_i) = E_i \cdot u_i$.   Moreover  $u_i$ forms a circuit and the turn $(u_i,\bar u_i)$ is   legal.
 \end{lemma}

\proof     If $H_i$ consists of periodic edges then the lemma follows from   (Periodic Edges).   Otherwise   (Rotationless), (Completely Split) and  (Vertices) imply that $H_i$ is a single edge $E_i$ and that there is a non-trivial closed path  $u_i$  such that $f(E_i) = E_i  u_i$ is completely split.  To prove that $f(E_i) = E_i \cdot u_i$ we must show that the first term $\sigma_1$ in the complete splitting of $E_iu_i$  is the single  edge $E_i$.  It is obviously not contained in a zero stratum and is not a Nielsen path by (\noneg\ Nielsen Paths).  It remains to show that $\sigma_1$ is not an exceptional path and for this there is no loss in assuming that $E_i$ is linear.   In the notation of (Linear Edges), $f(E_i) = E_i w_i^{d_i}$, no initial segment of which is an exceptional path by Remark~\ref{exceptional}.    This completes the proof that $f(E_i) = E_i \cdot u_i$.

The turn $(Df^{k-1}(\bar u_i), Df^k(u_i))$ is the $Df^{k-1}$ image of
the legal turn $(\bar E_i, u_i)$ and is therefore legal for all $k \ge
1$.  Since $f$  is rotationless and since $v$ is  principal by (Vertices), $Df^k(d)$ is independent
of $k$ for all directions $d$ based at $v$ and all sufficiently large $k$.  It
follows that $(Df^{k}(\bar u_i), Df^k(u_i))$ is legal for all
sufficiently large $k$ and hence that $(u_i,\bar u_i)$   is   legal.  In particular,  $(u_i,\bar u_i)$ is non-degenerate which implies that $u_i$ forms a circuit.
\endproof

 \begin{lemma}  \label{two parts of prop} Suppose that $E_i$ is the unique edge of height $i$ in a rotationless \rtt\ $\fG$,   that $f(E_i) = E_i \cdot u_i$ for some non-trivial closed path $u_i \subset G_{i-1}$ and that every periodic Nielsen path with height less than $i$ has period one. Suppose further that either there are no Nielsen paths of height $i$  or    $E_i$ is a linear edge and  all Nielsen paths of height $i$  have the form   $\sigma = E_i w_i^k \bar E_i$ where  $k \ne 0$ and where $w_i$ is root-free and $u_i =w_i^{d_i}$ for some $d_i \ne 0$.     Let $h :\Gamma_{i-1} \to \Gamma_{i-1}$ be as in Notation~\ref{defining h}.   Then
 \begin{itemize}
 \item $\Fix(h) = \emptyset$.
 \item  $E_i$ is a linear edge  if and only if  there is a covering translation $T:\Gamma_{i-1} \to \Gamma_{i-1}$ that commutes with $h$ and whose axis covers $u_i$.   
 \end{itemize}
\end{lemma}
 
 \proof   Let $\ti f : \Gamma \to \Gamma$ and $\ti E_i$ be as in  Notation~\ref{defining h}.  Thus $\ti f(\ti E_i) = \ti E_i \cdot\ti u_i$ where $\ti u_i \subset \Gamma_{i-1}$  is a lift of $u_i$ and $h$ maps the initial endpoint $\ti x_1$ of $\ti u_i$ to   the terminal endpoint $\ti x_2$ of $\ti u_i$.   If $\ti v \in \Fix(h)$  and $\ti \gamma$ is the path  from $\ti x_1$ to $\ti v$ then $\ti E_i \ti \gamma$ is a Nielsen path for $\ti f$.  But then $E_i\gamma$ is a Nielsen path of height $i$ for $f$ that is not of the form  $E_i w_i^k \bar E_i$.  This contradiction verifies the first item.
 
  If $u_i$ is a Nielsen path  and  $T:\Gamma_{i-1} \to \Gamma_{i-1}$ is the  covering translation that maps $\ti x_1$ to $\ti x_2$,  then $Th(\ti x_1) = hT (\ti x_1)$ is  the terminal endpoint of the lift of $u_i$ that begins at $\ti x_2$.   Thus $T$ commutes with $h$.  For the converse suppose that $h$ commutes with some covering translation $T:\Gamma_{i-1} \to \Gamma_{i-1}$.    Corollary~\ref{non-trivial fixed set} implies that $h$ is not a principal  lift of $f|G_{i-1}$ and hence that the endpoints of the axis of $T$ are the only fixed points in $\partial \Gamma_{i-1}$.    On the other hand, the ray $\ti u_i \cdot h_\#(\ti u_i) \cdot h^2_\#(\ti u_i)\dots$ converges to a fixed point in $\partial \Gamma_{i-1}$.  The end of this ray is therefore contained in the axis of $T$.  It follows that $u_i$ is a periodic Nielsen path  and hence a Nielsen path and that the  axis of $T$ covers $u_i$.
\endproof

Recall (Definition~4.1.3 of \cite{bfh:tits1}) that if $H_i$ is a \noneg\ strata with unique edge $E_i$ then paths of the form $E_i \gamma \bar E_i, E_i \gamma$ or $\gamma \bar E_i$ where $\gamma \subset G_{i-1}$ are called {\em basic paths of height $i$.}
 
 \begin{corollary}  \label{noneg properties} Suppose that $\fG$ is a \ct\ and that $H_i$ is a \noneg\ strata with unique edge $E_i$.   If $\sigma \subset G_i$ is a basic 
path of height $i$
 that does not split as a 
concatenation of two
basic paths of height $i$ or as a concatenation of a 
basic path of height
$i$ with a path contained in
$G_{i-1}$, then   either : (i) some
$f_\#^k(\sigma)$ splits into pieces, one of which
equals $E_i$ or $\bar E_i$;  or (ii) $u_i$ is
a Nielsen path and some $f_\#^k(\sigma)$ is
is an exceptional path of height $i$.
\end{corollary}

\proof  Lemma~\ref{two parts of prop} and Corollary~\ref{induced splitting}  imply that $\fG$ satisfies the 
conclusions of Proposition~5.4.3 of \cite{bfh:tits1}.  In conjunction with Corollary~\ref{eg irt},  Lemma~\ref{neg edges} and Lemma~\ref{period one} we see that $\fG$ satisfies the hypotheses of Proposition~5.5.1 of \cite{bfh:tits1}, from which the  corollary follows.
 \endproof
 
 We conclude this subsection with two additional properties of   \ct s.
 
  \begin{lemma}  \label{no zero strata}  Suppose that $\fG$ is a rotationless \rtt, that $H_r$ is an \eg\ stratum satisfying (\eg\ Nielsen Paths),  and that $\rho $ is   an \iNp\ of height $r$.  Then 
 \begin{enumerate}
  \item $H^z_r = H_r$.  
 \item If $\rho =a_1 b_1 \ldots b_la_{l+1}$ is the decomposition  into subpaths $a_i$ of height $r$    and maximal subpaths $b_i$ of height less than $r$ then each $b_i$ is a Nielsen path.
 \item  If $E$ is an edge  of $H_r$ then each   maximal subpath of $f(E)$ in $G_{r-1}$ is one of the  $b_i$'s from (2).  In particular $f(E)$ splits into edges in $H_r$ and Nielsen paths in $G_{r-1}$.
\end{enumerate}
 \end{lemma}
 
 \proof      The maps  $f_r, f_{r-1}$ and $\theta$  induce  bijections on the set of   components in the filtration element of height $r-1$.  It follows that  $f = \theta f_{r-1} f_r$ induces a bijection on the set of   components of $G_{r-1}$ and hence that each component of $G_{r-1}$ is non-wandering.  This proves (1).   
  
  For (2), let $(f_r)_\#(\rho)  =a'_1 b_1' \ldots a'_mb'_ma'_{m+1}$ be the decomposition into subpaths $a'_j$ of height $r$ and maximal subpaths $b'_j$ of height less than $r$.  It is an immediate consequence (see the proof of Lemma~5.3.3 of \cite{bfh:tits1})  of the definition of an extended fold that   the set of distinct $b'_j$'s is contained in the set of distinct $b_i$'s.   Now let  $(\theta f_{r-1})_\#(a'_1 b_1' \ldots a'_mb'_{m+1})  =c_1 d_1 \ldots c _pd_pc_{p+1}$ be the decomposition into subpaths $c_k$ of height $r$ and maximal subpaths $d_k$ of height less than $r$. Then for each $k$ there exists $j$ such that   $d_k = (\theta f_{r-1})_\#(b'_j)$. Combining this with the fact that   $a_1 b_1 \ldots b_la_{l+1} =c_1 d_1 \ldots c _pd_pc_{p+1}$,  we conclude that $f_\#$ permutes the $b_i$'s.  Since $f$ is rotationless, each $b_i$ is a Nielsen path.  

 If $E$ is an edge  of $H_r$ then, by construction,  each maximal subpath of $f_r(E)$ in $G_{r-1}$ is a $b_i$.  By (2), each $b_i$ is a Nielsen path for $f$ and hence for $\theta f_{r-1}$.  This completes the proof of (3).
  
\endproof

 \begin{lemma} \label{iterates split}   If  $\fG$ is a \ct \ and 
 $\sigma \subset G_r$ is a path  with endpoints at vertices then 
$f^k_\#(\sigma)$ 
is completely split for all sufficiently large $k$.
\end{lemma}

\proof  The proof is by induction on the height of $\sigma$.  The height zero 
case is vacuously true so suppose that $\sigma$ has height $j\ge 1$ and that 
the lemma holds for all paths of height less than $j$. By Lemmas~\ref{stays cs}, \ref{cs is unique} 
and the inductive hypothesis,  it suffices to show that some $f^k(\sigma)$ has a 
splitting into subpaths that are either \cs\ or contained in $G_{j-1}$.  This is 
immediate if $H_j$ is a zero stratum or if $H_j$ is a single fixed edge. If 
$H_j$ is \noneg\ then   $\sigma$ has a splitting into  basic paths of height $j$ and subpaths  in $G_{j-1}$ by Lemma 4.1.4  of    \cite{bfh:tits1}.  The desired splitting of $\sigma$ therefore follows from Lemma~\ref{noneg properties}.     If $H_j$ is \eg, then  Lemmas~4.2.6 and 
4.2.5 of \cite{bfh:tits1} imply that some $f^k(\sigma)$ 
splits into pieces, each of which is either $j$-legal or a Nielsen path and  
Lemma 4.2.1 of \cite{bfh:tits1}  implies that the $j$-legal paths in $G_j$ split 
into single edges in 
$H_j$ and subpaths in $G_{j-1}$.  
\endproof

  \subsection{A new move} \label{s:new move}    
We make use of a move that plays the same role for zero and \eg\
strata that sliding (section 5.4 of \cite{bfh:tits1}) does for \noneg\
strata.  See item (7) of Lemma~\ref{new move} below for its main
application.

\begin{definition} \label{remarking}
Suppose that $\fG$ is a rotationless \rtt\ satisfying the conclusions
of Theorem~\ref{rtt existence} with respect to the filtration $\filt$,
that $1 \le j \le N$, that every component of $G_j$ is
non-contractible and that $f$ fixes every vertex in $G_j$ whose link
is not contained in $G_j$.  Define a homotopy equivalence $g : G \to
G$ by $g|G_j = f|G_j$ and $g|(G \setminus G_j)=$ identity.

Define $G'$ from $G$ by changing the marking via $g$.  More precisely,
if $X$ is the underlying graph of $G$ and $\tau :R_n \to X$ is the
marking that defines $G$, then $g \tau:R_n \to X$ is the marking that
defines $G'$.  Since $G$ and $G'$ have the same underlying graph,
there is a natural identification of $G$ with $G'$ and we use this
when discussing edges and strata.

Define $f' : G' \to G'$ by $f'|G_j' = f|G_j$ and $f'(E) = (gf)_\#(E) $
for all edges $E$ in $H_i$ with $i > j$.  

We say that $f' : G' \to G'$ is obtained from $\fG$ by {\it changing
the marking on $G_j$ via $f$.}
\end{definition}

The following lemma is the analog of Lemma 5.4.1 of \cite{bfh:tits1}.

\begin{lemma} \label{new move}   Suppose that   $f' : G' \to G'$ is obtained from 
$\fG$ by changing the marking on $G_j$ via $f$.  Then :
\begin{enumerate}
\item  $f'|G_j = f|G_j$.   
\item  for every path $\sigma \subset G$ with endpoints   at 
vertices  and for every $k > 0$, $g_\# f^k_\# (\sigma) = (f')^k_\#g_\# 
(\sigma)$. 
\item  $f' : G' \to G'$ is a homotopy equivalence that determines the same 
element of $\Out(F_n)$ as $\fG$.
\item  there is a one-to-one  
correspondence between   Nielsen paths for $f$ and   Nielsen 
paths for $f'$. 
\item  $f' : G' \to G'$ is a rotationless \rtt\ satisfying the
conclusions of Theorem~\ref{rtt existence} with respect to $\filt$.
\end{enumerate}
\end{lemma}

\proof Item (1) is immediate from the definitions as is the fact that
$f'$ preserves the filtration $\filt$.  Also immediate are:
\begin{enumeratecontinue}
\item  If $f(x) \ne f'(x)$ then $x \not \in G_j$ and $f(x), f'(x)
\in G_j$.  In particular, $\Fix(f) = \Fix(f') \subset \Fix(g)$,
$\Per(f) = \Per(f') \subset \Per(g)$ and $Df$ and $Df'$ have the same
fixed and periodic directions.
\item Suppose that  $E$ is an edge in $H_i$ for $i > j$ and that $f(E) = 
\mu_1 \nu_1 \mu_2 
\dots \nu_{k-1} \mu_k$ where the $\nu_l$'s are the maximal subpaths in $G_j$ 
and where $\mu_1$ and  $\mu_k$ may be trivial.  Then  $f'(E) = \mu_1
f_\#(\nu_1) \mu_2 \dots f_\#(\nu_{k-1}) \mu_k$ and all the $f_\#(\nu_l)$'s are 
non-trivial. (The non-triviality follows from the fact that  $f$ fixes the 
 endpoints of each $\nu_l$.)
\end{enumeratecontinue}
 which implies
 \begin{enumeratecontinue}
 \item each stratum  $H_i$ has the same type (\eg, 
\noneg, zero) for $f$ as for $f'$.
\end{enumeratecontinue}

To verify (2), it suffices to assume that $k = 1$ and that $\sigma$ is
a single edge $E$.  If $E \subset G_j$ then $g_\#f_\#(E)=f_\#(f_\#(E))
= f'_\#g_\#(E)$. If $E \subset G_i$ for $i > j$ then $g_\#f_\#(E) =
f'_\#(E) = f'_\#g_\#(E)$. This completes the proof of (2) which
implies (3).

If $\rho' \subset G'$ is a path in $G$ with endpoints $P_1,P_2 \in
\Fix(f') = \Fix(f)$, then there is a unique path $\rho \subset G$ with
endpoints $P_1$ and $P_2$ such that $g_\#(\rho) = \rho'$.  Condition
(2) implies that $\rho'$ is fixed by $(f')_\#$ if and only if $\rho$
is fixed by $f_\#$. This proves (4).

To show that $f' : G' \to G'$ is a   \rtt\  it suffices by  (1), (6) and (7) to prove that if $H_i$ is an \eg\ stratum with $i > j$ and if $\sigma \subset G_{i-1}$ is a connecting path for $H_i$   then $(f')_\#(\sigma)$ is non-trivial.   If $\sigma$ is contained in a non-contractible component of $G_{i-1}$ then its endpoints are $f$-fixed, and hence $f'$-fixed, by Remark~\ref{rtt2}.   Non-triviality of  $(f')_\#(\sigma)$ therefore follows from the fact that   $f'$ is a homotopy equivalence.  If $\sigma$ is contained in a contractible component of $G_{i-1}$ then it is contained in a zero stratum in $H^z_i$ so $(f')_\#(\sigma) = g_\# f_\# (\sigma)$.  If  the connecting path $f_\# (\sigma)$ is contained in a non-contractible component of $G_{i-1}$ then  $g_\# f_\# (\sigma)$ is non-trivial by the same argument used in the previous case.  Otherwise, $g_\# f_\# (\sigma) =  f_\# (\sigma)$ and again we are done.  This completes the proof that $f' : G' \to G'$ is a   \rtt.

   Item (5) follows from  (4) and  (6).
   \endproof

\subsection{Existence Theorem}  \label{s:existence}

\begin{thm}  \label{comp sp exists}\label{cs exists}
Suppose that $\oone \in \Out(F_n)$ is forward rotationless and that
$\cal C$ is a nested sequence of $\oone$-invariant free factor
systems.  Then $\oone$ is represented by a  \ct\ $\fG$ and filtration $\F$ that realizes $\cal C$.
\end{thm}

\proof  We assume without loss
$\cal C$ is maximal with respect to $\sqsubset$.  Thus any filtration
that realizes $\cal C$ is reduced.     By Theorem~\ref{rtt existence} and  item (3) of Lemma~\ref{formerly remarks}   we may choose a \rtt\ $\fG$ that represent $\oone$ and    realizes $\cal C$ and such that each contractible component of a filtration element is a union of zero strata and the endpoints of all \iNp s of \eg\ height are vertices.    For the remainder of the proof all \rtt s are assumed to   satisfy these properties.

\vspace{.1in}

\noindent{\bf Step 1: (\eg\ Nielsen Paths)}         Let $N(f)$ be the number of \iNp s of  \eg\ height.   In the construction of  an IRT    in \cite{bfh:tits1}  it is assumed (see Definition~5.2.1 of  \cite{bfh:tits1})    that  $N(f)$ is as small as possible.       The \eg\ properties of an IRT are then established by contradiction:  the failure of these properties allows one to reduce $N(f)$ which is impossible.     In order to make our constructions more algorithmic, we drop the assumption that  $N(f)$ is minimal and argue inductively:  the failure of  (\eg\ Nielsen Paths) allows one to reduce $N(f)$ and since $N(f)$ is finite,  this process eventually terminates in an $\fG$ satisfying (\eg\ Nielsen Paths).           As we are no longer    assuming that $N(f)$ is minimal   we cannot  quote statements of results from \cite{bfh:tits1} but must instead refer to  their proofs.

\begin{lemma}  \label{partial} Suppose that $H_r$ is an \eg\ stratum of a \rtt\ $\fG$ and that $\rho$ is an \iNp\ of height $r$.  If the fold at the illegal turn of $\rho$ is partial then there is a \rtt\ $f':G' \to G'$ satisfying $N(f') < N(f)$.  
\end{lemma}

\proof   This follows from the proofs of  Lemmas~5.2.3 and 5.2.4 of  \cite{bfh:tits1}.   The latter constructs    a topological representative of $\phi$ with all the desired properties and the former shows that this topological representative can be made into a \rtt\ without losing any of the desired properties.  \endproof

\begin{lemma}   \label{proper}Suppose that $H_r$ is an \eg\ stratum of a \rtt\ $\fG$,   that $\rho$ is an \iNp\ of height $r$ and that  the fold at the illegal turn of $\rho$ is proper.  Let  $f':G' \to G'$ be the \rtt\  obtained from $\fG$ by folding $\rho$.   Then $N(f') = N(f)$ and there  is  a bijection $H_s \to H'_s$ between the \eg\ strata of $f$ and the \eg\ strata of $f'$ such that $H'_s$  and $H_s$ have the same number of edges for all $s$. 
\end{lemma}   
\proof  This follows from the definition of  $f':G' \to G'$ and the  proof of Lemma~5.3.3 of  \cite{bfh:tits1}.
 \endproof

\begin{lemma}   \label{improper}Suppose that $H_r$ is an \eg\ stratum of a \rtt\ $\fG$ and that $\rho$ is an \iNp\ of height $r$.  If the fold at the illegal turn of $\rho$ is improper then there  is an \rtt\ $f':G' \to G'$ and a bijection $H_s \to H'_s$ between the \eg\ strata of $f$ and the \eg\ strata of $f'$ with the following properties.
\begin{enumerate}
\item    $N(f') = N(f)$. 
 \item $H'_r$ has  fewer edges than $H_r$.
 \item If $s > r$ then $H'_s$  and $H_s$ have the same number of edges.
\end{enumerate}
\end{lemma}

\proof  This follows from Definition~5.3.4 and the  proof of Lemma~5.3.5 of  \cite{bfh:tits1}.
 \endproof

 \begin{lemma}  \label{all proper} If $H_r$ is an \eg\ stratum of $\fG$ and $\rho$ is an \iNp\ of height $r$ such that the illegal turn at each \iNp\ obtained by iteratively folding $\rho$ is  proper then    $H_r$ satisfies (\eg\ Nielsen Paths).
      \end{lemma}
  
  \proof       The conclusion of Lemma~5.3.6 of \cite{bfh:tits1}  is that  $H_r$ satisfies (\eg\ Nielsen Paths).   The proof of that lemma uses  only standard folding arguments, the hypotheses of our lemma and  uniqueness of  the illegal turn of height $r$, which  follows from  Lemma~\ref{still works bh}.   
    \endproof
  
 \begin{cor}  \label{condition for noneg NP} Suppose that $H_r$ is an \eg\ stratum of a \rtt\ $\fG$ and that $\rho$ is an \iNp\ of height $r$.   Then the illegal turn at each \iNp\ obtained by iteratively folding $\rho$ is  proper if and only if  $H_r$ satisfies (\eg\ Nielsen Paths).
\end{cor}

\proof  This is an immediate corollary of Lemmas~\ref{stable} and \ref{all proper}.
\endproof
  
        Our algorithm for modifying a given $\fG$ so that it satisfies  (\eg\ Nielsen Paths) is as follows.  If some \eg\ stratum  does not  satisfy   (\eg\ Nielsen Paths),  let  $H_r$  be the highest such stratum.   By Lemma~\ref{all proper},  there is a (possibly empty) sequence of proper folds leading to a \rtt\ and \iNp\ with either a partial fold or an improper fold.  Apply Lemma~\ref{partial} or Lemma~\ref{improper} respectively.    If the resulting \rtt\  does not  satisfy   (\eg\ Nielsen Paths) go back to the beginning and start again.      
     
     \begin{remark} Iteratively folding any $\rho$ in $H_r$ either determines $f_r$ as in  (\eg\ Nielsen Paths) or leads to a partial or improper fold in a predictable number of steps.
     \end{remark}
 
     Suppose that the algorithm does not terminate.  Denote the \rtt s that are produced by $f=f_0, f_1, f_2 \ldots$.  Since $N(f_i)$ is non-decreasing and is strictly decreasing when a partial fold occurs, there are only finitely many such occurrences and  we may assume without loss that all the folds are full.    We make use of the bijection $H_s(i) \to H_s(j)$ between \eg\ stratum for $f_i$ and \eg\ stratum of $f_j$ given by Lemmas~\ref{proper} and \ref{improper}.   Let $H_r$ be the highest stratum for which (\eg\ Nielsen Paths) is not satisfied by $f_k$ for all sufficiently large $k$.   Then the   number of edges of height $r$ is a  non-increasing function of $k$ that strictly   decreases when an    improper fold of height $r$ occurs.  These folds do not therefore occur for sufficiently large $k$.  But this contradicts Lemma~\ref{all proper} and the choice of $r$. This proves that the algorithm terminates at a \rtt\ (still called) $\fG$ satisfying  (\eg\ Nielsen Paths). 
     
     In the steps that follow the number of edges in each \eg\ strata and the number of \iNp s of \eg\ height are not increased.  If after some modification, (\eg\ Nielsen Paths)  fails then we can return to step 1 and start again.  By the above argument this  terminates after finitely many repetitions.  (In fact, it is never necessary to return to step 1 but this requires an additional argument.)
\vspace{.1in}

\noindent{\bf Step 2: (Theorem~\ref{rtt existence})}       Apply steps (1) through (6) of the proof of Theorem~\ref{rtt existence}    to  produce a new $\fG$ satisfying the conclusions of that Theorem.     As noted above,  we may assume by  Remark~\ref{same Nielsen paths}  that  (\eg\ Nielsen Paths) is still satisfied.
 \vspace{.1in}

\noindent{\bf Step 3: ((Rotationless), (Filtration) and  (Zero Strata))}   Items (Rotationless) and   (Filtration)  follow from Proposition~\ref{rotationless} and Theorem~\ref{rtt existence}(F).    
   To achieve   (Zero Strata) it suffices, by  item (Z) of  Theorem~\ref{rtt existence}, to arrange that every edge in a zero stratum $H_i$ is $r$-taken.    Each edge $E$ in $H_i$ is contained in an $r$-taken path $\sigma \subset H_i$.  If $E$ is not $r$-taken, replace $E$ by a path that has the same endpoints as $\sigma$ and is marked by $\sigma$.  After finitely many such tree replacements, (Zero Strata)  is satisfied.  
    
 \vspace{.1in}
      
\noindent{\bf Step 4: (Periodic Edges)}  Suppose at first that no component $C$ of  of $\Per(f)$ is topologically a circle with each
point in $C$ having exactly two periodic directions.    Then   the endpoints of any periodic edge are principal, each periodic edge is fixed and each periodic stratum $H_r$ has a single edge $E_r$.  If  $G_{r-1}$ is not a core graph that contains both endpoints of $E_r$ then one could collapse $E_r$ without changing the free factor systems realized by the filtration elements, in violation of item (P) of Theorem~\ref{rtt existence}).  Thus (Periodic Edges) is satisfied.

     For the general case, it suffices to assume   that some component $C$ of $\Per(f)$ is topologically a circle with each
point in $C$ having exactly two periodic directions  and     modify $\fG$ to reduce the number of such
components.     

 Lemma~\ref{no iterates necessary}-(1) implies that $C$ is $f$-invariant and that $g= f|C$ is
orientation preserving.    By (Zero Strata) and the fact that   there are  no periodic directions based in $C$
and pointing out of $C$, every edge $E_j$ not in $C$ that has an endpoint in $C$ is
non-periodic, NEG and intersects $C$ in exactly its terminal endpoint.  Since all non-periodic vertices   are contained in \eg\ strata, no vertex in the complement of  $C$ maps into $C$.     Also, $C$ is a component of some $G_l$ by item (NEG) of  Theorem~\ref{rtt existence}. 
Let $E_m$ be the first   non-periodic \noneg\ edge $E_j$ that has terminal endpoint in $C$  and note that $f(E_m) = E_mC^d$.   We modify $\fG$ near $C$ in two steps as follows.

In the first step we make $C \subset \Fix(f)$.  Extend the rotation
$g^{-1}:C \to C$ to a map $h:G \to G$ that has
support on a small neighborhood of $C$, that is homotopic to the
identity and such that $h(E_j) \subset E_j \cup C$ for each
non-periodic \noneg\ edge $E_j$ that has terminal endpoint in $C$.
Redefine $f$ on each edge $E$ to be $h_\#f_\#(E)$.  The filtration is
unchanged.  Edges in $C$ are now fixed.  If $f(E_j) = E_ju_j$ then the
new $u_j$ and the old $u_j$ agree with the possible exception of
initial and terminal segments in $C$.   The $f$-image of all other edges is unchanged.    In fact, $f_\#(\sigma)$ is unchanged for any path
$\sigma$ with the property that endpoints of $f(\sigma)$ are not in
the support of $h$.  It     is straightforward to check that  $\fG$ is a \rtt\  and that all of the properties that we have established to date are preserved  with the
possible exception of item (P) of Theorem~\ref{rtt existence}, which   fails if one or more of the $E_j$'s is now a fixed edge that should be collapsed.  If there is no such edge then proceed to the next paragraph.   If there is such an edge, collapse it as in step 3 of the proof of  Theorem~\ref{rtt existence}.   That step is described very explicitly and we leave it to the reader, here and later in the proof, to check  that this operation does not undo previously established properties.   After finitely many such collapses,  we have  $C \subset \Fix(f)$ and  all previously established properties are preserved.  If $C$ now has outward pointing periodic directions we have finished our modifications of $C$.  Otherwise proceed to the next paragraph.
    
Recall that if $E_m$ is the first   non-periodic \noneg\ edge $E_j$ that has terminal endpoint in $C$, then$f(E_m) =E_m C^d$ for some $d \in \Z$.  In this
second step we arrange that $d=0$.  Choose $h' : G \to G$ that is the
identity on $C$, that satisfies $h'(E_j) = E_jC^{-d}$ for all $E_j$ and that has
support in a small neighborhood of $C$.  This map is homotopic to the
identity since we can simply unwind the twisting on $C$.  Redefine $f$
on each edge $E$ to be $h_\#f_\#(E)$ and note that $C \cup E_m \subset
\Fix(f)$ so the component of $\Fix(f)$ containing $C$ is no longer a topological circle.  The filtration is unchanged.  If necessary,   collapse fixed edges with an endpoint in $C$ and repeat this second step.   

 \vspace{.1in}
      
 \noindent{\bf Step 5: (Induction: the \noneg\ case)}  It remains to establish (Completely Split) and the items related to non-fixed \noneg\ edges.   We do this by induction up the filtration making use of sliding and the new move described in Section~\ref{s:new move}.  
 
  Let $NI$ be the
number of irreducible strata in the filtration and for each $0 \le m
\le NI$, let $G_{i(m)}$ be the smallest filtration element containing
the first $m$ irreducible strata.  We will prove by induction on $m$ 
  that for
all $0 \le m \le NI$, one can modify $f$ to arrange that $f|G_{i(m)}$ is a \ct.  The $m=0$ case is vacuously true so we assume that $f|G_{r}$ is a \ct\  for $r= i(m)$ and make  modifications to
arrange that $f|G_{s}$ is a \ct\   for $s=i(m+1)$.  In this step we assume that  $H_s$ is
\noneg\ and is hence a single edge $E_s$ satisfying $f(E_s) = E_s u_s$
for some path $u_s \subset G_{s-1}$. 

 By (Zero Strata), $r = s-1$.
The sliding operation described in section~\ref{s:modification}
(complete details in section~5.4 of \cite{bfh:tits1}) allows us to
modify $E_s$ and $u_s$ by choosing a path $\tau \subset G_{s-1}$ with
initial endpoint equal to the terminal vertex of $E_s$ and \lq
sliding\rq\ the terminal end of $E_s$ to the terminal vertex of
$\tau$.    As   noted in step 2, we may assume that sliding preserves (\eg\ Nielsen Paths).   

As a first case suppose that after sliding along $\tau$ we have $E_s \subset
\Fix(f)$.  For future reference note that by Lemma~\ref{sliding
prelim} this is equivalent to $[\bar \tau u_s f_\#(\tau)]$ being
trivial and hence equivalent to $ f_\#(E_s \tau) = E_s[u_s f_\#(\tau)]
= E_s\tau$; i.e., to $E_s \tau$ being a Nielsen path.

If both  endpoints of $E_s$ are contained in $G_{s-1}$ then
(Periodic Edges) is satisfied as are   all of   the conclusions of Theorem~\ref{rtt
existence} and the three properties established in step 3.   The remaining
properties of a \ct\   follow from the inductive
hypothesis.
   
If either endpoint of $E_s$ is not contained in $G_{s-1}$  then  collapse $E_s$ to a point as in
step 4.   None of the previously achieved properties are lost   and the remaining
properties of a \ct\   follow from the inductive
hypothesis.
   This completes the inductive step in the case that $E_s
\subset \Fix(f)$ is trivial after sliding.

We assume now that there is no choice of $\tau$ such that $E_s \tau$
is a Nielsen path. The following proposition is a combination of
Propositions~5.4.3 and 5.5.1 of \cite{bfh:tits1}.

\begin{proposition}   \label{good slide}
Suppose that
\begin{itemize}
\item [($i$)] $\fG$  is a \rtt\ that satisfies (\eg\ Nielsen Paths) 
\item  [($ii$)] $f|G_{s-1}$ satisfies (Completely Split)
\item   [($iii$)] $H_s$ is an \noneg\ stratum with single edge $E_s$ for which there does not exists a path
$\mu \subset G_{s-1}$ such that $E_s \mu$ is a Nielsen path. 
\end{itemize}
 Then
there exists a path $\tau \subset G_{s-1}$ with initial endpoint equal
to the terminal endpoint of $E_s$ such after performing the slide
associated to $\tau$ the following conditions are satisfied.
\begin{enumerate}
\item $f(E_s) = E_s \cdot u_s$ is a non-trivial splitting.
\item If $\sigma$ is a circuit or path with endpoints at vertices and
if $\sigma$ has height $s$ then there exists $k \ge 0$ such that
$f^k_\#(\sigma)$ splits into subpaths of the following type.   
\begin{enumerate} 
\item $E_s$ or $\bar E_s$
\item an exceptional path of height $s$  
\item a subpath of $G_{s-1}$
\end{enumerate}
\item $u_s$ is completely split and its initial vertex is principal.
\item    $f|G_{s}$ satisfies (Linear Edges).
\end{enumerate}
 \end{proposition}

\proof     The construction of a  path \ $\tau$ along which to slide is carried out in the proof of Proposition 5.4.3 of \cite{bfh:tits1}.  We assume that $\tau$ has been chosen to satisfy the conclusions of that proposition.  In particular,  $f(E_s) = E_s \cdot u_s$ is a splitting that is non-trivial by   ($iii$).   Thus (1) is satisfied. (The statement of Proposition 5.4.3 of \cite{bfh:tits1} allows the
possibility that $G$ is subdivided at a periodic point and that the
terminal endpoint of $E_s$ is one of the new periodic vertices.  By
the end of the construction, we will have shown that the terminal
endpoint of $E_s$ is principal and hence fixed.  At that point we can
undo the subdivision.)

For (3) we must make use of facts that are explicitly stated and used
in the proof of Proposition 5.4.3 of \cite{bfh:tits1} but are not
contained in its statement.  The first is that by a further slide one
can replace $u_s$ with $f^k_\#(u_s)$ for any $k \ge 1$.  Since
$f|G_{s-1}$ satisfies (Completely Split) we may assume by Lemma~\ref{iterates split} that $u_s$ is
completely split.  The second is that if $u_s = \alpha\cdot \beta$ is
a coarsening of the complete splitting of $u_s$, then by a further
slide we may assume that the terminal endpoint of (the new) $E_s$ is
the terminal endpoint of $\alpha$.  Thus to complete the proof of (3)
we need only show that the endpoint of some term in the complete
splitting of $u_s$ is principal.  The only way that this could fail
would be if $u_s$ has height $r'$ where $H_{r'}$ is \eg\ and if each
height $r'$ term in the  complete splitting of $u_s$ is a single
edge. After replacing $u_s$ with a sufficiently high iterate, we may
assume that $u_s$ has such a long ${r'}$-legal segment that every edge
in $H_{r'}$ occurs as a term in the complete splitting of $u_s$.
Lemma~\ref{some essential} then completes the proof of (3).

    If $E_s$ is a linear edge,  choose a root-free Nielsen path $w_s$ and $d_s \ne 0$ so that $u_s =w_s^{d_s}$.   If  $E_t \subset G_{s-1}$ is a linear edge with the same axis as $E_s$   then after reversing the orientation on $w_s$   we may assume that $w_t$ and $w_s$ agree as
oriented loops.  After a further slide as in the proof of (3) we may
assume that $w_s = w_t$.  Item ($iii$) implies that  $E_s \bar E_t$ is not a Nielsen path and hence that
$d_s \ne d_t$.  This completes the proof of item (4). 

Lemma~4.1.4 of \cite{bfh:tits1} states that if $\sigma$ is a height $s$ circuit or path with endpoints at vertices then $\sigma$  splits into
subpaths that are either contained in $G_{s-1}$ or are  basic
paths of height $s$  meaning that they, or their inverse has   the form $E_s \gamma$
  or $E_s \gamma \bar E_s$ for some $\gamma \subset
G_{s-1}$.  By Corollary~\ref{eg irt},  $\fG$   satisfies the hypotheses, and hence the conclusions, of
Proposition~5.5.1 of \cite{bfh:tits1}.  These conclusions address both types of basic paths of height $s$ and verify (2).
 \endproof

We assume now that we have performed the slide move of
Proposition~\ref{good slide}.   Since $u_s$ is non-trivial, $f|G_{s}$ satisfies (Periodic Edges) and all of the properties achieved in the first four steps of our construction.  Items (Completely Split), (Vertices), (\noneg\
Strata), (Linear Edges) and (Nielsen Paths) for $f|G_s$ follows from
Proposition~\ref{good slide} and these properties for $f|G_{s-1}$.
This completes the proof of the inductive step in the case that $H_s$
is \noneg.
\vspace{.1in}

 \vspace{.1in}
      
 \noindent{\bf Step 6: (Induction: the \eg\ case)}  Suppose  now that $H_s$ is \eg.
 items (Vertices),  (\noneg\ Strata), (Linear Edges) and (Nielsen Paths) for
 $f|G_s$ follows from these properties for $f|G_{s-1}$.
 
For each edge $E \subset H_{s}$, there is a decomposition $f(E) =
\mu_1 \cdot \nu_1 \cdot \mu_2 \cdot\ldots\cdot \nu_{m-1}\cdot \mu_m$
where the $\nu_l$'s are the maximal subpaths in $G_{r}$. Let
$\{\nu_l\}$ be the collection of all such paths that occur as $E$
varies over the edges of $H_{s}$.  By (RTT-$ii$), $f^k_\#(\nu_l)$ is
non-trivial for each $k$ and $l$.  By  Lemma~\ref{iterates split} we may choose $k$ so large that each
$f^k_\#(\nu_l)$ is \cs.  We may also assume that the endpoints of
$f^k_\#(\nu_l)$ are periodic and hence principal.  There are finitely
many connecting paths $\sigma$ contained in the strata (if any)
between $G_r$ and $H_s$.  Each $f(\sigma)$ is either a connecting path
or a non-trivial path in $G_r$ with fixed endpoints.  We may therefore
assume that $f^k_\#(\sigma)$ is \cs\ for each such $\sigma$.  After $k$
applications of  Lemma~\ref{new move} with $j = r$ (see in
particular item (7) of that lemma) we have that $f|G_{s}$ is \cs.
 This completes the induction step and so also the
proof of the theorem.  \endproof

\subsection{Further properties of a \ct} \label{cs properties} \label{s:extra}

The following lemma is an extension of Lemma~\ref{new iterates to}.

\begin{lemma} \label{iterates to}
Assume that $\fG$ is a \csirt. The following properties hold for every
principal lift $\ti f : \Gamma \to \Gamma$.
\begin{enumerate}
\item If $\ti v \in \Fix(\ti f)$ and a non-fixed edge $\ti E$
determines a fixed direction at $\ti v$, then $\ti E\ \subset \ti
f_\#(\ti E) \subset \ti f_\#^2(\ti E) \subset \ldots$ is an increasing
sequence of paths whose union is a ray $\ti R$ that converges to some
$P \in \Fix_N(\hat f)$ and whose interior is fixed point free.
\item For every isolated $P \in \Fix_N(\hat f)$ there exists $\ti E$
and $\ti R$ as in (1) that converges to $P$.  The edge $E$ is
non-linear.
\end{enumerate}
\end{lemma}

\proof For $\ti E$ as in (1) and for each $m > 0$, Lemma~\ref{stays cs} implies that 
 $\ti E \subset \ti f_\#(\ti E) \subset \ti f_\#^2(\ti E)\subset
\ldots$ is a nested sequence of completely split paths.      This
increasing sequence of paths defines a ray $\ti R'$ that converges to
some non-repelling fixed point $P \in \Fix_N(\hat f)$ and that,  by Corollary~\ref{induced splitting}, intersects $\Fix(\ti f)$ only in its initial endpoint.  This completes the proof of (1).

If $P \in \Fix_N(\hat f)$ is isolated then $\ti f$ moves points that
are sufficiently close to $P$ toward $P$ by Lemma~\ref{l: second from bk3}.  We may therefore choose a ray $\ti R$ that converges to
$P$ and that intersects $\Fix(\ti f)$ only in its initial endpoint.
Moreover, the initial edge $\ti E$ of $\ti R$ determines a fixed
direction by Lemma~\ref{l:maps over} and so extends to a fixed point
free ray $\ti R$ converging to some $Q \in \Fix_N(\hat f)$ by
(1).   Lemma~\ref{l:maps over} implies that $P =
Q$.  Since $P$ is isolated it is not an endpoint of the axis of a
covering translation and $E$ is not a linear edge.  \endproof

\begin{notn} If $\ti E$ and $P$ are as in Lemma~\ref{iterates to}(1) then we 
say that {\em $\ti E$ iterates to} $P$ and that $P$ is {\em associated to $\ti 
E$}.
\end{notn}

The following lemma is used in section~\ref{recog} and also in
\cite{fh:abeliansubgroups}.

\begin{lemma} \label{isolated for egs}  Suppose that $\otwo \in \ofn$ is forward
rotationless and that $P \in \Fix_N(\hat \atwo)$ for some $\atwo \in
\PA(\otwo)$.  Suppose further that $\Lambda$ is an attracting
lamination for some element of $\ofn$, that $\Lambda$ is
$\otwo$-invariant and that $\Lambda$ is contained in the accumulation
set of $P$.  Then $\PF_{\Lambda}(\otwo) \ge 0$ and
$\PF_{\Lambda}(\otwo) > 0$ if and only if $P$ is isolated in
$\Fix_N(\hat \atwo)$.
\end{lemma}

\proof Let $g : G \to G$ be a \csirt\ representing $\otwo$ , let $\ti
g : \Gamma \to \Gamma$ be the lift corresponding to $\atwo$ and let
$\ti R \subset \Gamma$ be a ray converging to $P$.  Choose a generic
leaf $\gamma \subset G$ of the realization of $\Lambda$ in $G$.  Then
every finite subpath of $\gamma$ lifts to a subpath of every subray of
$\ti R$.  If $P$ is not isolated in $\Fix_N(\hat{\atwo})$ then
Lemma~\ref{l: second from bk3} implies that every finite subpath of
$\gamma$ extends to a Nielsen path for $g$. By  (\eg\ Strata) and
(Nielsen Paths),  this extension
can be done with a uniformly bounded number of edges.  It is an
immediate consequence of the definition of the expansion factor
(Definition 3.3.2 of \cite{bfh:tits1}) that $\PF_{\Lambda}(\otwo) =
0$.

Assume now that $P$ is isolated in $\Fix_N(\hat{\atwo})$ and let
$\Sigma$ be the set of finite paths $\sigma \subset G$ with endpoints
at vertices and with the property that every finite subpath of
$\gamma$ is contained in $g^m_\#(\sigma)$ for some $m>0$.
Lemma~\ref{iterates to} and the assumption that $\Lambda$ is contained
in the accumulation set of $P$ imply that $\Sigma$ contains a path
that is a single edge and in particular is non-empty.

Let $\sigma \in \Sigma$ be an element of minimal height, say $k$. Then
$\sigma$ decomposes as a concatenation of edges $\mu_i \subset H_k$
and subpaths $\nu_i \subset G_{k-1}$ and we let $K$ be the number of
elements in this decomposition.  Choose a nested sequence of subpaths
$\gamma_j$ of $\gamma$ whose union equals $\gamma$.  Since the
$\nu_i$'s are not in $\Sigma$, there exists $J > 0$ so that $\gamma_j$
is not contained in any $g_\#^m(\nu_i)$ for $j > J$.  Since $\gamma$
is generic, it is birecurrent.  Choose $j_0 > J$ and $j > j_0$ so that
$\gamma_j$ contains at least $3K$ disjoint copies of $\gamma_{j_0}$.
There exists $m > 0$ such that $\gamma_j \subset g^m_\#(\sigma)$.  It
follows that $\gamma_{j_0} \subset g^m_\#(\mu_i)$ for some $i$.  There
is a choice of $i$ that works for all choices of $j_0$ and this proves
that $\mu_i \in \Sigma$.

Let $E$ be the single edge in $\mu_i$.  We assume that $E$ is \noneg\
and argue to a contradiction.  There is a path $u \subset G_{k-1}$
such that $g_\#^m(E) = E\cdot u\cdot\ldots\cdot g_\#^{m-1}(u)$ for all
$m > 0$.  Lemma~3.1.16 of \cite{bfh:tits1} states that $\gamma$ is not
a circuit.  It follows that $u$ is not a Nielsen path and hence that
the length of $g_\#^l(u)$ goes to infinity with $l$.  The birecurrence of
$\gamma$ and the fact that $E \in \Sigma$ imply that for every
$\gamma_j$ there exits $p >0$ such that $\gamma_j \subset g_\#^{p-
1}(u)g_\#^{p}(u)=g_\#^{p-1}(ug_\#(u))$ in contradiction to the
assumption that no element of $\Sigma$ has height less than $k$.

We now know that $H_k$ is EG.  Since $E \in \Sigma$, $\gamma$ is a
leaf in the attracting lamination $\Lambda_k$ associated to $H_k$.
There is a splitting of $g(E)$ into subpaths in $H_k$ and subpaths in
$G_{k-1}$.  If $\gamma$ were contained in $G_{k-1}$ then one of the
subpaths in $G_{k-1}$ would be contained in $\Sigma$ in contradiction
to our choice of $k$.  Thus $\gamma$ is not entirely contained in
$G_{k-1}$ and Lemma~3.1.15 of \cite{bfh:tits1} implies that $\gamma$
is a generic leaf of $\Lambda_k$.  In other words, $\Lambda =
\Lambda_k$.  It follows that $\PF_{\Lambda}(\otwo) > 0$ which
completes the proof.  \endproof

Assume that $\oone$ is forward rotationless and that $\fG$ is a
\csirt\ representing $\oone$.  Following the notation of
\cite{bfh:tits3} we say that an unoriented conjugacy class $\mu$ of a
root-free element of $F_n$ is an {\em axis for $\oone$} if for some
(and hence any) representative $c \in F_n$ there exist distinct
$\aone_1, \aone_2 \in \PA(\oone)$ that fix $c$.  Equivalently
$\Fix_N(\hat \aone_1) \cap \Fix_N(\hat \aone_2)$ is the endpoint set
of the axis $ A_c$ for $T_c$.   It is a consequence of Lemma~\ref{axes and lifts} below
that an unoriented conjugacy class $\mu$ is an axis for $\phi$ if and only if it is an axis for
a linear edge in some (every)   \ct\ representing $\phi$.

\begin{remark}  
In the context of the mapping class group, a conjugacy class is an
axis if and only if it is represented by a reducing curve in the
minimal Thurston normal form.
\end{remark}

Lemma~\ref{no iterates necessary} implies that the oriented conjugacy
class of $c$ is $\oone$-invariant.  By Lemmas~4.1.4 and 4.2.6 of
\cite{bfh:tits1}, the circuit $\gamma$ representing $c$ splits into a
concatenation of  subpaths $\alpha_i$, each of which is either a  fixed edge or an \iNp.   (\noneg\ Nielsen Paths) and   Corollary~\ref{eg irt} imply that each turn $(\bar \alpha_i, \alpha_{i+1})$ is legal.   Item  1 of Lemma~\ref{cs is unique} therefore implies that this splitting is the complete splitting of $\gamma$ .  

There is
an induced complete splitting of $A_c$ into subpaths $\ti \alpha_i$ that
project to either fixed edges or \iNp s. The lift $\ti f_0 : \Gamma
\to \Gamma$ that fixes the endpoints of each $\ti \alpha_i$ is a
principal lift by Corollary~\ref{new essential is principal} and
commutes with $T_c$.  We say that $\ti f_0$ and the corresponding
$\aone_0 \in \PA(\oone)$ are the {\em base lift} and {\em base
principal automorphism} associated to $\mu$ and the choices of $T_c$
and $\fG$.  (If $\mu$ is not represented by a basis element then
$\aone_0$ is independent of the choice of $\fG$ but otherwise it is
not.  Indeed, if $\mu$ is not represented by a basis element, then the
smallest free factor $F$ that carries $\mu$ has rank greater than one.
Since $\mu$ is not an axis for $\oone|F$ there is a unique principal
automorphism $\aone|F$ that fixes $c$ and $\aone$ must be the
extension of $\aone|F$.)     Item  2 of Lemma~\ref{cs is unique} implies  that for each $\ti \alpha_i$ and for each $\ti x \in \ti \alpha_i$, the nearest point to $\ti f_0(\ti x)$ in $A_c$ is contained in $\ti \alpha_i$.  It follows that 
  $\Fix(T_c^j\ti f_0) = \emptyset$ for all $j \ne 0$ and hence that $\ti f_0$ is the only lift that commutes with $T_c$ and has fixed
points in $A_c$.

\begin{lemma}  \label{axes and lifts}   Suppose that $\phi$ is forward rotationless and that the unoriented conjugacy class $\mu$ is an axis for $\phi$.  
Assume notation as above.  There is a bijection between the set of principal
lifts [principal automorphisms] $\ti f_j \ne \ti f_0$ [respectively
$\aone_j \ne \aone_0 \in \PA(\oone)$] that commute with $T_c$ [fix
$c$] and the set of  linear edges $\{E_j\}$ with axis equal to  $\mu$.
Moreover, if $f(E_j) = E_j w_j^{d_j}$ then $\ti f_j = T_c^{d_j} \ti
f_0$ [$\aone_j = i_c^{d_j}\aone_0$].
\end{lemma}

\proof The $w_j$'s in question are equal by (Linear Edges) and we
label this path $w$.  There is a lift $\ti w$ that is a fundamental
domain of $A_c$ and that is a Nielsen path for $\ti f_0$.  Let $\ti
E_j$ be the lift of $E_j$ that terminates at the initial endpoint $\ti
v$ of $\ti w$ and let $\ti f_j$ be the lift that fixes the initial
endpoint of $\ti E_j$.  Then $\ti f_j$ is a principal lift that
commutes with $T_c$ and satisfies $\ti f_j(\ti v) = T_c^{d_j}(\ti v) =
T_c^{d_j}(\ti f_0\ti v)$ which implies that $\ti f_j = T_c^{d_j} \ti
f_0$.

Conversely, if $\ti f \ne \ti f_0$ is a principal lift that commutes
with $T_c$ then $\ti f = T_c^d \ti f_0$ for some $d \ne 0$.  In
particular, $ A_c$ is disjoint from $\Fix(\ti f)$ and there is a ray
$\ti R_1$ that intersects $\Fix(\ti f)$ in exactly its initial
endpoint and that terminates at the endpoint $P$ of $A_c$ that is the
limit of the forward $\ti f$ orbit of $\ti v$.  Let $\ti E$ be the
initial edge of $\ti R_1$.  Lemma~\ref{l:maps over} implies that $\ti
E$ determines a fixed direction and also that $\ti R_1$ must be the
ray constructed from the initial edge $\ti E$ of $\ti R_1$ by
Lemma~\ref{iterates to}.  If $\ti f(\ti E) = \ti E \cdot \ti u$ then
$\ti R_1 = E\cdot \ti u \cdot \ti f_\#(\ti u) \cdot \ti f^2_\#(\ti
u)\cdot \ldots$.  Since $\ti R_1$ has a common infinite end with
$A_c$, it follows that $f^k_\#(u)$ is a periodic, hence fixed, Nielsen
path for sufficiently large $k$ and for $u$ equal to the projected
image of $\ti u$.  In particular, $u$ and $f_\#(u)$ have the same
$f_\#^k$-image, and since they have the same endpoints, they must be
equal.  In other words, $u$ is a Nielsen path.  This proves that $\ti
E$ is the lift of a linear edge $E$ whose associated axis is $\mu$.
By (Linear Edges),   $E=E_j$ and $d = d_j$ for some $j$ and, after
translating $\ti E$ by some iterate of $T_c$ if necessary, that $\ti
v$ is the terminal endpoint of $\ti E$.  \endproof

\begin{remark} \label{must be trivial}
Suppose that $\fG$ is a \csirt, that $C$ is a component of some
filtration element $G_s$, that $C$ has no valence one vertices and
that $\oone|\F(C)$ is the trivial outer automorphism.  Then $f|C$ is
the identity.  To see this, let $H_r$ be the first non-fixed stratum
in $C$.  It can not be EG because the identity element has no
attracting laminations.  If it were NEG it would have to be linear
because $f|G_{r-1}$ is the identity and it cannot be linear because
the identity element has no axes.
\end{remark}

 We conclude this section by showing that every element of  $\Out(F_n)$ has a uniformly bounded iterate that is forward rotationless.
 
\begin{lemma}  \label{uniform bound}
For all $n \ge 1$ there exists $K_n > 1$ so that $\oone^{K_n}$ is
forward rotationless for all $\oone \in \Out(F_n)$.
\end{lemma}

\proof   Given $\oone \in \Out(F_n)$, let $\fG$ be a \ct\ representing some  forward rotationless iterate $\psi =\oone^N$ of $\oone$.    By Corollary~\ref{non-trivial fixed set} and Lemma~\ref{equivalence classes}, the number of isogredience classes  of principal  lifts of $\psi$ is the same as the number of Nielsen classes for $ \fG$.  If $x$ is a  principal  vertex that has valence two and that is isolated in $\Fix(f)$ then  $x$ is  either 
 the initial endpoint of a non-fixed \noneg\ edge or  an endpoint of an \iNp\ of \eg\ height. 
By (Vertices) and Corollary~\ref{eg irt}, there   is a uniform (i.e. depending only on $n$)   upper bound to the number of   isolated fixed principal vertices.   By (Periodic Edges) there  is also a uniform upper bound  to the number of components of $\Fix(f)$ that contain at least one edge.    It follows that there is a uniform upper bound to  the number of   Nielsen classes for $ \fG$ and to the  number of edges based at principal vertices.  From the former we conclude that   the number of isogredience classes  of principal  lifts $\ti f : \Gamma \to \Gamma$ of $f$   is uniformly bounded.

      Since $\oone$ commutes with $\psi$,  it acts on the set of   isogredience classes  of principal automorphisms representing $\psi$.   After replacing $\phi$ with a uniformly bounded iterate, we may assume that $\oone$ fixes each isogredience class.   Thus, if $\Psi$ is a principal automorphism representing $\psi$ then there exists an automorphism $\Phi$ representing $\phi$ such that $\Phi$ commutes with $\Psi$.   In particular,   $\mathbb F := \Fix(\Psi)$ is $\Phi$-invariant.  By construction,  the outer automorphism determined by $\Phi|\mathbb F$ has finite order and so is represented by a homeomorphism of a graph with no valence one or valence two vertices.  Since the rank of $\mathbb F$ is uniformly bounded, the period of the outer automorphism determined by $\Phi|\mathbb F$ is uniformly bounded.  After replacing $\phi$ with a further uniformly bounded iterate, we may assume that  $\mathbb F \subset \Fix(\Phi)$.   Thus $\Fix (\hat \Phi)$ contains each non-isolated point of $\Fix(\hat \Psi)$ by Lemma~\ref{l: second from bk3}. 
      
        By   Lemma~\ref{iterates to}(2), the number of isolated points in $\Fix(\hat \Psi)$, up to the action of $\mathbb F$, is bounded above by the number of edges based at principal fixed points for $f$ and so is uniformly bounded.     We may therefore assume that if $P$ is an isolated point in $ \Fix(\hat \Psi)$ then $\hat \Phi(P) =  \hat T_a(P)$  for some  $a=a_P  \in \mathbb F$,  from which  it follows that    $\hat \Phi^N(P) =  \hat T_a^N(P)$.

    The proof now divides into cases.  If $\mathbb F$ is trivial then $\Fix(\hat \Psi) \subset \Fix(\hat \Phi)$.  If   $\mathbb F$ has rank at least two  then $\Phi^N = \Psi$.  It follows  that   $T_a$ is trivial and again $\Fix(\hat \Psi) \subset \Fix(\hat \Phi)$.    The final case is that $\mathbb F$ has rank one.       After replacing $\Phi$ with $T_a^{-1} \Phi$ we may assume that $P \in \Fix(\hat \Phi)$.  Since  $\Fix(\hat \Phi)$ and $ \Fix(\hat \Psi)$ have at least three points in common, $\Phi^N = \Psi$.  As in the higher rank case, it follows that $\Fix(\hat \Psi) \subset \Fix(\hat \Phi)$ in this case as well.    As this holds for each principal automorphism representing $\psi$, $\phi$ is forward rotationless. 
              \endproof

\section{Recognition Theorem} \label{recog}
In this section we specify invariants that uniquely determine a
forward rotationless $\oone$.  As a warm-up to the general theorem, we
consider the special case, essentially proved in \cite{bfh:tits0},
that $\oone$ is irreducible, meaning that there are no non-trivial
proper $\oone$-invariant free factor systems.  It follows that a
\csirt\ $\fG$ representing $\oone$ has a single stratum and that the
stratum is EG.  In particular, $\oone$ has infinite order and
$\L(\oone)$ has exactly one element.  Lemma~\ref{no iterates
necessary}(3) implies that all iterates of $\oone$ are irreducible.

\begin{lemma} \label{bk0}  If $\oone \in \Out(F_n)$ is irreducible and forward 
rotationless, then $\oone$ has infinite order and is determined by its
unique attracting lamination $\Lambda$ and the expansion factor
$\PF_{\Lambda}(\oone)$.  More precisely, if $\oone$ and $\otwo$ are
forward rotationless and irreducible and if they have the same unique
attracting lamination and the same expansion factor then $\oone =
\otwo$.
\end{lemma}

\proof As noted above, $\oone$ and $\otwo$ have infinite order and all
iterates of $\oone$ and $\otwo$ are irreducible.  Theorem 2.14 of
\cite{bfh:tits0} implies that $\otwo^{-1} \oone$ has finite order and
that $\otwo^k=\oone^k$ for some $k \ge 1$.  By Lemma~\ref{some
essential} and Lemma~\ref{lam2} there exists $\aone \in \PA(\oone)$
such that $\Fix_N(\aone)$ contains at least three points $P_1,P_2$ and
$P_3$, each of whose accumulation set equals $\Lambda$.  The
$\Fix_N$-preserving bijections between $\PA(\oone)$ and $\PA(\oone^k)$
and between $\PA(\otwo)$ and $\PA(\otwo^k)$ induce a
$\Fix_N$-preserving bijection between $\PA(\oone)$ and $\PA(\otwo)$.
Thus there exists $\atwo \in \PA(\otwo)$ such that $\Fix_N(\hat \atwo)
= \Fix_N(\hat \aone)$.

Choose a finite order homeomorphism $f:G\to G$ of a marked graph $G$
representing $\otwo^{-1} \oone$, let $\ti f : \Gamma \to \Gamma$ be
the lift corresponding to $\atwo^{-1}\aone$ and note that $P_1,P_2,P_3
\in \Fix(\hat f)$.  The line with endpoints $P_1$ and $P_2$ and the
line with endpoints $P_1$ and $P_3$ are $\ti f_\#$-invariant and since
$\ti f$ is a homeomorphism they are $\ti f$-invariant.  The
intersection of these lines is an $\ti f$- invariant, and hence $\ti
f$-fixed, ray $\ti R$ that terminates at $P_1$.  The lamination
$\Lambda$ is carried by the subgraph $G_0 \subset \Fix(f)$ that is the
image of $\ti R$.  Example~2.5(1) of \cite{bfh:tits0} implies that $G_0
= G$; thus $f$ is the identity and $\otwo = \oone$.  \endproof
 
If $\aone_1, \aone_2 \in \PA(\oone)$ and if there exists a non-trivial
indivisible $a \in \Fix(\aone_1) \cap \Fix(\aone_2)$, then
$\aone_2\aone_1^{- 1} = i_a^d$ for some $d \ne 0$.  We think of $d$ as
a {\em twist coefficient} for the ordered pair $(\aone_1, \aone_2)$
relative to $a$.  In our next example we show that an elementary
linear outer automorphism is determined by the fixed subgroups of its
principal automorphisms and by a twist coefficient.

\begin{ex} \label{linear Ex} Let $x_1,x_2,\dots,x_n$ be a basis for $F_n$ and 
let $F_{n-1}= \langle x_1,\dots, x_{n-1}\rangle$.  Define $\aone_1$ by
$\aone_1|F_{n-1} =$ identity and $\aone_1(x_n) = x_n a^d$ for some
non-trivial root-free $a \in F_{n-1}$ and some $d > 0$.  Define
$$
\aone_2 = i_{x_n}^{-1}\aone_1 i_{x_n} = i_{\bar x_n \aone_1(x_n)}\aone_1 = 
i_{a}^d \aone_1.
$$
 Then  $\Fix( \aone_2) = i_{\bar x_n}\Fix( \aone_1)$ and $\Fix(\aone_1) \cap 
\Fix(\aone_2) =    \langle a \rangle$.  Since $\Fix(\aone_1)$ and 
$\Fix(\aone_2)$ have rank greater than one,     $\aone_1, \aone_2 \in 
\PA(\oone)$.

We claim that if  there exist $\atwo_1,\atwo_2 \in \PA(\otwo)$ such that 
$\atwo_2 = i_{a}^d \atwo_1$ and such that $\Fix(\atwo_i) = \Fix(\aone_i)$ for $i 
=1,2$,   then $\otwo = \oone$.   
It is obvious that $\atwo_1|F_{n-1} = $ identity.  
Moreover, 
\begin{eqnarray*}
\Fix(i_{a^d}\atwo_1) & = & \Fix(i_{a^d}\aone_1) = \Fix(i_{ x_n}^{-1}
\aone_1 i_{x_n}) = i_{\bar x_n} \Fix(\aone_1) \\ & = & i_{\bar x_n}
\Fix(\atwo_1) = \Fix(i_{x_n}^{-1} \atwo_1 i_{x_n}) = \Fix(i_{\bar
x_n\atwo_1(x_n) }\atwo_1).
\end{eqnarray*}
Since $i_{a^d}\atwo_1$ and $i_{\bar x_n\atwo_1(x_n) }\atwo_1$
represent the same outer automorphism and have a common fixed subgroup
of rank greater than one, they are equal.  Thus $a^d = \bar
x_n\atwo_1(x_n) $ or equivalently $\atwo_1(x_n)= x_n a^d$. This proves
that $\atwo_1 = \aone_1$ and $\oone = \otwo$.
\end{ex}

We now turn to the general case.

\begin{thm}{\bf (Recognition Theorem)} \label{t:recognition} Suppose that 
$\oone,\otwo\in\Out(\f)$ are forward rotationless and that
\begin{enumerate}
\item  $\PF_{\Lambda}(\oone)=\PF_{\Lambda}(\otwo)$, for all $\Lambda\in 
\L(\oone)=\L(\otwo)$.
\item there is bijection $B : \PA(\oone) \to \PA(\otwo)$ such that:  
\begin{itemize} 
\item  [(i)]   {\bf (fixed sets preserved)} $\Fix_N(\hat \aone) = 
\Fix_N(\widehat{B(\aone)})$
\item [(ii)]    {\bf (twist coordinates preserved)} If $ w \in \Fix(\aone)$  and 
$\aone, i_w\aone  \in \PA(\oone)$, then $B(i_w\aone)=  i_wB(\aone)$.
\end{itemize}
\end{enumerate}
Then  $\oone=\otwo$.
\end{thm}

\begin{remark} \label{equivariance} The bijection $B$ is necessarily equivariant 
in the sense that $B(i_c \aone i_c^{-1}) = i_c B(\aone)i_c^{-1}$ for
all $c \in F_n$.  This follows from the fact that $\Fix(i_c \aone
i_c^{-1})= i_c(\Fix(\aone_1))$ and from Remark~\ref{distinct fixed
point sets}.  Thus $B$ is determined by its value on one
representative from each of the finitely many equivalence classes in
$\PA(\oone)$ and 2(i) can be verified by checking finitely many cases.
Similarly, 2(ii) can be verified by checking finitely many cases.  The
$w$'s to which 2(ii) apply have the form $w = a^d$ where $a$
represents a common axis of $\oone$ and $\otwo$.  The values of $d$
can be read off from \rtt s as in Lemma~\ref{axes and lifts}.
\end{remark}

\begin{remark}  The assumption in (1) that $ \L(\oone)=\L(\otwo)$ is redundant. It follows from 
Lemma~\ref{new iterates to} and 2(i). We include it in the statement
of the theorem for clarity.
\end{remark}

\proof The proof is by induction on $n$.  By convention, all forward
rotationless outer automorphisms are the identity when $n= 1$ so we
may assume that the theorem holds for all ranks less than $n$ and
prove it for $n$.

The case that both $\oone$ and $\otwo$ are irreducible is proved in
Lemma~\ref{bk0} so we may assume that at least one of these, say
$\oone$, is reducible and so admits a proper non-trivial invariant
free factor system. Since this free factor system is realized by a
filtration element in a \rtt\ representing $\oone$, some proper free
factor carries either an attracting lamination $\Lambda$ for $\oone$
or a $\oone$-periodic conjugacy class $[c]$.  Lemma~\ref{no iterates
necessary} implies that the elements of $\L(\oone) = \L(\otwo)$ are
invariant by both $\oone$ and $\otwo$, that $\oone$ and $\otwo$ have
the same periodic conjugacy classes and that all these conjugacy
classes are fixed.  The smallest free factor that carries
$\Lambda$ or $[c]$ is both $\oone$-invariant and $\otwo$-invariant by
Corollary~\ref{ffs invariance}.  This proves the existence of
non-trivial proper free factors that are that both $\oone$-invariant
and $\otwo$-invariant.

Among all proper free factor systems, each of whose elements is both
$\oone$-invariant and $\otwo$-invariant, choose one $\F
=\{[[F^1]],\dots,[[F^k]]\}$ that is maximal with respect to inclusion.
We claim that $\oone|F^i = \otwo|F^i$ for each $i$.  If $F^i$ has rank
one then this follows from Lemma~\ref{no iterates necessary}.  If
$F^i$ has rank at least two then principal automorphisms representing
$\oone|F^i$ and $\otwo|F^i$ extend uniquely to principal automorphisms
representing $\oone$ and $\otwo$.  Thus $\oone|F^i$ and $\otwo|F^i$,
which are forward rotationless by Corollary~\ref{restriction of
rotationless}, satisfy the hypothesis of Theorem~\ref{t:recognition}
and the inductive hypothesis implies that $\oone|F^i =
\otwo|F^i$. This verifies the claim.  Let $\fG$ be a \csirt\
representing $\oone$ with $[\pi_1(G_r)] = \F$ for some filtration
element $G_r$, which we may assume without loss has no valence one
vertices.  Then $f|G_r$ represents both $\oone|\F$ and $\otwo|\F$.
The proof now divides into two cases.  The arguments are sufficiently
elaborate that we treat the cases in separate subsections.

\subsection{The NEG Case} \label{pg recognition}
In this subsection we complete the proof of
Theorem~\ref{t:recognition} in the case that there exists $s > r$ such
that $G_s$ is not homotopy equivalent to $G_r$ and such that $H_i$ is
NEG for all $r < i \le s$.  After reordering the $H_i$'s if necessary,
we may assume that $G_s$ is obtained from $G_r$ as a topological space
by either adding a disjoint circle or by attaching an arc $E$ with
both endpoints in $G_r$. In the former case, $[\pi_1(G_s)]$ is both
$\oone$-invariant and $\otwo$-invariant in contradiction to the
assumption that $\F$ is maximal and the fact that $G_s$ is
disconnected.  Thus $G_s =G_r\cup E$ where $f(E)= \bar u_1Eu_2$ for
some closed paths $u_1,u_2 \subset G_r$.  If $E$ is a single edge of
$G$ then $s=r+1$ and at least one of $u_1$ or $u_2$ is trivial.
Otherwise $E$ is made up of two edges and $s = r+2$.

Choose lifts $\ti E, \ti u_1, \ti u_2 \subset \Gamma$ and $\ti f :
\Gamma \to \Gamma$ such that $\ti f(\ti E) = \ti u_1^{-1} \ti E \ti
u_2$.  Denote the component of $G_r$ that contains $u_i$ by $C^i$ and
the copy of the universal cover of $C^i$ that contains $\ti u_i$ by
$\Gamma_r^i$.  If $u_1$ or $u_2$ is trivial then (Remark~\ref{attaching vertex}) at least one of the endpoints of $\ti E$ is a  principal vertex that is fixed by $\ti f$.    Otherwise $\ti E$ subdivides into two
\noneg\ edges  whose common initial vertex  is principal  and is fixed by $\ti f$.
  Corollary~\ref{new essential is principal} therefore
implies that $\ti f$ is a principal lift.  Lemma~\ref{new iterates to}
implies that there is a line $\ti \gamma \subset \Gamma$ that crosses
$\ti E$ and has endpoints in $\Fix_N(\hat f)$.  The smallest free
factor system that carries $[\pi_1(G_r)]$ and $\gamma$ is both
$\oone$-invariant and $\otwo$-invariant.  Since $\F$ is maximal, $G=
G_s$.

Corollary~3.2.2 of \cite{bfh:tits1} implies that $\otwo$ is
represented by $g : G \to G$ such that $g|G_r = f|G_r$ and such that
$g(E) = \bar w_1 E w_2$ for some closed paths $w_1,w_2 \subset
G_r$. It suffices to prove that $u_i = w_i$. The cases are symmetric
so we show that $u_1 = w_1$.

Suppose at first that $C^1$ has rank one and hence is a topological
circle that is contained in $\Fix(f)$.  By (Periodic Edges), the vertices in $C^1$ are principal.   Thus at least one of $E$ or $\bar E$ determines a  
direction based in $C^1$ that is fixed by $Df$.  If $u_1$ is non-trivial then it must be
that $\bar E$ determines a fixed direction based in $C^1$.  In this
case $C^1 \cup E$ is a component of $G$ and hence equal to $G$.  We
conclude that $n =2$ and that there is a basis $\{x_1,x_2\}$ for $F_2$
and $d \ne 0$ such that $x_1 \mapsto x_1$ and $x_2 \mapsto x_2x_2^d$
defines an automorphism representing $\oone$.  This is a special case
of Example~\ref{linear Ex} and so $u_1 = w_1$ in this case.  We may
therefore assume that $u_1$ is trivial.  The symmetric argument with
$g$ replacing $f$ reduces us to the case that $u_1$ and $w_1$ are both
trivial and so equal.  We may now assume that $C^1$ has rank at least
two.
 
The principal lift $\ti g : \Gamma \to \Gamma$ that corresponds to
$\ti f$ under the bijection $B$ satisfies $\Fix_N(\hat g) =
\Fix_N(\hat f)$.  Since $\hat g$ fixes the endpoints of $\ti \gamma$
and $\ti E$ is the only edge in $\ti \gamma$ that does not project
into $G_r$, it follows that $\ti g(\ti E) = \ti w_1 \ti E \ti w_2$.
Let $\ti v \in \Gamma_r^1$ be the initial endpoint of $\ti E$.  Then
$\ti u_1$ and $\ti w_1$ are the paths in $\Gamma_r^1$ connecting $\ti
v$ to $\ti f(\ti v)$ and $\ti v$ to $\ti g(\ti v)$ respectively.  It
therefore suffices to show that $\ti f|\Gamma_r^1 = \ti g|\Gamma_r^1$.

We know that $\ti f|\Gamma_r^1$ and $\ti g|\Gamma_r^1$ are both lifts
of $f|C^1$ and that $\hat f|\partial \Gamma_r^1$ and $\hat g|\partial
\Gamma_r^1$ have a common fixed point $P$.  If $P$ is not an endpoint
of the axis of some covering translation $T_c$ of $\Gamma_r^1$, then
there is at most one lift of $f|C^1$ that fixes $P$ and we are done.
Suppose then that $P \in \{T_c^{\pm}\}$.  By Remark~\ref{at least
one}, there exists a principal lift of $f|C^1$ that commutes with
$T_c$.  This lift extends over $\Gamma$ to principal lifts $\ti f'$
and $\ti g'$ of $f$ and $g$ respectively.  Since $\Fix_N(\hat
f|\partial \Gamma_{r-1}) \subset \Fix_N(\widehat{f'}) \cap
\Fix_N(\widehat{g'})$ contains at least three points, $\ti g' = B(\ti
f')$.  Condition 2(ii) therefore implies that $\ti g = T_c^d \ti g'$
and $\ti f = T_c^d \ti f'$ for some $d \in \Z$.  We conclude that $\ti
f|\Gamma_r^1 = \ti g|\Gamma_r^1$ as desired.

\subsection{The EG Case} \label{eg recognition} 
In this subsection we prove Theorem~\ref{t:recognition} assuming that
there exists $s > r$ where $H_s$ is exponentially growing and where
the union of the non-contractible components of $G_{s-1}$ is homotopy
equivalent to $G_r$.  In light of subsection~\ref{pg recognition},
this completes the proof of Theorem~\ref{t:recognition}.  Since
$G_{s-1}$ and $G_r$ carry the same elements of $\L(\oone)$, all
irreducible strata between $G_r$ and $G_s$ are NEG.  Our choice of
$\F$ guarantees that $G_s =G$.

Denote $ \otwo^{-1}\oone$ by $\theta$.  We must show that $\theta$ is
trivial.  By construction, $\theta|[[F^l]]$ is trivial for each
$[[F^l]] \in \F$ and the attracting lamination $\Lambda$ associated to
$H_s$ is $\theta$-invariant with expansion factor one.  Moreover, for
any principal lift $\aone$ of $\oone$ there is a unique lift $ \Theta$
of $\theta$ such that $\Fix(\hat \Theta) \supset \Fix_N(\hat \aone)$.

Each $F^l$ corresponds to a non-contractible component $D_l$ of
$G_{s-1}$.  Let $\ti D_l$ be the component of the full pre-image of
$D_l$ whose accumulation set in $\partial F_n$ is $\partial F^l$.
Suppose that $v \in H_s \cap D_l$ and that $\ti v \in \ti D_l$ is a
lift of $v$.  Then $\ti v$ is principal by Remark~\ref{attaching vertex} and the lift $\ti
f_{\ti v}$ of $f$ that fixes $\ti v$ is principal by
 Lemma~\ref{new essential is
principal}. The link of $\ti v$ contains edges that project to $H_s$
and determine fixed directions for $\ti f_{\ti v}$.  Lemma~\ref{new
iterates to} implies that any such edge extends to a fixed point free
ray that terminates at a point $P \in \Fix_N(\hat f_{\ti v})$ whose
accumulation set is $\Lambda$.  Let $\mathcal P_l$ be the union of
such $P$ for all $v \in H_s \cap D_l$ and all lifts $\ti v \in \ti
D_l$.

\begin{lemma} \label{fixes P}  There is a lift $\Theta$ of $\theta$ such that 
$\partial F^l \cup \mathcal P_l \subset \Fix(\hat \Theta)$.
\end{lemma}

\proof Assume at first that $F^l$ has rank at least two.  Let $
\Theta$ be the unique lift of $\theta$ such that $\partial F^l \subset
\Fix(\hat \Theta)$. If $P \in \mathcal P_l$ corresponds to $\ti v$ as
above and if $\aone_1 \in \PA(\oone)$ corresponds to $\ti f_{\ti v}$
then $\Fix_N(\hat \aone_1)$ contains $P$ and intersects $\partial F^l$
non-trivially.  There exists $\Theta_1$ such that $\Fix(\hat \Theta_1)
\supset \Fix_N(\hat \aone_1)$.  If $\Fix(\hat \aone_1) \cap \partial
F^l \ne \{T_c^{\pm}\}$ for some covering translation $T_c$, then
$\Theta_1 = \Theta$ and we are done.  Suppose then that $\Fix(\hat
\aone_1) \cap \partial F^l = T_c^{\pm}$.  By Remark~\ref{at least one}
there is a principal lift $\aone_2$ such that $\Fix(\aone_2)$ contains
$T_c^{\pm}$ and such that $\aone_2|F_l$ is a principal lift of
$\oone|F_l$.  In particular, $\aone_2 = i_c^d \aone_1$ for some $d \ne
0$.  By hypothesis, there are principal lifts $\atwo_1$ and $\atwo_2$
such that $\Fix_N(\hat \atwo_i) = \Fix_N(\hat \aone_i)$ and such that
$\atwo_2 = i_c^d \atwo_1$.  Thus $\Theta = \atwo_2^{-1} \aone_2=
\atwo_1^{-1} \aone_1 = \Theta_1$ where the first equality comes from
the fact that $\Fix(\hat \atwo_2^{-1} \hat \aone_2)$ contains at least
three points in $\partial F^l$.

It remains to consider the case that $F^l$ has rank one.  For each $P
\in \mathcal P_l$, there exist $\aone_P$ and $\atwo_P$ such that
$\Fix_N(\hat \aone_P) = \Fix_N(\hat \atwo_P)$ contains $P \cup
\partial F^l$.  Define $\Theta_P = \atwo_P^{-1} \aone_P$.  For any $Q
\in \mathcal P_l$ there exists $w$ such that $\aone_P = i_w\aone_Q $
and $\atwo_P = i_w\atwo_Q$.  It follows that $\Theta_P =\atwo_P^{-1}
\aone_P = \atwo_Q^{-1} \aone_Q = \Theta_Q$ as desired. \endproof

\begin{corollary} If $\theta$ has finite order then $\theta$ is trivial.
\end{corollary}

\proof If $\theta$ has finite order then  \cite{mc:finite}  there
is a marked graph $X$, a subgraph $X_0$ such that $\F(X_0) = \F$ and a
homeomorphism $h : X \to X$ that represents $\theta$ and is the
identity on $X_0$.  By Lemma~\ref{fixes P} there is an
$h_\#$-invariant $R$ whose initial endpoint is in $X_0$ and whose
accumulation set contains $\Lambda$.  No proper free factor system
carries $\F$ and $\Lambda$, so $R$ crosses every edge in $X \setminus
X_0$.  Since $h$ is a homeomorphism $R \subset \Fix(h)$ and we
conclude that $ h$ is the identity.  \endproof

We now assume that $\theta$ has infinite order and argue to a
contradiction. There is no loss in replacing $\theta$ by an iterate,
so we may assume that both $\theta$ and $\theta^{-1}$ are forward
rotationless.  There is a \csirt\ $h : G' \to G'$ representing
$\theta$ and there exists $r' < s'$ such that $G'_{s'} = G'$,
such that $\F(G'_{r'}) = \F$ and such that $h|G'_{r'} = $ identity
(see Remark~\ref{must be trivial}).

\begin{lemma} \label{not isolated} Suppose that $\mathcal P_l$ and $ \Theta$ are 
as in Lemma~\ref{fixes P} and that $P \in \mathcal P_l$.  Then $P$ is
not isolated in $\Fix(\hat \Theta)$.
\end{lemma}

\proof Suppose at first that $P$ is an attractor for $\hat \Theta$.
Let $\ti h$ be the lift of $h$ corresponding to $\Theta$. By
Lemma~\ref{iterates to} there is an edge $\ti E$ that iterates to $P$;
let $\ti R$ be the ray connecting $\ti E$ to $P$.  If $E$ belongs to
an EG stratum, then $\Lambda$, which is the accumulation set of $P$,
is an attracting lamination for $\theta$ by Lemma~\ref{isolated for
egs}.  This contradicts the fact that $\theta$ acts on $\Lambda$ with
expansion factor one.  If $E$ is NEG, then $\Lambda$ is carried by $G'
\setminus E$ in contradiction to the fact that no proper free factor
can carry $\F$ and $\Lambda$.  This proves that $P$ is not an
attractor for $\hat \Theta$.

The symmetric argument using a \rtt\ for $\theta^{-1}$ proves that $P$
is not a repeller so Lemma~\ref{l: second from bk3} completes the
proof.  \endproof

\begin{corollary}  \label{extension}
If $\gamma' \subset G'$ is a finite subpath of either:
\begin{enumerate}
\item   a  leaf of the realization of $\Lambda$ in $G'$ or 
\item the projection of the line in $\Gamma'$ connecting a pair of points 
$P_1,P_2 \in \mathcal P_l$
\end{enumerate}
then $\gamma'$  extends to a Nielsen path for $h$.
\end{corollary}

\proof Let $\Theta$ be as in Lemma~\ref{fixes P} and let $\ti h :
\Gamma' \to \Gamma'$ be the lift corresponding to $\Theta$. For
case (1), let $\ti R' \subset \Gamma'$ be a ray converging to $P \in
\mathcal P_l$.  There are lifts $\ti \gamma' \subset \ti R'$ of
$\gamma'$ that are arbitrarily close to $P$.  Lemma~\ref{not isolated}
and Lemma~\ref{l: second from bk3} therefore imply that $\ti \gamma'$
extends to a Nielsen path for $\ti h$.  In case (2), $P_1,P_2 \in
\Fix(\hat h)$.  Lemma~\ref{not isolated} and Lemma~\ref{l: second from
bk3} imply that any finite subpath of the line connecting $P_1$ to
$P_2$ extends to a Nielsen path for $\ti h$.  \endproof

It is well known that if $\oone$ acts trivially on conjugacy classes
in $F_n$ then $\oone$ is the trivial element.  This can be proved by
induction up the strata of $\fG$ representing $\oone$ or directly as
in Lemma~3.3 of \cite{fh:commensurator}.  The following lemma
therefore completes the proof of Theorem~\ref{t:recognition}.

\begin{lemma} \label{large overlap}
$\theta$ fixes each conjugacy class $[c]$ in $F_n$.
\end{lemma}

\proof If $v$ is a vertex in $G$ whose link $lk(v)$ is contained in
$H_s$, then the local stable Whitehead graph $SW_v$ is defined to be
the graph with one vertex for each oriented edge based at $v$ whose initial
direction is fixed by $Df$ and an edge connecting the vertices
corresponding to $E_1$ and $E_2$ if there is an edge $E$ of $H_s$ and
$k \ge 1$ so that the path $f^k_\#(E)$ contains $\bar E_1 E_2$ or
$\bar E_2 E_1$ as a subpath.  By Lemma~\ref{lam2} this is equivalent
to $\bar E_1 E_2$ or $\bar E_2 E_1$ being a subpath of a generic leaf
of $\Lambda$.  If $SW_v$ is not connected then then one can blow up
$v$ to an edge $E$ as in Proposition 5.4 of \cite{bh:tracks} to obtain
a proper free factor that carries $\F$ and $\Lambda$.  Since this is
impossible, $SW_v$ is connected.

Choose a positive integer $M$ such that $Df^M$ maps every direction in
$H_s$ to a fixed direction in $H_s$.  At one point in the proof we
need a way to choose partial edges and for this we subdivide the edges
of $H_s$ at the full $f^M$-pre-image of the set of vertices.  Edges in
this subdivision will be called {\em edgelets}.  Thus an edgelet maps
by $f^M$ to an edge.

 Let $g : G \to G'$ be a homotopy equivalence that respects the
marking and that satisfies $g(G_r) = G_{r'}'$.  We show below that
there is a positive integer $N$ so that for all circuits $\sigma
\subset G$ the conjugacy class in $F_n$ determined by
$g_\#f^{MN}_\#(\sigma) \subset G'$ is fixed by $\theta$.  Since every
conjugacy class in $F_n$ is realized in this manner by some $\sigma$,
this completes the proof of the lemma.

Since $\theta$ acts by the identity on $\F(G'_{r'})$ we may assume
without loss that $\sigma$ crosses at least one edge in $H_s$.  The
proof involves choosing a closed curve that is homotopic to
$f^{MN}_\#(\sigma)$ and a covering of that curve by subpaths with
large overlap.

To begin, choose a cyclic ordering of the $m$ edges of $H_s$ in
$\sigma$.  Define $\sigma_1$ to be first edge of $H_s$ in $\sigma$,
$\sigma_3$ to be second edge of $H_s$ in $\sigma$ and $\sigma_2$ to be
the subpath of $\sigma$ that begins with the last edgelet in
$\sigma_1$ and ends with the first edgelet in   $\sigma_3$.  Define
$\sigma_5$ to be third edge of $H_s$ in $\sigma$ and $\sigma_4$ to be
the subpath of $\sigma$ that begins with the last edgelet in
$\sigma_3$ and ends with the first edgelet in the $\sigma_5$.
Continue in this manner stopping with $\sigma_{2m}$ that begins with
the last edgelet in $\sigma_{2m-1}$ and ends with the first edgelet in
$\sigma_1$.

Let $\rho=f^{M}_\#(\sigma)$ and let $\rho_i = f^{M}_\#(\sigma_i)$.
Then each $\rho_i$ is an edge path in $G$ whose initial and terminal
edges are in $H_s$ and whose initial and terminal directions are fixed
by $Df$.  Suppose that $\rho_{2j} = \bar E E'$ where the link of the
common initial endpoint $v$ of $E$ and $E'$   is entirely contained
in $H_s$.  Since $SW_v$ is connected, there are edges $E=
E_1,E_2,\dots, E_l = E'$ in $lk(v)$ with initial directions fixed by
$Df$ such that each $\bar E_pE_{p+1}$ is a subpath of a generic leaf
of $\Lambda$.  Replace $\rho_{2j}$ by the concatenation $(\bar E_1
E_2) \cdot (\bar E_2 E_3) \cdot \dots \cdot (\bar E_{l-1}E_l)$.

After adjusting the indices, we have produced paths $\rho_1, \dots,
\rho_t$ with the following properties:
\begin{description}
\item [(a)] The initial edge $E_i^1$ of $\rho_i$ and the terminal edge
$E_i^2$ of $\rho_i$ are contained in $H_s$ and $E_i^2$ equals
$E_{i+1}^1$ up to a possible change of orientation.
\item [(b)] For all $k \ge 0$, $f_\#^k(\rho_i)$ is a finite subpath of either:
\begin{enumerate}
\item   a  generic leaf of the realization of $\Lambda$ in $G$ or 
\item the projection of the line in $\Gamma$ connecting a pair of points 
$P_1,P_2 \in \mathcal P_l$.
\end{enumerate}
\end{description}
 Suppose that each $E_i^1$ as been decomposed into proper subpaths
$E_i^1 =a_ib_i$.  The
equality $E^2_i=E^1_{i+1}$ or $\bar E^1_{i+1}$ determines a
corresponding decomposition of $E^2_i$  .Define $\tau_i$ from $\rho_i$ by deleting the
initial $a_i$ segment of $E_i^1$ and by deleting the terminal segment
of $E_i^2$ determined from (a) as follows. If $E_{i}^2 = E_{i+1}^1$
then remove the terminal $ b_{i+1}$ ; if $E_{i}^2 = \bar E_{i+1}^1$
then remove the terminal $ a_{i+1}$.
\begin{description}
\item [(c)] For any $\{a_i\}$ as above, $\rho$ is homotopic to the loop 
determined by  the concatenation of the $\tau_i$'s.  
\end{description}

Choose $K$ greater than the number of edges with height $s'$ in any
\iNp\ for $h$.  By \cite{co:bcc} there is a positive constant $C$ so
that if $\beta_1 \subset \beta_2$ are finite subpaths in $G$ then
$g_\#(\beta_2) \subset G'$ contains the subpath of $g_\#(\beta_1)$
obtained by removing the initial and terminal segments of edge length
$C$.  Since generic leaves of $\Lambda$ are birecurrent and since the
realization of $\Lambda$ in $G'$ can not be contained in $G'_{r'}$,
there is a subpath $\gamma'$ of a generic leaf of the realization of
$\Lambda$ in $G'$ that contains $2K+3C$ edges of $H'_{s'}$.  Choose a
subpath $\gamma$ of a generic leaf of the realization of $\Lambda$ in
$G$ whose $g_\#$ image contains $\gamma'$.  There exists $N > 0$ so
that $f^N_\#(E)$ contains $\gamma$ as a subpath for each edge $E$ of
$H_s$.  It follows that the path $g_\#(f^N_\#(E))$ contains at least
$2K+C$ edges of $H'_{s'}$ for each edge $E$ of $H_s$.

The subpath $\nu_i'$ of $g_\#(f^N_\#(\rho_i))$ defined by removing
initial and terminal segments with exactly $C$ edges of $H'_{s'}$ is
contained in either the realization of a leaf of $\Lambda$ in $G'$ or
the projection of a line in $\Gamma'$ connecting a pair of points
$P_1,P_2 \in \mathcal P_l$.  Corollary~\ref{extension} implies that
$\nu_i'$ extends to a Nielsen path $\mu_i'$ for $h$.  Let $\mu_i' =
\mu'_{i,1} \cdot \mu'_{i,2} \cdot \ldots \cdot \mu'_{i,m_i}$ be the
complete splitting of $\mu_i'$.  There is no loss in assuming that
$\mu'_{i,2} \subset \nu'_{i}$.  There are at most $K+C$ edges of
$H'_{s'}$ in $g_\#(f^N_\#(\rho_i))$ that precede $\mu_{i,2}'$.
Without changing this estimate we may assume that $\mu_{i,2}'$ has
height $s'$.  Note that $\mu'_{i,2} \subset g_\#(f^N_\#(E^1_i))$ and
hence that $\mu'_{i+1,2} \subset g_\#(f^N_\#(E^2_i))$.  Let $a_i$ be
an initial segment of $E_i^1$ such that $g_\#f^N_\#(a_i)$ is the
initial segment of $g_\#(f^N_\#(\rho_i))$ that precedes $\mu'_{i,2}$.
Define the $\tau_i$'s as in (c).

Lemma~\ref{cs is unique}(3) implies that there exists $l \le m_i$ such
that $\mu'_{i,l} = \mu'_{i+1,2}$ up to a change of orientation.  Thus
$g_\#f^N(\tau_i) = \mu'_{i,2} \cdot \dots \cdot \mu'_{i,l}$ is a
Nielsen path for $h$.  Property (c) implies that the conjugacy class
determined by $g_\#f^N_\#(\rho) = g_\#f^{MN}(\sigma)$ is
$\theta$-invariant as desired.  \endproof
  
\begin{corollary} \label{equal without iteration}
If $\oone$ and $\otwo$ are forward rotationless and if $\oone^m =
\otwo^m$ for some $m > 0$ then $\oone = \otwo$.
\end{corollary}

\proof Since $\oone$ and $\otwo$ are forward rotationless there are
$\Fix_N$-preserving bijections between $\PA(\oone^m)$ and $\PA(\oone)$
and between $\PA(\otwo^m)$ and $\PA(\otwo)$.  By assumption,
$\PA(\oone^m) = \PA(\otwo^m)$ so there is a $\Fix_N$-preserving
bijection between $\PA(\oone)$ and $\PA(\otwo)$.  The lemma now
follows from the Recognition theorem and the fact that expansion
factors and twist coefficients for $\oone^m$ are $m$ times those of
$\oone$ and similarly for $\otwo$.  \endproof

\bibliographystyle{plain}

\bibliography{ref}

\end{document}